\documentclass[11pt]{article}
\usepackage{latexsym,amssymb,amsmath}
\begin{document}

\baselineskip=20pt

\newcommand{\rd}{\mbox{Res}}
\newcommand{\kn}{\mbox{ker}}
\newcommand{\psp}{\vspace{0.4cm}}
\newcommand{\pse}{\vspace{0.2cm}}
\newcommand{\ptl}{\partial}
\newcommand{\dlt}{\delta}
\newcommand{\Dlt}{\Delta}
\newcommand{\sgm}{\sigma}
\newcommand{\cta}{\theta}
\newcommand{\al}{\alpha}
\newcommand{\be}{\beta}
\newcommand{\G}{\Gamma}
\newcommand{\g}{\gamma}
\newcommand{\lmd}{\lambda}
\newcommand{\td}{\tilde}
\newcommand{\vf}{\varphi}
\newcommand{\ad}{\mbox{ad}}
\newcommand{\stl}{\stackrel}
\newcommand{\ol}{\overline}
\newcommand{\es}{\epsilon}
\newcommand{\vsi}{\varsigma}
\newcommand{\ves}{\varepsilon}
\newcommand{\la}{\langle}
\newcommand{\ra}{\rangle}
\newcommand{\vt}{\vartheta}
\newcommand{\wt}{\mbox{wt}\:}
\newcommand{\sym}{\mbox{sym}}
\newcommand{\for}{\mbox{for}}
\newcommand{\rta}{\rightarrow}
\newcommand{\der}{\mbox{Der}}
\newcommand{\mbb}{\mathbb}
\newcommand{\dg}{\dag}
\newcommand{\tr}{\mbox{tr}\:}
\newcommand{\llra}{\Longleftrightarrow}
\newcommand{\Z}{{\cal Z}}
\newcommand{\vcm}{|0\ra}

\begin{center}{\Large \bf Representations of Centrally-Extended Lie Algebras }\end{center}
\begin{center}{\Large \bf over Differential Operators and Vertex Algebras}\footnote{1991 Mathematical Subject Classification. Primary 17B10,
17B69; Secondary 81Q40}\end{center} \vspace{0.2cm}

\begin{center}{\large Xiaoping Xu}\footnote{Research Supported
 by China NSF 10371121}\footnote{ACKNOWLEDGEMENT:
Part of this work was done during the author's visit to The
University of Sydney, under the financial support from Prof.
Ruibin Zhang's ARC research grant. The author thanks Prof. Zhang
for his invitation and hospitality.}
\end{center}
\begin{center}{Institute of mathematics, Academy of Mathematics \&
System Sciences}\end{center}
\begin{center}{Chinese Academy of Sciences, Beijing 100080, P. R.
China}\end{center}

\vspace{2cm}

 \begin{center}{\Large\bf Abstract}\end{center}
\vspace{1cm}

{\small We construct irreducible modules of centrally-extended
classical Lie algebras over left ideals of the algebra of
differential operators on the circle, through certain irreducible
modules of centrally-extended classical Lie algebras of infinite
matrices with finite number of nonzero entries. The structures of
vertex algebras associated with the vacuum representations of
these algebras are determined. Moreover, we prove that under
certain conditions, the highest weight irreducible modules of
centrally-extended classical Lie algebras of infinite matrices
with finite number of nonzero entries naturally give rise to the
irreducible modules of the simple quotients of these vertex
algebras. Our results are natural generalizations of the
well-known WZW models in conformal field theory associated with
affine Kac-Moody algebras. They can also be viewed as quadratic
generalizations of free field theory. In the very special case of
$W_{1+\infty}$, our results with integral central charges are more
direct and explicit than those of Kac and Radul.} \vspace{0.8cm}

\newpage

\section{Introduction}

The well-known $W_{1+\infty}$ Lie algebra is the
centrally-extended Lie algebra of differential operators on the
circle. It serves as the symmetry algebra of the famous
KP-hierarchy in integrable systems. In fact, it is the Lie algebra
of both rank-one charged quadratic free bosonic fields and
rank-one charged quadratic fermionic fields. Kac and Radul [KR1]
gave a classification of quasifinite highest weight irreducible
modules of the $W_{1+\infty}$ algebra. Frenkel, Kac, Radul and
Wang [FKRW] proved that the categories of irreducible modules of
the vertex operator algebras associated to ${\cal W}_{1+\infty}$
and ${\cal W}(gl_N)$ with integral central charge $N$ are
equivalent. Moreover, Kac and Radul [KR2] used free quadratic
bosonic fields to study the irreducible representations of the
vertex operator algebras associated to ${\cal W}_{1+\infty}$ with
negative integral central charge. The results in [KR1] were
extended to certain Lie subalgebras of ${\cal W}_{1+\infty}$ by
Kac,  Wang and Yan [KWY].

In general, the Lie algebras of  charged free quadratic bosonic
fields and fermionic fields of any rank are centrally-extended
general linear Lie algebras over differential operators on the
circle. They are the symmetry algebras of the multi-component
KP-hierarchies in integrable systems (cf. [V]). Boyallian and
Liberati [BL1] studied certain classical Lie subalgebras of the
Lie algebra of matrix differential operators on the circle. Ma
[M1, M2] systematically investigated the conformal algebra
structures and two-cocycles of centrally-extended classical Lie
superalgebras over left ideals of the algebra of differential
operators on the circle.

In this paper, we give constructions of irreducible modules of
centrally-extended classical Lie algebras over left ideals of the
algebra of differential operators on the circle, through certain
irreducible modules of centrally-extended classical Lie algebras
of infinite matrices with finite number of nonzero entries. The
structures of vertex algebras associated with the vacuum
representations of these algebras are determined. Moreover, we
prove that under certain conditions, the highest weight
irreducible modules of centrally-extended classical Lie algebras
of infinite matrices with finite number of nonzero entries
naturally give rise to the irreducible modules of the simple
quotients of these vertex operator algebras. It turns out that our
results resemble those for the vertex operator algebras associated
with affine Kac-Moody algebras, which are called WZW models in
conformal field theory. With special examples of the Hecke
algebras and group algebras borne in mind, we present general
constructions of representations of certain Lie subalgebras of the
centrally-extended  Lie algebra of the tensor algebra of any
associative algebra with the algebra of differential operators on
the circle. Our representation theory can be viewed as quadratic
generalizations of free field theory. Below we give more detailed
and technical introduction.

Throughout this paper, all the variables are formal and commute
with each other. All the vector spaces are assumed over $\mbb{C}$,
the field of complex numbers. Denote by $\mbb{Z}$ the ring of
integers and by $\mbb{N}$ the additive semigroup of nonnegative
integers.

Denote $\ptl_t=d/dt.$ Let
$$\mbb{A}=\sum_{i=0}^{\infty}\mbb{C}[t,t^{-1}]\ptl_t^i\eqno(1.1)$$
be the algebra of differential operators on the circle. Let
$M_{n\times n}(\mbb{A})$ be the algebra of $n\times n$ matrices
with entries in $\mbb{A}$. Denote by $E_{i,j}$ the $n\times n$
matrix with 1 as its $(i,j)$-entry and 0 as the others. Define the
vector space $\widehat{gl}(n,\mbb{A})=M_{n\times
n}(\mbb{A})\oplus\mbb{C}\kappa$ and its  Lie bracket:
\begin{eqnarray*} \hspace{1cm}& &[t^{m_1}\ptl_t^{r_1}E_{i_1,j_1}+\mu_1\kappa,t^{m_2}
\ptl_t^{r_2}E_{i_2,j_2}+\mu_1\kappa]\\
&=&\dlt_{j_1,i_2}t^{m_1}\ptl_t^{r_1}\cdot
t^{m_2}\ptl_t^{r_2}E_{i_1,j_2}-\dlt_{i_1,j_2}t^{m_2}\ptl_t^{r_2}\cdot
t^{m_1}\ptl_t^{r_1}E_{i_2,j_1}\\& &
+(-1)^{r_1}\dlt_{i_1,j_2}\dlt_{j_1,i_2}
\dlt_{r_1+r_2,m_1+m_2}r_1!r_2!\left(\!\!\begin{array}{c}m_1\\
r_1+r_2+1\end{array}
\!\!\right)\kappa\hspace{4.3cm}(1.2)\end{eqnarray*} for
$i,j\in\ol{1,n}$ and $m_1,m_2,r_1,r_2\in\mbb{N}$. Fixed an element
$\vec\ell=(\ell_1,...,\ell_n)\in\mbb{N}^{\:n}$. The subspace
$$\widehat{gl}(\vec\ell,\mbb{A})=\sum_{i,j=1}^n\mbb{A}\ptl_t^{\ell_j}E_{i,j}+\mbb{C}
\kappa\eqno(1.3)$$ of $\widehat{gl}(n,\mbb{A})$ forms a Lie
subalgebra. The well-known Lie algebra $W_{1+\infty}$ is the
special case of $\widehat{gl}(\vec\ell,\mbb{A})$ when $n=1$ and
$\vec\ell=0$.

Denote $\ol{1,n}=\{1,2,...,n\}$, and set
$i^\ast=n+1-i\;\for\;i\in\ol{1,n}.$ We fix $\es\in\{0,1\}$ and
take
$$\vec\ell=(\ell_1,\ell_2,...,\ell_n)\in\mbb{N}^{\:n}\;\;\mbox{such that}\;\;\{\ell_1,\ell_2,...,\ell_n\}\subset
2\mbb{N}+\es\eqno(1.4)$$ and $\ell_i=\ell_{i^\ast}$ for $
i\in\ol{1,n}.$ The subspace
$$\hat{o}(\vec\ell,\mbb{A})=\sum_{i,j=1}^n\:\sum_{r=0}^{\infty}
\sum_{m\in\mbb{Z}} \mbb{C}(t^m\ptl_t^{r+\ell_j}E_{i,j}
-(-1)^\es(-\ptl)^rt^m\ptl_t^{\ell_i}E_{j^\ast,i^\ast}) +
\mbb{C}\kappa\eqno(1.5)$$ of $\widehat{gl}(n,\mbb{A})$ forms a Lie
subalgebra. Next we suppose that $n=2n_0$ is an even positive
integer. Moreover, we define the parity of indices:
$$p(i)=0,\qquad
p(n_0+i)=1\qquad\for\;\;i\in\ol{1,n_0}.\eqno(1.6)$$ The subspace
$$\widehat{sp}(\vec\ell,\mbb{A})=\sum_{i,j=1}^n\:\sum_{r=0}^{\infty}
\sum_{m\in\mbb{Z}}\mbb{C}(t^m\ptl_t^{r+\ell_j}E_{i,j}-(-1)^{p(i)+p(j)+\es}(-\ptl_t)^rt^m
\ptl_t^{\ell_i}E_{j^\ast,i^\ast}) + \mbb{C}\kappa \eqno(1.7)$$ of
$\widehat{gl}(n,\mbb{A})$ forms a Lie subalgebra.

Note that the Lie algebras $\hat{o}(\vec\ell,\mbb{A})$ and
$\widehat{sp}(\vec\ell,\mbb{A})$ are in general not graded by
conformal weights. One of the main objectives in this paper is to
construct irreducible modules of
$\widehat{gl}(\vec\ell,\mbb{A}),\;\hat{o}(\vec\ell,\mbb{A})$ and
$\widehat{sp}(\vec\ell,\mbb{A})$ through weighted irreducible
modules (may not be highest weight type) of centrally-extended
classical Lie algebras of infinite matrices with finite number of
nonzero entries.

 For $i,j\in\ol{1,n}$ and
$r\in\mbb{N}$, we denote
$$E_{i,j}(r,z)=\sum_{l\in\mbb{Z}}^{\infty}t^l\ptl_t^rE_{i,j}z^{-l-1}.\eqno(1.8)$$
 The {\it vacuum module} ${\cal V}_\chi$ of
$\widehat{gl}(n,\mbb{A})$ is a module generated by a vector
$\vcm$, called {\it vacuum}, such that $\kappa|_{{\cal
V}_\chi}=\chi\mbox{Id}_{{\cal V}_\chi}$,
$$E_{i,j}(r,z)\vcm=\sum_{m=0}^{\infty}t^{-m-1}\ptl_t^rE_{i,j}\vcm z^m\eqno(1.9)$$
 for $i,j\in\ol{1,n},\;r\in\mbb{N}$, and any other
 $\widehat{gl}(n,\mbb{A})$-module with the same property must be a
 quotient module of ${\cal V}_\chi$. Denote by $U(\cdot)$ the
 universal enveloping algebra of a Lie algebra.
 Suppose that ${\cal G}$ is one of the Lie algebras
$\widehat{gl}(\vec\ell,\mbb{A})$, $\hat{o}(\vec\ell,\mbb{A})$ or
$\widehat{sp}(\vec\ell,\mbb{A})$. The ${\cal G}$-module
$${\cal V}_{\chi}({\cal G})=U({\cal G})\vcm\eqno(1.10)$$
is called the {\it vacuum module} of ${\cal G}$ and the
corresponding representation is called the {\it vacuum
representation} of ${\cal G}$. Our first main result in this paper
is as follows:\psp

{\bf Theorem 1.1}. {\it The module ${\cal V}_\chi({\cal G})$ is
irreducible if $\chi\not\in\mbb{Z}$. When $\chi\in\mbb{Z}$, the
module ${\cal V}({\cal G})$ has a unique maximal proper submodule
$\bar{\cal V}_\chi({\cal G})$, and the quotient $V_\chi({\cal
G})={\cal V}({\cal G})/\bar{\cal V}_\chi({\cal G})$ is an
irreducible ${\cal G}$-module. When $n>1$ and $\chi\in\mbb{N}$,
the submodule
$$\bar{\cal V}_\chi({\cal G})=
U({\cal G})(t^{-1}\ptl^{\ell_1}E_{n,1})^{\chi+1}\vcm\eqno(1.11)$$
for ${\cal G}=\widehat{gl}(\vec\ell,\mbb{A}),\;{\cal
G}=\hat{o}(\vec\ell,\mbb{A})$ with $\es=1$ and ${\cal
G}=\widehat{sp}(\vec\ell,\mbb{A})$ with $\es=0$.

If $n>3$, the submodule
$$\bar{\cal V}_\chi({\cal G})=
U({\cal G}
)(t^{-1}(\ptl^{\ell_1}E_{n-1,1}-\ptl_t^{\ell_2}E_{n,2}))^{\chi+1}
\vcm\eqno(1.12)$$ for ${\cal G}=\hat{o}(\vec\ell,\mbb{A})$ with
$\es=0$ and ${\cal G}=\widehat{sp}(\vec\ell,\mbb{A})$ with
$\es=1$.}\psp

Suppose that $\chi$ is a positive integer. Note
$(t^{-1}E_{n,1})^{\chi+1}\vcm$ is a singular vector generating the
maximal proper submodules of the vacuum modules at level $\chi$ of
the affine Kac-Moody algebras $\widehat{sl}(n,\mbb{C})$ and
$\widehat{sp}(n,\mbb{C})$, respectively. Moreover,
$(t^{-1}(E_{n-1,1}-E_{n,2}))^{\chi+1}\vcm$ is  a singular vector
generating the maximal proper submodules of the vacuum module at
level $\chi$ of the affine Kac-Moody algebra $\hat{o}(n,\mbb{C})$
when $n$ is even. Our above results are exactly analogous in
higher order differential operators to those of the affine
Kac-Moody algebras.

On ${\cal V}_\chi$, there exists a unique vertex algebra structure
whose structure map $Y(\cdot,z)$ satisfying
$Y(\vcm,z)=\mbox{Id}_{{\cal V}_\chi}$ and
$$Y(t^{-m-1}\ptl^rE_{i,j},z)=\frac{1}{m!}\frac{d^m}{dz^m}E_{i,j}(r,z)\eqno(1.13)
$$
for $i,j\in\ol{1,n}$ and $r,m\in\mbb{N}$. Let ${\cal G}$ be one of
the Lie algebras $\widehat{gl}(\vec\ell,\mbb{A})$,
$\hat{o}(\vec\ell,\mbb{A})$ or $\widehat{sp}(\vec\ell,\mbb{A})$.
The family $({\cal V}_\chi({\cal G}),Y(\cdot,z))$ forms a vertex
subalgebra . If $\chi\not\in \mbb{Z}$, the vertex algebra $({\cal
V}_\chi({\cal G} ),Y(\cdot,z))$ is simple. When $\chi\in\mbb{Z}$,
the quotient space $V_\chi({\cal G})$ forms a simple vertex
algebra. If $\ell_i\in\{0,1\}$ for $i\in\ol{1,n}$, the above
vertex algebras have a Virasoro element, and thus they are vertex
operator algebras.

\newcommand{\T}{{\cal T}}
\newcommand{\E}{{\cal E}}

Denote $\Z=\mbb{Z}+1/2.$
 Let $\ol{gl}(\infty)$ be a vector space with a basis $\{\E_{l,k}\mid l,k\in\Z\}$
 and multiplication:
$$ \E_{l_1,l_2}\cdot
\E_{k_1,k_2}=\dlt_{l_2+k_1,0}\E_{l_1,k_2}\qquad\for
\;\;l_1,l_2,k_1,k_2\in\Z.\eqno(1.14)$$ Then $\ol{gl}(\infty)$ is
isomorphic to the associative algebra of infinite matrices with
finite number of nonzero entries. Define the step function $H$ on
$\Z$ by
$$H(l)=\left\{\begin{array}{ll}1&\mbox{if}\;\;l>0,\\
0&\mbox{if}\;\;l<0\end{array}\right.\qquad\for\;\;l\in\Z.\eqno(1.15)$$
Set $\td{gl}(\infty)=\ol{gl}(\infty)\oplus\mbb{C}\kappa_0,$ where
$\kappa_0$ is symbol for base element. We have the following the
Lie bracket on $\td{gl}(\infty)$:
\begin{eqnarray*}\hspace{1cm}& &[\E_{l_1,l_2}+\mu_1\kappa_0,\E_{k_1,k_2}+\mu_2\kappa_0]
=\E_{l_1,l_2}\E_{k_1,k_2}- \E_{k_1,k_2}\E_{l_1,l_2}
\\ & &+\dlt_{l_1+k_2,0}\dlt_{l_2+k_1,0}[H(l_1)H(l_2)-H(k_1)H(k_2)]\kappa_0\hspace{5.6cm}
(1.16)\end{eqnarray*} for $l_1,l_2,k_1,k_2\in\Z$ and
$\mu_1,\mu_2\in\mbb{C}$. Moreover, the subspace
$\T=\sum_{l\in\Z}\mbb{C}\E_{l,-l}+\mbb{C}\kappa_0$ forms a toral
Cartan subalgebra of $\td{gl}(\infty)$. In fact, the root
structures of $\ol{gl}(\infty)$ and $\td{gl}(\infty)$ are the
same. Thus they have exactly the same representation theory.

Given $\mu\in\mbb{C}$, we set $\la\mu\ra_0=1$ and
$\la\mu\ra_m=\mu(\mu-1)\cdots(\mu-(m-1))$ for $ 0<m\in\mbb{N}.$
Denote by $\T^\ast$ the space of linear functions on $\T$. For
$\lmd\in\T^\ast$, we define
$$\mbox{supp}\:\lmd=\{l\in\Z\mid \lmd(\E_{l,-l})\neq
0\}.\eqno(1.17)$$ Pick a weight $\lmd$  such that
$\lmd(\kappa_0)=\chi$ and $\mbox{supp}\:\lmd$ is a finite set. Let
$\cal M$ be the highest weight irreducible
$\td{gl}(\infty)$-module with weight $\lmd$, whose highest weight
vector is annihilated by the subalgebra
$\mbox{Span}\:\{\E_{l,k}\mid l,k\in\Z;\;l+k>0\}.$ Fix a constant
$\iota\in\mbb{C}$. We construct a module structure of the vertex
algebra ${\cal V}_\chi$ whose structure map $Y^\iota_{\cal
M}(\cdot,z)$ satisfying $Y^\iota_{\cal M}(\vcm,z)=\mbox{Id}_{\cal
M}$ and
$$Y(t^{-1}\ptl^rE_{i,j},z)\equiv\sum_{l,k\in\mbb{Z}}\la
\iota-k\ra_r\E_{ln+i-1/2,kn-j+1/2}z^{-l-k-r-1}\qquad(\mbox{mod}\;\mbb{C}\kappa_0)\eqno(1.18)$$
if $\iota\not\in\mbb{Z}$, and
\begin{eqnarray*}  Y(t^{-1}\ptl^rE_{i,j},z)&\equiv&
\sum_{l,k=0}^n[\la
-k-1\ra_r\E_{(l-\iota)n+i-1/2,(k+\iota+1)n-j+1/2}z^{-l-k-\ell_i-r-2}\\
& &+\la k+\ell_j\ra_
r\E_{(l-\iota)n+i-1/2,(\iota-k)n-j+/2}z^{-l+k+\ell_j-\ell_i-r-1}\\
& &+\la k+\ell_j \ra_r\E_{-(l+\iota+1)n+i-1/2,(\iota-k)n-j+1/2}z^{l+k+\ell_j-r}\\
& & +\la-k-1\ra_r
\E_{-(l+\iota+1)n+i-1/2,(k+\iota+1)n-j+1/2}z^{l-k-r-1}]\;\;\;(\mbox{mod}\;\mbb{C}\kappa_0)
\hspace{0.6cm} (1.19)\end{eqnarray*} if $\iota\in\mbb{Z}$.

 Let $m$ is a positive integer. Set ${\cal
S}_m=\{\{3/2-r,5/2-r,...,(2m+1)/2-r\}\mid r\in\ol{1,m+1}\}$ the
set of intervals around $0$ of length $m$ in $\Z$. Define
\begin{eqnarray*}\hspace{1.5cm}\G^m&=&\{\lmd\in\T^\ast\mid\lmd(\kappa_0)=-m,\;
-s^{-1}|s|\lmd(\E_{s,-s})\in\mbb{N}\\ & &\for\;s\in\Z;\;
\mbox{supp}\:\lmd\subset S\;\mbox{for some}\;S\in{\cal S}_m\}.
\hspace{5cm}(1.20)\end{eqnarray*} \vspace{0.1cm}

{\bf Theorem 1.2}. {\it The family $({\cal M} ,Y^\iota_{\cal
M}(|_{{\cal V}_\chi({\cal G})} ,z))$ forms an irreducible module
of the vertex algebra $({\cal V}_\chi({\cal G}), Y(\cdot,z))$ for
${\cal G}=\widehat{gl}(\vec\ell,\mbb{A})$ and ${\cal
G}=\hat{o}(\vec\ell,\mbb{A}),\widehat{sp}(\vec\ell,\mbb{A})$ if
$\iota\not\in\mbb{Z}/2$.

Suppose that $\chi$ is a positive integer and
$$\lmd(\E_{1/2,-1/2}-\E_{-1/2,1/2}+\kappa_0),\lmd(\E_{l+1,-l-1}-\E_{l,-l})\in\mbb{N}\eqno(1.21)$$
for $-1/2\neq l\in \Z$. The family $({\cal M} ,Y^\iota_{\cal
M}(|_{{\cal V}_\chi({\cal G})} ,z))$ induces an irreducible module
of the quotient simple vertex algebra $ (V_\chi({\cal G}),
Y(\cdot,z))$ for ${\cal G}=\widehat{gl}(\vec\ell,\mbb{A})$ and
${\cal
G}=\hat{o}(\vec\ell,\mbb{A}),\widehat{sp}(\vec\ell,\mbb{A})$ if
$\iota\not\in\mbb{Z}/2$.

Assume that $\chi=-m$ is a negative integer and $\lmd\in\G^m$. The
family $({\cal M} ,Y^\iota_{\cal M}(|_{{\cal V}_\chi({\cal G})}
,z))$ induces an irreducible module of the quotient simple vertex
algebra $(V_{-m}({\cal G}), Y(\cdot,z))$ for ${\cal
G}=\widehat{gl}(\vec\ell,\mbb{A})$ and ${\cal
G}=\hat{o}(\vec\ell,\mbb{A}),\widehat{sp}(\vec\ell,\mbb{A})$ if
$\iota\not\in\mbb{Z}/2$.} \psp

We remark that the module ${\cal M}$ is an unitary module if
(1.21) is satisfied or $\lmd\in\G^m$. When $\iota\in\mbb{Z}/2$, we
construct irreducible modules of the vertex algebras $({\cal
V}_{\chi}(\hat{o}(\vec\ell,\mbb{A})),Y(\cdot,z))$ and $({\cal
V}_{\chi}(\widehat{sp}(\vec\ell,\mbb{A})),Y(\cdot,z))$ through
highest weight modules of non-standard centrally-extended the
other types of classical Lie algebras of infinite matrices with
finite number of nonzero entries. Similar conclusions for their
quotient simple vertex algebras hold. When $\chi$ is a positive
integer, our results on modules are natural generalizations of
those for the simple vertex operators algebras associated with
affine Lie algebras. When $n>3$ and $\chi$ is a positive integer,
the condition (1.21) can be relaxed to obtain certain non-unitary
irreducible modules of the simple vertex algebras
$V_\chi(\widehat{gl}(\vec\ell,\mbb{A})),\;V_\chi(\hat{o}(\vec\ell,\mbb{A}))$
and $V_\chi(\widehat{sp}(\vec\ell,\mbb{A}))$. In the case
$\vec\ell=\vec 0$, our theory coincide with the charged quadratic
free bosonic field theory if $\chi=-1$, and with the charged
quadratic free fermionic field theory if $\chi=1$. If
$\vec\ell=(1,1,...,1)$, our irreducible modules of the vertex
algebra $(V_{-1}(\hat{o}(\vec\ell,\mbb{A})),Y(\cdot,z))$ include
those studied by Dong and Nagatomo [DN1,DN2].

      The paper is organized as follows. In Section 2, we present the frame works of
constructing representations of certain Lie subalgebras of the
centrally-extended  Lie algebra of the tensor algebra of any
associative algebra with the algebra of differential operators on
the circle. Section 3 is devoted to the constructions of
irreducible modules of the Lie algebras
$\widehat{gl}(\vec\ell,\mbb{A})$, $\hat{o}(\vec\ell,\mbb{A})$ or
$\widehat{sp}(\vec\ell,\mbb{A})$ from weighted irreducible modules
of the centrally-extended  general linear Lie algebra of infinite
matrices with finite number of nonzero entries. In Section 4, we
give detailed constructions of irreducible modules with the
parameter $\iota\in \mbb{Z}+1/2$ of the Lie algebras
$\hat{o}(\vec\ell,\mbb{A})$ and $\widehat{sp}(\vec\ell,\mbb{A})$
 from weighted irreducible modules of certain  central extensions
of  the Lie algebras  of infinite skew  matrices with finite
number of nonzero entries. The cases with the parameter $\iota\in
\mbb{Z}$ are handled in Section 5. In Section 6, we study the
vacuum representation of the Lie algebra
$\widehat{gl}(\vec\ell,\mbb{A})$, its vertex algebra structure and
vertex algebra irreducible representations. We deal with the cases
for the Lie algebras $\hat{o}(\vec\ell,\mbb{A})$ and
$\widehat{sp}(\vec\ell,\mbb{A})$ in Section 7.

\section{General Frames}

In this section, we construct irreducible representations of
certain Lie subalgebras of the centrally-extended  Lie algebra of
the tensor algebra of any associative algebra with the algebra of
differential operators on the circle, from  representations of the
certain Lie subalgebras of a centrally-extended  Lie algebra of
the tensor algebra of the associative algebra with the algebra of
infinite  matrices with finite number of nonzero entries. When the
associative algebra is a finite matrix algebra, the later Lie
algebras are exactly centrally-extended classical Lie algebra of
infinite matrices with finite number of nonzero entries.

Recall
$$\ptl_t=\frac{d}{dt}\eqno(2.1)$$
and the algebra of differential operators on the circle:
$$\mbb{A}=\sum_{i=0}^{\infty}\mbb{C}[t,t^{-1}]\ptl_t^i\eqno(2.2)$$
 Note
$$ f(t)\ptl_t^i\cdot g(t)\ptl_t^j=\sum_{r=0}^i\left(\!\!\begin{array}{c}i\\ r\end{array}
\!\!\right)f(t)\frac{d^rg}{dt^r}(t)\ptl_t^{i+j-r}\qquad\for\;\;f(t),g(t)\in\mbb{C}[t,t^{-1}],
\;i,j\in\mbb{N}.\eqno(2.3)$$ Let ${\cal A}$ be an associative
algebra with an identity element $1_{\cal A}$ and a linear  map
$\tr:{\cal A}\rta\mbb{C}$ such that
$$\tr 1_{\cal A}\neq 0,\qquad\tr ab=\tr
ba\qquad\for\;\;a,b\in{\cal A}.\eqno(2.4)$$ Such a map ``$\tr$" is
called a {\it trace map}. Set

$$\hat{\cal A}={\cal A}\otimes_{\mbb{C}} \mbb{A}\oplus\mbb{C}\kappa,\eqno(2.5)$$
where $\kappa$ is a base element. According to [L] (also cf.
[M1]), we have the following Lie bracket on $\hat{\cal A}$:
\begin{eqnarray*} \qquad& &[a\otimes t^{i_1}\ptl_t^{j_1}+\mu_1\kappa,b\otimes
t^{i_2}\ptl_t^{j_2}+\mu_1\kappa]\\ &=&ab\otimes
t^{i_1}\ptl_t^{j_1}\cdot t^{i_2}\ptl_t^{j_2}-ba\otimes
t^{i_2}\ptl_t^{j_2}\cdot
t^{i_1}\ptl_t^{j_1}\\ & &+(-1)^{j_1}\dlt_{i_1+i_2,j_1+j_2}j_1!j_2!\left(\!\!\begin{array}{c}i_1\\
j_1+j_2+1\end{array} \!\!\right)(\tr
ab)\kappa\hspace{5.6cm}(2.6)\end{eqnarray*} for $a,b\in{\cal
A},\;i_1,i_2\in\mbb{Z},\;j_1,j_2\in\mbb{N}$ and
$\mu_1,\mu_2\in\mbb{C}$.

For two vector spaces $V_1$ and $V_2$, we denote by $LM(V_1,V_2)$
the space of linear maps from $V_1$ to $V_2$.  We  also use the
following operator of taking residue:
$$\rd_z(z^n)=\dlt_{n,-1}\qquad\for\;\;n\in \Bbb{Z}.\eqno(2.7)$$
Furthermore, all the binomials are assumed to be expanded in the
nonnegative powers of the second variable.

 A {\it conformal algebra} $R$ is a  $\Bbb{C}[\ptl]$-module with a linear map
 $Y^+(\cdot,z):\;R\rightarrow
  LM(R,R[z^{-1}]z^{-1})$ satisfying:
$$Y^+(\ptl u,z)={dY^+(u,z)\over dz}\qquad\for\;\;u\in R;\eqno(2.8)$$
$$Y^+(u,z)v=\rd_x{e^{x\ptl}Y^+(v,-x)u\over z-x},\eqno(2.9)$$
$$Y^+(u,z_1)Y^+(v,z_2)-Y^+(v,z_2)Y^+(u,z_1)=\rd_x{Y^+(Y^+(u,z_1-x)v,x)\over z_2-x}
 \eqno(2.10)$$
for $u,v\in R$. We denote by $(R,\ptl,Y^+(\cdot,z))$ a conformal
algebra.

Define
$$\hat{R}({\cal A})={\cal
A}[\vsi_1,\vsi_2]\oplus\mbb{C}{\bf 1},\eqno(2.11)$$ where
$\vsi_1,\vsi_2$ are indeterminates and ${\bf 1}$ is a base
element. For convenience, we denote
$$a[m_1,m_2]=a\vsi_1^{m_1}\vsi_2^{m_2}\qquad\for\;\;a\in{\cal
A},\;m_1,m_2\in\mbb{N}.\eqno(2.12)$$ We define the action of
$\mbb{C}[\ptl]$ on $\hat{R}({\cal A})$ by:
$$\ptl({\bf 1})=0,\;\;\ptl(a[m_1,m_2])=(m_1+1)a[m_1+1,m_2]+(m_2+1)a[m_1,m_2+1],\eqno(2.13)$$
and a linear map $Y^+(\cdot,z):\hat{R}({\cal A})\rta
LM(\hat{R}({\cal A}),\hat R({\cal A})[z^{-1}]z^{-1})$ by
$$Y^+(w,z){\bf 1}=Y^+({\bf 1},z)w=0\qquad \for\;\;w\in \hat R({\cal A})\eqno(2.14)$$
and
\begin{eqnarray*}& &
Y^+(a[m_1,m_2],z)b[n_1,n_2]\\
&=&\left(\!\!\begin{array}{c} -n_1-1\\ m_2\end{array}\!\!\right)
\sum_{p=0}^{m_1+m_2+n_1}\left(\!\!\begin{array}{c}p\\ m_1
\end{array}\!\!\right)ab[p,n_2]z^{p-m_1-m_2-n_1-1}\\
&&-\left(\!\!\begin{array}{c}-n_2-1\\
m_1\end{array}\!\!\right)\sum_{q=0}^{m_1+m_2+n_2}\left(\!\!\begin{array}{c}
q\\
m_2\end{array}\!\!\right)ba[n_1,q]z^{q-m_1-m_2-n_2-1}\\ &
&+\left(\!\!\begin{array}{c} -n_1-1\\
m_2\end{array}\!\!\right)\left(\!\!\begin{array}{c}-n_2-1\\
m_1\end{array}\!\!\right)(\tr ab){\bf
1}z^{-m_1-m_2-n_1-n_2-2}\hspace{5.3cm}(2.15)\end{eqnarray*} for
$a,b\in{\cal A}$ and $m_1,m_2,n_1,n_2\in\mbb{N}.$ Then $(\hat
R({\cal A}),\ptl,Y^+(\cdot,z))$ forms a conformal algebra (cf.
Section 7.3 in [X1]).

 Let
$$\hat{\cal A}[[z^{-1},z]]=\{\sum_{m\in\mbb{Z}}u_mz^m\mid
u_m\in\hat{\cal A}\}\eqno(2.16)$$ be the space of formal power
series with coefficients in $\hat{\cal A}$. Define a linear map
$Y(\cdot,z):\hat{R}({\cal A})\rta\hat{\cal A}[[z^{-1},z]]$ by
$$Y({\bf
1},z)=\kappa,\;\;Y(a[m_1,m_2],z)=\frac{1}{m_1!m_2!}\sum_{i\in\mbb{Z}}a\otimes
(-\ptl_t)^{m_1}t^i \ptl_t^{m_2}z^{-i-1}\eqno(2.17)$$ for
$a\in{\cal A}$ and $m_1,m_2\in\mbb{N}$. By Lemma 3.1 in [M1], we
have: \psp

{\bf Lemma 2.1}. {\it The Lie bracket (2.6) on $\hat{\cal A}$ is
equivalent to:
$$[Y(u,z_1),Y(v,z_2)]=\rd_{z_0}z_1^{-1}\dlt\left(\frac{z_2+z_0}{z_1}\right)
Y(Y^+(u,z_0)v,z_2)\eqno(2.18)$$ for $u,v\in\hat{R}({\cal A})$.}
\psp

Recall
$${\cal Z}=\mbb{Z}+\frac{1}{2}\eqno(2.19)$$
and the step function $H$ on $\Z$:
$$H(l)=\left\{\begin{array}{ll}1&\mbox{if}\;\;l>0,\\
0&\mbox{if}\;\;l<0\end{array}\right.\qquad\for\;\;l\in\Z.\eqno(2.20)$$
Set
$$\check{\cal A}=\sum_{i,j\in\Z}{\cal
A}t_1^it_2^j\oplus \mbb{C}\kappa_0,\eqno(2.21)$$ where $t_1,t_2$
are indeterminates and $\kappa_0$ is a base element. For
convenience, we denote
$$a(i,j)=at_1^it_2^j\qquad\for\;\;a\in{\cal
A},\;i,j\in\Z.\eqno(2.22)$$ According to Section 7.3 in [X2], we
have the following Lie bracket on $\check{\cal A}$:
\begin{eqnarray*} \qquad& &[a(l_1,l_2)+\mu_1\kappa_0,b(k_1,k_2)+\mu_2\kappa_0]
=\dlt_{l_2+k_1,0}ab(l_1,k_2)- \dlt_{l_1+k_2,0}ba(k_1,l_2)\\ &
&+\dlt_{l_1+k_2,0}\dlt_{l_2+k_1,0}[H(l_1)H(l_2)-H(k_1)H(k_2)](\tr
ab)\kappa_0\hspace{4.9cm}(2.23)
\end{eqnarray*}
for $a,b\in{\cal A},\;l_1,l_2,k_1,k_2\in\Z$ and
$\mu_1,\mu_2\in\mbb{C}$. The algebra $\check{\cal A}$ is
isomorphic to the corresponding central extension of the
commutator Lie algebra of ${\cal A}\otimes\ol{gl}(\infty)$ (cf.
(1.14) and the above).

Set
$$ \check{{\cal A}}^m=\mbox{Span}\:\{a(i,j)\mid a\in{\cal
A},\;i,j\in\Z;\;m<|i|,|j|;\;ij<0\}\qquad\for\;\;0<m\in\mbb{N}.\eqno(2.24)$$
 It can be verified that $\check{{\cal A}}^m$ is a Lie subalgebra of
$\check{\cal A}$. Suppose that ${\cal M}$ is an $\check{\cal
A}$-module
$$\mbox{generated by a subspace}\;\;{\cal M}_0\;\;\mbox{such
that}\;\;\check{\cal A}^m({\cal M}_0)=\{0\}\;\;\mbox{for
some}\;\;m\in\mbb{N}.\eqno(2.25)$$ Fixed a constant
$\iota\in\mbb{C}$. Define
\begin{eqnarray*}&
&\sum_{r_1,r_2=0}^{\infty}\Im_{r_1,r_2}x^{r_1}y^{r_2}z^{-r_1-r_2-1}=\frac{1}{x-y}\left(\left(
\frac{z+y}{z+x}\right)^\iota-1\right)\\
&=&-\left[\sum_{r=2}^{\infty}\left(\!\!\begin{array}{c}\iota\\ r
\end{array}\!\!\right)\frac{x^{r-1}+x^{r-2}y+\cdots
y^{r-1}}{z^r}+\frac{\iota}{z}\right]\left[\sum_{s=0}^{\infty}\left(\!\!\begin{array}{c}-\iota\\
s\end{array}\!\!\right)\left(\frac{x}{z}\right)^s\right].\hspace{3.5cm}(2.26)\end{eqnarray*}
Motivated from the construction of the twisted modules of spinor
vertex operator algebras in [X1], we define a linear map
$Y^\iota_{\cal M}(\cdot,z):\hat{R}({\cal A})\rta LM({\cal M},{\cal
M}[[z^{-1},z]])$ by
$$Y^\iota_{\cal M}({\bf 1},z)=\kappa_0,\eqno(2.27)$$
\begin{eqnarray*} Y^\iota_{\cal M}(a[r_1,r_2],z)&=&\sum_{i,j\in\Z}\left(\!\!\begin{array}{c}
-i-\iota-1/2\\
r_1\end{array}\!\!\right)\left(\!\!\begin{array}{c}-j+\iota-1/2\\
r_1\end{array}\!\!\right)a(i,j)z^{-i-j-r_1-r_2-1}\\ & &
+\Im_{r_1,r_2}(\tr
a)\kappa_0z^{-r_1-r_2-1}\hspace{7.4cm}(2.28)\end{eqnarray*} for
$a\in{\cal A}$ and $r_1,r_2\in\mbb{N}$. The above expression make
sense because of (2.25). Denote
$$\la\mu\ra_0=1,\;\;\la\mu\ra_m=\mu(\mu-1)\cdots(\mu-(m-1))\qquad\for\;\;\mu\in\mbb{C},\;
0<m\in\mbb{N}.\eqno(2.29)$$ \vspace{0.1cm}

{\bf Theorem 2.2}. {\it On the $\check{\cal A}$-module $\cal M$,
we have
$$[Y_{\cal M}^\iota(u,z_1),Y_{\cal M}^\iota(v,z_2)]=\rd_{z_0}z_1^{-1}
\dlt\left(\frac{z_2+z_0}{z_1}\right) Y_{\cal
M}^\iota(Y^+(u,z_0)v,z_2)\eqno(2.30)$$ for $u,v\in\hat{R}({\cal
A})$. In particular, Lemma 2.1  and (2.30) imply that ${\cal M}$
provides the following representation $\sgm_{\cal M}^\iota$ of the
Lie algebra $\hat{\cal A}$:
$$\sgm_{\cal M}^\iota(\kappa)=\kappa_0,\eqno(2.31)$$
$$\sgm_{\cal M}^\iota(a\otimes
t^m\ptl_t^r)=\sum_{l\in\Z}\la-l+\iota-1/2\ra_r
a(m-r-l,l)+r!\dlt_{r,m}\Im_{0,r}(\tr a)\kappa_0\eqno(2.32)$$ for $
m\in\mbb{Z},\;r\in\mbb{N}$ and $a\in{\cal A}$.}

{\it Proof}. We shall use generating functions to prove the lemma.
Set
$$
c[x,y]=\sum_{r_1,r_2=0}^{\infty}c[r_1,r_2]x^{r_1}y^{r_2}\qquad\for\;\;c\in{\cal
A}.\eqno(2.33)$$ Fix $a,b\in{\cal A}$ and let
$$u=a[x_1,x_2],\qquad v=b[y_1,y_2].\eqno(2.34)$$
Then
\begin{eqnarray*}Y^+(u,z_0)v&=&\frac{1}{z_0+x_2-y_1}ab[z_0+x_1,y_2]-\frac{1}{z_0+x_1-y_2}
ba[y_1,z_0+x_2]\\ & &+ \frac{1}{(z_0+x_2-y_1)(z_0+x_1-y_2)}(\tr
ab){\bf 1}\hspace{5.8cm}(2.35)\end{eqnarray*} by (2.15).

 On the other hand, we define
$$c(z_1,z_2)=\sum_{i,j\in\Z}c(i,j)z_1^{-i-\iota-1/2}z_2^{-j+\iota-1/2}\qquad\for\;\;c\in{\cal
A}.\eqno(2.36)$$ Then
$$Y^\iota_{\cal M}(u,z_1)=a(z_1+x_1,z_1+x_2)+\frac{1}{x_1-x_2}\left(\left(
\frac{z_1+x_2}{z_1+x_1}\right)^\iota-1\right)(\tr
a)\kappa_0,\eqno(2.37)$$ $$ Y^\iota_{\cal
M}(v,z_1)=b(z_2+y_1,z_2+y_2)+\frac{1}{y_1-y_2}\left(\left(
\frac{z_2+y_2}{z_2+y_1}\right)^\iota-1\right)(\tr
a)\kappa_0.\eqno(2.38)$$ Hence
\begin{eqnarray*}& &[Y^\iota_{\cal M}(u,z_1),Y^\iota_{\cal
M}(v,z_1)]\\
&=&[a(z_1+x_1,z_1+x_2),b(z_2+y_1,z_2+y_2)]\\
&=&(z_2+y_1)^{-1}\left(\frac{z_1+x_2}{z_2+y_1}\right)^{\iota}\dlt\left(\frac{z_1+x_2}
{z_2+y_1}\right) ab(z_1+x_1,z_2+y_2)\\ & &
-(z_2+y_2)^{-1}\left(\frac{z_2+y_2}{z_1+x_1}\right)^{\iota}
\dlt\left(\frac{z_1+x_1}{z_2+y_2}\right) ba(z_2+y_1,z_1+x_2)\\
&&+\left(\frac{(z_1+x_2)(z_2+y_2)}{(z_1+x_1)(z_2+y_1)}\right)^\iota[\frac{1}
{(z_1+x_1-z_2-y_2)(z_1+x_2-z_2-y_1)}\\ & &
-\frac{1}{(z_2+y_1-z_1-x_2)(z_2+y_2-z_1-x_1)}](\tr
ab)\kappa_0.\hspace{5.7cm}(2.39)\end{eqnarray*} Note that
\begin{eqnarray*}& &\left(\frac{(z_1+x_2)(z_2+y_2)}{(z_1+x_1)(z_2+y_1)}\right)^\iota
[\frac{1}
{(z_1+x_1-z_2-y_2)(z_1+x_2-z_2-y_1)}\\& &-\frac{1}{(z_2+y_1-z_1-x_2)(z_2+y_2-z_1-x_1)}]\\
&=&\frac{1}{x_1+y_1-x_2-y_2}\left(\frac{(z_1+x_2)(z_2+y_2)}{(z_1+x_1)(z_2+y_1)}\right)^\iota
[\frac{1}{z_1+x_2-z_2-y_1}\\ &
&-\frac{1}{z_1+x_1-z_2-y_2}+\frac{1}{z_2+y_1-z_1-x_2}-\frac{1}{z_2+y_2-z_1-x_1}]\\
&=&\frac{z_1^{-1}}{x_1+y_1-x_2-y_2}\left(\frac{(z_1+x_2)(z_2+y_2)}{(z_1+x_1)
(z_2+y_1)}\right)^\iota
\\ & &\times \left[\dlt\left(\frac{z_2+y_1-x_2}{z_1}\right)-\dlt\left(\frac{z_2+y_2-x_1}{z_1}
\right)\right]
\\&=&\frac{z_1^{-1}}{x_1+y_1-x_2-y_2}[\left(\frac{z_2+y_2}{z_2+x_1+y_1-x_2}\right)^\iota
\dlt\left(\frac{z_2+y_1-x_2}{z_1}\right)\\ & &
-\left(\frac{z_2+x_2+y_2-x_1}{z_2+y_1}\right)^\iota\dlt\left(\frac{z_2+y_2-x_1}{z_1}\right)
].\hspace{6.3cm}(2.40)\end{eqnarray*}

Observe that
\begin{eqnarray*}(z_2+y_1)^{-1}\dlt\left(\frac{z_1+x_2}{z_2+y_1}\right)&=&
\frac{1}{z_1+x_2-z_2-y_1} +\frac{1}{z_2+y_1-z_1-x_2}\\ &
=&z_1^{-1}\dlt\left(\frac{z_2+y_1-x_2}{z_1}\right).\hspace{6cm}(2.41)\end{eqnarray*}
 Hence
\begin{eqnarray*}& &\rd_{z_0}z_1^{-1}\dlt\left(\frac{z_2+z_0}{z_1}\right)\frac{1}{z_0+x_2-y_1}
ab(z_2+z_0+x_1,z_2+y_2)\\
&=&\rd_{z_0}\frac{(z_2+y_1)^{-\iota}}{z_0+x_2-y_1}z_1^{-1}\dlt\left(\frac{z_2+z_0}{z_1}\right)
[(z_2+y_1)^\iota ab(z_2+z_0+x_1,z_2+y_2)]\\
&=&\rd_{z_0}\frac{(z_2+y_1)^{-\iota}}{z_0+x_2-y_1}z_1^{-1}\dlt\left(\frac{z_2+z_0}{z_1}\right)
[(z_1-z_0+y_1)^\iota ab(z_1+x_1,z_2+y_2)]\\ &=&
\left(\frac{z_1+x_2}{z_2+y_1}\right)^\iota
z_1^{-1}\dlt\left(\frac{z_2+y_1-x_2}{z_1}\right)
ab(z_1+x_1,z_2+y_2)\\
&=&(z_2+y_1)^{-1}\left(\frac{z_1+x_2}{z_2+y_1}\right)^\iota\dlt\left(\frac{z_1+x_2}
{z_2+y_1}\right)
ab(z_1+x_1,z_2+y_2).\hspace{4.4cm}(2.42)\end{eqnarray*} Similarly,
we have
\begin{eqnarray*}&
&\rd_{z_0}z_1^{-1}\dlt\left(\frac{z_2+z_0}{z_1}\right)\frac{1}{z_0+x_1-y_2}ba(z_2+y_1,
z_2+z_0+x_2)\\
&=&(z_2+y_2)^{-1}\left(\frac{z_2+y_2}{z_1+x_1}\right)^{\iota}
\dlt\left(\frac{z_1+x_1}{z_2+y_2}\right)
ba(z_2+y_1,z_1+x_2).\hspace{4.5cm}(2.43)\end{eqnarray*} Moreover,
\begin{eqnarray*}& &\rd_{z_0}z_1^{-1}\dlt\left(\frac{z_2+z_0}{z_1}\right)\frac{1}{(z_0+x_2-y_1)
(z_0+x_1-y_2)}\\ & &\times
\left[\left(\frac{z_2+y_2}{z_2+z_0+x_1}\right)^\iota
+\left(\frac{z_2+z_0+x_2}{z_2+y_1}\right)^\iota-1\right]\\
&=&\frac{z_1^{-1}}{x_1+y_1-x_2-y_2}\{\rd_{z_0}\frac{1}{z_0+x_2-y_1}
\dlt\left(\frac{z_2+z_0}{z_1}\right)\\ & &\times
\left[\left(\frac{z_2+y_2}{z_2+z_0+x_1}\right)^\iota
+\left(\frac{z_2+z_0+x_2}{z_2+y_1}\right)^\iota-1\right]\\ &
&-\rd_{z_0}\frac{1}{z_0+x_1-y_2}\dlt\left(\frac{z_2+z_0}{z_1}\right)
\left[\left(\frac{z_2+y_2}{z_2+z_0+x_1}\right)^\iota
+\left(\frac{z_2+z_0+x_2}{z_2+y_1}\right)^\iota-1\right]\}\\ &=&
\frac{z_1^{-1}}{x_1+y_1-x_2-y_2}[\dlt\left(\frac{z_2+y_1-x_2}{z_1}\right)
\left(\frac{z_2+y_2}{z_2+x_1+y_1-x_2}\right)^\iota\\ & &-
\dlt\left(\frac{z_2+y_2-x_1}{z_1}\right)\left(\frac{z_2+x_2+y_2-x_1}{z_2+y_1}\right)^\iota
].\hspace{6.9cm}(2.44)\end{eqnarray*}

By (2.35), (2.39), (2.40) and (2.42)-(2.44), we have
\begin{eqnarray*}& &\rd_{z_0}z_1^{-1}
\dlt\left(\frac{z_2+z_0}{z_1}\right) Y_{\cal
M}^\iota(Y^+(u,z_0)v,z_2)\\ &=&\rd_{z_0}z_1^{-1}
\dlt\left(\frac{z_2+z_0}{z_1}\right)[\frac{1}{z_0+x_2-y_1}Y_{\cal
M}^\iota(ab[z_0+x_1,y_2],z_2) \\ & &-\frac{1}{z_0+x_1-y_2}Y_{\cal
M}^\iota(ba[y_1,z_0+x_2],z_2)\\ & &+
\frac{1}{(z_0+x_2-y_1)(z_0+x_1-y_2)}(\tr ab)\kappa_0]\\ &=&
\rd_{z_0}z_1^{-1}
\dlt\left(\frac{z_2+z_0}{z_1}\right)\{\frac{1}{z_0+x_2-y_1}[ab(z_2+z_0+x_1,z_2+y_2)
\\ & &+\frac{1}{z_0+x_1-y_2}
\left[\left(\frac{z_2+y_2}{z_2+z_0+x_1}\right)^\iota-1\right](\tr
ab)\kappa_0]\hspace{9cm}\end{eqnarray*}
\begin{eqnarray*} & &-\frac{1}{z_0+x_1-y_2}[ba(z_2+y_1,z_2+z_0+x_2)-\frac{1}{z_0+x_2-y_1}\\
& &\times
\left[\left(\frac{z_2+z_0+x_2}{z_2+y_1}\right)^\iota-1\right](\tr
ab)\kappa_0]+ \frac{1}{(z_0+x_2-y_1)(z_0+x_1-y_2)}(\tr
ab)\kappa_0\}\\ &=&\rd_{z_0}z_1^{-1}
\dlt\left(\frac{z_2+z_0}{z_1}\right)\frac{1}{z_0+x_2-y_1}ab(z_2+z_0+x_1,z_2+y_2)\\
& &-\rd_{z_0}z_1^{-1}
\dlt\left(\frac{z_2+z_0}{z_1}\right)\frac{1}{z_0+x_1-y_2}ba(z_2+y_1,z_2+z_0+x_2)\\
&
&+\rd_{z_0}z_1^{-1}\dlt\left(\frac{z_2+z_0}{z_1}\right)\frac{1}{(z_0+x_2-y_1)
(z_0+x_1-y_2)}\\ & &\times
\left[\left(\frac{z_2+y_2}{z_2+z_0+x_1}\right)^\iota
+\left(\frac{z_2+z_0+x_2}{z_2+y_1}\right)^\iota-1\right]\\
&=&[Y_{\cal M}^\iota(u,z_1),Y_{\cal
M}^\iota(v,z_2)],\hspace{10cm}(2.45)\end{eqnarray*} \vspace{0.1cm}
that is (2.30) holds. Moreover, (2.32) follows from (2.17), Lemma
2.1 and (2.30).$\qquad\Box$ \psp

 Suppose that the associative algebra ${\cal A}$ has $n$ left
ideals $\{{\cal B}_1,{\cal B}_2,...,{\cal B}_n\}$ such that
$${\cal A}=\bigoplus_{i=1}^n{\cal B}_i.\eqno(2.46)$$
Take
$$\vec\ell=(\ell_1,\ell_2,...,\ell_n)\in\mbb{N}^{\:n}\eqno(2.47)$$
and set
$$\hat{\cal A}_{\vec\ell}=\sum_{i=1}^n{\cal
B}_i\otimes\mbb{A}\ptl_t^{\ell_i}+\mbb{C}\kappa.\eqno(2.48)$$ Then
$\hat{\cal A}_{\vec\ell}$ forms a Lie subalgebra of $\hat{\cal A}$
(cf. (2.5) and (2.6)). Let  ${\cal M}$ be an $\check{\cal
A}$-module satisfying (2.25). Consider the restricted
representation $\sgm_{\cal M}^\iota|_{\hat{\cal A}_{\vec\ell}}$ of
$\hat{\cal A}_{\vec\ell}$. However, when $\iota\in\mbb{Z}$ and
$\vec\ell\neq\vec 0$, the summation in (2.32) has redundant terms
in terms of constructing irreducible modules. We want to represent
the representation in reduced form. Set
$$a(r,z)=\sum_{m\in\mbb{Z}}a\otimes
t^m\ptl_t^rz^{-m-1}\qquad\mbox{for}\;\;a\in{\cal
A},\;r\in\mbb{N}.\eqno(2.49)$$ Suppose that $\iota\in\mbb{Z}$  and
each
$${\cal B}_j=\bigoplus_{i=1}^n{\cal B}_{i,j}\eqno(2.50)$$
has a subspace decomposition such that
$$ {\cal B}_{i,j}\cdot {\cal
B}_{l,k}=\{0\}\qquad\mbox{if}\;\;j\neq l.\eqno(2.51)$$ By deleting
the redundant terms and changing indices, we obtain: \psp

{\bf Theorem 2.3}. {\it On the $\check{\cal A}$-module ${\cal M}$,
we have the following representation $\sgm_{\vec\ell,{\cal M}}$ of
$\hat{\cal A}_{\vec\ell}$: $\sgm_{\vec\ell,{\cal
M}}(\kappa)=\kappa_0$ and
\begin{eqnarray*} & & \sgm_{\vec\ell,{\cal M}}(a(r,z))\\&=&
\sum_{0<l_1,l_2\in\Z}[\la -l_2-1/2\ra_r
a(l_1-\iota,l_2+\iota)z^{-l_1-l_2-\ell_i-r-1}+\la
l_2+\ell_j-1/2\ra_ r a(l_1-\iota,\iota-l_2)\\ & &\times
z^{-l_1+l_2+\ell_j-\ell_i-r-1}+\la-l_2-1/2\ra_r
a(-l_1-\iota,l_2+\iota) z^{l_1-l_2-r-1}\\ &&+\la l_2+\ell_j
-1/2\ra_
ra(-l_1-\iota,\iota-l_2)z^{l_1+l_2+\ell_j-r-1}]+r!\Im_{0,r}(\tr
a)\kappa_0z^{-r-1}\hspace{2.3cm}(2.52)\end{eqnarray*} for $a\in
{\cal B}_{i,j}$ and $r\in\mbb{N}$.} \psp

Suppose that $\tau$ is a linear transformation on ${\cal A}$ such
that
$$\tau^2=\mbox{Id}_{\cal A},\;\;\tau(1_{\cal A})=1_{\cal
A},\;\;\tau(ab)=\tau(b)\tau(a)\qquad\for\;\;a,b\in{\cal
A},\eqno(2.53)$$
$$\tr \tau=\tr,\qquad\tau({\cal B}_{i,j})\subset
{\cal B}_{\pi(j),\pi(i)}\eqno(2.54)$$ for $i,j\in\ol{1,n}$, where
$$\pi\;\;\mbox{is a permutation
on}\;\;\{1,2,...,n\}.\eqno(2.55)$$ Let
$$\es\in\{0,1\}\eqno(2.56)$$
and take $\vec\ell\in\mbb{N}^{\:n}$ such that
$$\{\ell_1,\ell_2,...,\ell_n\}\subset 2\mbb{Z}+\es\eqno(2.57)$$
and
$$\ell_i=\ell_{\pi(i)}\qquad\for\;\;i\in\ol{1,n}.\eqno(2.58)$$
 Set
\begin{eqnarray*}\qquad\hat{\cal A}_{\vec\ell}^\tau&=&\mbox{Span}\:\{a\otimes
t^m\ptl_t^{r+\ell_j}-(-1)^\es\tau(a)\otimes
(-\ptl_t)^rt^m\ptl_t^{\ell_i}\mid i,j\in\ol{1,n};\\ & &a\in{\cal
B}_{i,j};\;
r\in\mbb{N},\;m\in\mbb{Z}\}+\mbb{C}\kappa.\hspace{7.7cm}(2.59)\end{eqnarray*}
It can be verified that $\hat{\cal A}_{\vec\ell}^\tau$ forms a Lie
subalgebra of $\hat{\cal A}$ (cf. (2.6), [M1]). For any
$\check{\cal A}$-module ${\cal M}$ satisfying (2.25), we have the
restricted representation $\sgm_{\cal M}^\iota|_{\hat{\cal
A}_{\vec\ell}^\tau}$. For convenience, we set
$$a^\tau_{\vec{\ell}}(r,z)=\sum_{m\in\mbb{Z}}(a\otimes
t^m\ptl_t^{r+\ell_j}-(-1)^\es\tau(a)\otimes
(-\ptl_t)^rt^m\ptl_t^{\ell_i})z^{-m-1}\eqno(2.60)$$ Then for $
a\in{\cal B}_{i,j}$  and $r\in\mbb{N},$ we have
\begin{eqnarray*}&&\sgm_{\cal
M}^\iota(a^\tau_{\vec{\ell}}(r,z))=Y^\iota_{\cal
M}((r+\ell_j)!a[0,r+\ell_j]-(-1)^\es r!\ell_i! \tau(a)[r,\ell_i],z)\\
&=&\sum_{l_1,l_2\in\Z}\la-l_2+\iota-1/2\ra_{
r+\ell_j}a(l_1,l_2)z^{-l_1-l_2-r-\ell_j-1}+(r+\ell_j)!\Im_{0,r+\ell_j}(\tr
a)\kappa_0z^{-r-\ell_j-1}\\ &&-(-1)^\es
\sum_{l_1,l_2\in\Z}\la-l_1-\iota-1/2\ra_ r\la-l_2+\iota-1/2\ra_{
\ell_i}\tau(a)(l_1,l_2)z^{-l_1-l_2-r-\ell_i-1}\\ & &-(-1)^\es
r!\ell_i! \Im_{r,\ell_i}(\tr
a)\kappa_0z^{-r-\ell_i-1}\hspace{10cm}\end{eqnarray*}
\begin{eqnarray*}
&=&\sum_{m\in\mbb{Z}}\:\sum_{l\in\Z}[\la l+r+\ell_j+\iota-1/2\ra_
{r+\ell_j}a(m+l,-l-r-\ell_j)-(-1)^\es\la l+r-\iota-1/2\ra_r\\ &
&\times\la-m-l+\ell_i+\iota-1/2\ra_{
\ell_i}\tau(a)(-l-r,m+l-\ell_i)]z^{-m-1}\\
&&+(r+\ell_j)!\Im_{0,r+\ell_j}(\tr
a)\kappa_0z^{-r-\ell_j-1}-(-1)^\es r!\ell_i!\Im_{r,\ell_i}(\tr
a)\kappa_0z^{-r-\ell_i-1}\hspace{2.7cm}(2.61)\end{eqnarray*} for
$a\in{\cal B}_{i,j}$  and $r\in\mbb{N}$ by (2.17) and (2.18). Thus
\begin{eqnarray*}& &\sgm_{\cal M}^\iota(a\otimes
t^m\ptl_t^{r+\ell_j}-(-1)^\es\tau(a)\otimes(-\ptl_t)^rt^m\ptl_t^{\ell_i})\\
&=&\sum_{l\in\Z}[\la l+r+\ell_j+\iota-1/2\ra_{r+\ell_j}
a(m+l,-l-r-\ell_j)\\ & &-(-1)^\es\la l+r-\iota-1/2\ra_r
\la-m-l+\ell_i+\iota-1/2\ra_{\ell_i}\tau(a)(-r-l,m+l-\ell_i)]\\
& &+[(r+\ell_j)!\dlt_{m,r+\ell_j}\Im_{0,r+\ell_j}-(-1)^\es
r!\ell_i!\dlt_{m,r+\ell_i}\Im_{r,\ell_i}](\tr
a)\kappa_0.\hspace{4.3cm}(2.62)\end{eqnarray*}

Suppose $\iota\in\mbb{Z}/2$. Then the restricted representation of
$\hat{\cal A}_{\vec\ell}^\tau$ is in general reducible. In (2.62),
\begin{eqnarray*}& &\sum_{l\in\Z}[\la l+r+\ell_j+\iota-1/2\ra_{r+\ell_j}
a(m+l,-l-r-\ell_j)\\ & &-(-1)^\es\la l+r-\iota-1/2\ra_r
\la-m-l+\ell_i+\iota-1/2\ra_{\ell_i}\tau(a)(-r-l,m+l-\ell_i)]\\
&=&\sum_{l\in\Z}[\la l+r+\ell_j+\iota-1/2\ra_{r+\ell_j}
a(m+l,-l-r-\ell_j)\\ & &-(-1)^\es\la l+r+\iota-1/2\ra_r
\la-m-l+\ell_i-\iota-1/2\ra_{\ell_i}\tau(a)(-r-l-2\iota,m+l+2\iota-\ell_i)]
\\ &=&\sum_{l\in\Z}\la l+r+\iota-1/2\ra_r[\la l+r+\ell_j+\iota-1/2\ra_{\ell_j}
a(m+l,-l-r-\ell_j)\\ & &-(-1)^\es
\la-m-l+\ell_i-\iota-1/2\ra_{\ell_i}\tau(a)(-r-l-2\iota,m+l+2\iota-\ell_i)]
\\ &=&\sum_{l\in\Z}\la l+r+\iota-1/2\ra_r[(-1)^\es\la -l-r-\iota-1/2\ra_{\ell_j}
a(m+l,-l-r-\ell_j)\\ & &- \la
m+l+\iota-1/2\ra_{\ell_i}\tau(a)(-r-l-2\iota,m+l+2\iota-\ell_i)].
\hspace{4.8cm}(2.63)\end{eqnarray*}
 Set
$$a^{\tau,\iota}_{\vec\ell}(l_1,l_2)=(-1)^\es\la
l_2+\iota-1/2\ra_{\ell_j}a(l_1,l_2+2\iota-\ell_j)-\la
l_1+\iota-1/2\ra_{\ell_i}
\tau(a)(l_2,l_1+2\iota-\ell_i)\eqno(2.64)$$ for $a\in{\cal
B}_{i,j}$ and $ l_1,l_2\in\Z$. Then
$$a^{\tau,\iota}_{\vec\ell}(l_1,l_2)=-(-1)^\es\tau(a)^{\tau,\iota}_{\vec\ell}(l_2,l_1).
\eqno(2.65)$$ Moreover, (2.61) becomes
\begin{eqnarray*}\sgm_{\cal
M}^\iota(a^\tau_{\vec{\ell}}(r,z))&=&\sum_{l_1,l_2\in\Z}\la
-l_2-\iota-1/2\ra_ra^{\tau,\iota}_{\vec\ell}(l_1,l_2)z^{-l_1-l_2-2\iota-r-1}
+(r+\ell_j)!\Im_{0,r+\ell_j}(\tr a)\kappa_0\\ & &\times
z^{-r-\ell_j-1} -(-1)^\es r!\ell_i!\Im_{r,\ell_i}(\tr
a)\kappa_0z^{-r-\ell_i-1}.\hspace{4.4cm}(2.66)\end{eqnarray*} For
$l,k\in\Z$ and $\ell\in\mbb{N}$, we have
$$\la l-1/2\ra_\ell=(-1)^{\ell}\la
k-1/2\ra_\ell\qquad\mbox{if}\;\;l+k=\ell.\eqno(2.67)$$ Given
$a\in{\cal B}_{i_1,j_1},\;b\in{\cal B}_{i_2,j_2}$   and
$l_1,l_2,k_1,k_2\in\Z$, (2.23) and (2.64) imply
\begin{eqnarray*}&
&[a^{\tau,\iota}_{\vec\ell}(l_1,l_2),b^{\tau,\iota}_{\vec\ell}(k_1,k_2)]\\
&=&[(-1)^\es\la
l_2+\iota-1/2\ra_{\ell_{j_1}}a(l_1,l_2+2\iota-\ell_{j_1})-\la
l_1+\iota-1/2\ra_{\ell_{i_1}} \tau(a)(l_2,l_1+2\iota-\ell_{i_1}),\\
& &(-1)^\es\la
k_2+\iota-1/2\ra_{\ell_{j_2}}b(k_1,k_2+2\iota-\ell_{j_2})-\la
k_1+\iota-1/2\ra_{\ell_{i_2}} \tau(b)(k_2,k_1+2\iota-\ell_{i_2})]\\
&\equiv &\la k_1+\iota-1/2\ra_{\ell_{i_2}}[
\dlt_{j_1,i_2}\dlt_{l_2+k_1+2\iota,\ell_{j_1}}(ab)^{\tau,\iota}_{\vec\ell}(l_1,k_2)+
\dlt_{i_2,i_1}\dlt_{l_1+k_1+2\iota,\ell_{j_1}}(\tau(b)a)^{\tau,\iota}_{\vec\ell}(k_2,l_2)]
\\&&-\la l_1+\iota-1/2\ra_{\ell_{i_1}}
[\dlt_{j_2,i_1}\dlt_{l_1+k_2+2\iota,\ell_{j_2}}(ba)^{\tau,\iota}_{\vec\ell}(k_1,l_2)\\
&
&+\dlt_{i_1,i_2}\dlt_{l_1+k_1+2\iota,\ell_{j_2}}(\tau(a)b)^{\tau,\iota}_{\vec\ell}(l_2,k_2)]
\;\;(\mbox{mod}\;\mbb{C}\kappa_0).\hspace{6cm}(2.68)\end{eqnarray*}

We define
$$\check{\cal
A}_{\vec\ell}^{\tau,\iota}=\mbox{Span}\:\{a^{\tau,\iota}_{\vec\ell}(l_1,l_2)\\
 \mid i,j\in\ol{1,n};\;a\in{\cal
B}_{i,j},\;l_1,l_2\in\Z\}+\mbb{C}\kappa_0.\eqno(2.69)$$ Then
 $\check{\cal A}_{\vec\ell}^{\tau,\iota}$ is a Lie subalgebra of
$\check{\cal A}$ by (2.68). Set
$$\td\ell=|2\iota|+\mbox{max}\:\{\ell_1,\ell_2,...,\ell_n\}.\eqno(2.70)$$
Define
$$\check{\cal
A}_{\vec\ell,m}^{\tau,\iota}=\mbox{Span}\:\{a^{\tau,\iota}_{\vec\ell}(l_1,-l_2)\\
 \mid i,j\in\ol{1,n};\;a\in{\cal
B}_{i,j},\;m<l_1,l_2\in\Z\}\qquad\for\;\;\td\ell
<m\in\mbb{N}.\eqno(2.71)$$ By (2.23), (2.64), (2.67) and (2.68),
$\check{\cal A}_{\vec\ell,m}^{\tau,\iota}$ forms a Lie subalgebra
of $\check{\cal A}_{\vec\ell}^{\tau,\iota}$.\psp

{\bf Theorem 2.4}. {\it Let $\iota\in\mbb{Z}/2.$ Let $M$ be a
$\check{\cal A}_{\vec\ell}^{\tau,\iota}$-module
$$\mbox{\it generated by a subspace}\;\;M_0\;\;\mbox{\it such
that}\;\;\check{\cal
A}_{\vec\ell,m}^{\tau,\iota}(M_0)=\{0\}\;\;\mbox{\it for
some}\;\;\td\ell<m\in\mbb{N}.\eqno(2.72)$$ Then we have the
following representation $\sgm_M$ of $\hat{\cal
A}_{\vec\ell}^\tau$ on $M$: $\sgm_M(\kappa)=\kappa_0$ and
\begin{eqnarray*}\sgm(a^\tau_{\vec{\ell}}(r,z))&=& \sum_{l_1,l_2\in\Z}\la
-l_2-\iota-1/2\ra_ra^{\tau,\iota}_{\vec\ell}(l_1,l_2)z^{-l_1-l_2-2\iota-r-1}
+(r+\ell_j)!\Im_{0,r+\ell_j}(\tr a)\kappa_0\\ & &\times
z^{-r-\ell_j-1} -(-1)^\es r!\ell_i!\Im_{r,\ell_i}(\tr
a)\kappa_0z^{-r-\ell_i-1}\hspace{5cm}(2.73)\end{eqnarray*}for
$a\in{\cal B}_{i,j}$  and $r\in\mbb{N}.$ }\psp

Next we assume $\iota\in\mbb{Z}$. When $\vec\ell\neq\vec 0$, the
summation in (2.73) has redundant terms in terms of constructing
irreducible modules. We want to represent the representation in
reduced form. It can be verified that the subspace $${\cal
I}=\mbox{Span}\:\{a^{\tau,\iota}_{\vec\ell}(l_1-\iota,l_2-\iota)\mid
i,j\in\ol{1,n};\;a\in{\cal
B}_{i,j},\;l_1,l_2\in\Z,\;0<\l_1<\ell_i\;\mbox{or}\;0<l_2<\ell_j\}\eqno(2.74)$$
forms an ideal of $\check{\cal A}_{\vec\ell}^{\tau,\iota}$.

Denote
$$\Z_i=\Z\setminus\{-\iota+1/2,-\iota+3/2...,-\iota+\ell_i-1/2\}\qquad\for\;\;i\in\ol{1,n}.
\eqno(2.75)$$ The the subspace
$$\check{\cal
D}_{\vec\ell}^{\tau,\iota}=\mbox{Span}\:\{a^{\tau,\iota}_{\vec\ell}(l_1,l_2)\mid
i,j\in\ol{1,n};\;a\in{\cal
B}_{i,j},\;l_1\in\Z_i,\;l_2\in\Z_j\}+\mbb{C}\kappa_0\eqno(2.76)$$
forms a Lie subalgebra of $\check{\cal
A}_{\vec\ell}^{\tau,\iota}$. Moreover,
$$\check{\cal A}_{\vec\ell}^{\tau,\iota}=\check{\cal
D}_{\vec\ell}^{\tau,\iota}\oplus {\cal I}.\eqno(2.77)$$ Set
$$\check{\cal
D}_{\vec\ell,m}^{\tau,\iota}=\mbox{Span}\:\{a^{\tau,\iota}_{\vec\ell}(l_1,-l_2)\mid
i,j\in\ol{1,n};\;a\in{\cal
B}_{i,j},\;m<l_1\in\Z_i,\;m<l_2\in\Z_j\}\eqno(2.78)$$ for $
\td\ell<m\in\mbb{N}.$ By (2.23), (2.64) and (2.68), $\check{\cal
D}_{\vec\ell,m}^{\tau,\iota}$ forms a Lie subalgebra of
$\check{\cal D}_{\vec\ell}^{\tau,\iota}$. According to (2.73), we
have: \psp

{\bf Theorem 2.5}. {\it Suppose $\iota\in\mbb{Z}$. Let ${\cal N}$
be a $\check{\cal D}_{\vec\ell}^{\tau,\iota}$-module
$$\mbox{\it generated by a subspace}\;\;{\cal N}_0\;\;\mbox{\it such
that}\;\;\check{\cal D}_{\vec\ell,m}^{\tau,\iota}({\cal
N}_0)=\{0\}\;\;\mbox{\it for some}\;\;\td\ell <
m\in\mbb{N}.\eqno(2.79)$$ Then we have the following
representation $\sgm_{\cal N}$ of $\hat{\cal A}_{\vec\ell}^\tau$
on $\cal N$: $\sgm_{\cal N}(\kappa)=\kappa_0$ and
\begin{eqnarray*}\sgm_{\cal N}(a^\tau_{\vec{\ell}}(r,z))&=&\sum_{l_1\in\Z_i,l_2\in\Z_j}\la
-l_2-\iota-1/2\ra_ra^{\tau,\iota}_{\vec\ell}(l_1,l_2)z^{-l_1-l_2-2\iota-r-1}
+(r+\ell_j)!\Im_{0,r+\ell_j}(\tr a)\kappa_0\\ & &\times
z^{-r-\ell_j-1} -(-1)^\es r!\ell_i!\Im_{r,\ell_i}(\tr
a)\kappa_0z^{-r-\ell_i-1}\hspace{4.8cm}(2.80)\end{eqnarray*} for
$a\in{\cal B}_{i,j}$ and $r\in\mbb{N}.$ }\psp

{\bf Example 2.1}. Let $k>1$ be integer. The {\it Hecke algebra}
${\cal H}_k$ is an associative algebra generated by
$\{T_1,...,T_{k-1}\}$ with the following defining relations
$$T_iT_j=T_jT_i\qquad\mbox{whenever}\;\;|i-j|\geq 2,\eqno(2.81)$$
$$T_iT_{i+1}T_i=T_{i+1}T_iT_{i+1},\;\;T_i^2=(q-1)T_i+q\eqno(2.82)$$
for $i,j\in \ol{1,k-1}$, where $0\neq q\in \Bbb{C}$. Let $\zeta\in
\Bbb{C}$ be a fixed constant.  According to Section 5 in [HKW],
there exists a unique trace map ``$\tr$'' of ${\cal H}_k$ such
that
$$\tr(e)={1\over 2},\;\;\tr(a T_mb)=\zeta \tr(ab)\qquad\mbox{for}\;\;a,b\in H_m\eqno(2.83)$$
with $m\in \ol{1,k-2}$. This trace map is the key to define the
well-known ``Jones polynomials'' of knots (e.g., cf. [HKW]).

We define a linear transformation  $\tau$ on ${\cal H}_k$ by
$$\tau(T_{i_1}T_{i_2}\cdots T_{i_{r-1}}T_{i_r})=T_{i_r}
T_{i_{r-1}}\cdots T_{i_2}T_{i_1}.\eqno(2.84)$$ Then $\tau$
satisfies (2.53) and (2.54) by (2.81)-(2.83) with certain choice
of $\{{\cal B}_{i,j}\}$, say $n=1$ and ${\cal B}_{1,1}={\cal H}_k$
.\psp

{\bf Example 2.2}. Suppose that $G$ is a group and $1$ is its
identity element. Let $\mbb{C}[G]$ be the vector space with a
basis $\{\varpi(g)\mid g\in\G\}$, and multiplication:
$$\varpi(g_1)\cdot\varpi(g_2)=\varpi(g_1g_2)\qquad\for\;\;g_1,g_2\in\G.\eqno(2.85)$$
Then $\mbb{C}[G]$ forms an associative algebra with the identity
element $\varpi(1)$, which is called the {\it group algebra} of
$G$. Define the ``trace map" of $\mbb{C}[G]$ by
$$\tr \varpi(g)=\dlt_{1,g}\qquad\for\;\;g\in\G\eqno(2.86)$$
and a linear transformation $\tau$ on $\mbb{C}[G]$ by
$$\tau(\varpi(g))=\varpi(g^{-1})\qquad\for\;\;g\in\G.\eqno(2.87)$$
It can be verified that $\tau$ satisfies (2.53) and (2.54) by
(2.85)-(2.87) with certain choice of $\{{\cal B}_{i,j}\}$, say
$n=1$ and ${\cal B}_{1,1}=\mbb{C}[G]$.\psp

The  representations of $\hat{\cal A}_{\vec\ell}$ and $\hat{\cal
A}_{\vec\ell}^\tau$ with ${\cal A}={\cal H}_k,\mbb{C}[G]$ and
their related vertex algebras will be studied in our future works.
In the rest of this paper, we will deal with the case when ${\cal
A}$ is a matrix algebra.

\section{Modules Related to General Linear Algebras}

In this section, we give detailed constructions of irreducible
modules of the Lie algebras $\hat{\cal A}_{\vec\ell}$ in (2.48)
and $\hat{\cal A}_{\vec\ell}^\tau$ in (2.59) when ${\cal A}$ is
the $n\times n$ matrix algebra, from weighted irreducible modules
of a central extension of  the general linear  Lie algebra of
infinite matrices with finite number of nonzero entries.

Recall that $\ol{gl}(\infty)$ be a vector space with a basis
$\{\E_{l,k}\mid l,k\in\Z\}$ and the multiplication:
$$ \E_{l_1,l_2}\cdot
\E_{k_1,k_2}=\dlt_{l_2+k_1,0}\E_{l_1,k_2}\qquad\for
\;\;l_1,l_2,k_1,k_2\in\Z.\eqno(3.1)$$ Moreover, $\ol{gl}(\infty)$
is isomorphic to the associative algebra of infinite matrices with
finite number of nonzero entries.  Let $M_{n\times n}(\mbb{C})$ be
the algebra of $n\times n$ matrices with entries in $\mbb{C}$.
Then we have
$$M_{n\times n}(\mbb{C})\otimes_{\mbb{C}}\ol{gl}(\infty)\cong
\ol{gl}(\infty)\;\;\mbox{as assciative algebras}.\eqno(3.2)$$

Take
$${\cal A}=M_{n\times n}(\mbb{C})\;\;\mbox{in last
section}\eqno(3.3)$$ with the sum of diagonal entries as the trace
map ``$\tr$''. Let $M_{n\times n}(\mbb{A})$ be the algebra of
$n\times n$ matrices with entries in $\mbb{A}$ (cf. (2.1)-(2.3)).
Again $E_{i,j}$  is the $n\times n$ matrix with 1 as its
$(i,j)$-entry and 0 as the others. Recall the Lie algebra
$$\widehat{gl}(n,\mbb{A})=M_{n\times
n}(\mbb{A})\oplus\mbb{C}\kappa\eqno(3.3)$$ with the Lie bracket:
\begin{eqnarray*} & &[t^{m_1}\ptl_t^{r_1}E_{i_1,j_1}+\mu_1\kappa,t^{m_2}\ptl_t^{r_2}
E_{i_2,j_2}+\mu_1\kappa]\\
&=&\dlt_{j_1,i_2}t^{m_1}\ptl_t^{r_1}\cdot
t^{m_2}\ptl_t^{r_2}E_{i_1,j_2}-\dlt_{i_1,j_2}t^{m_2}\ptl_t^{r_2}\cdot
t^{m_1}\ptl_t^{r_1}E_{i_2,j_1}\\ & &+(-1)^{r_1}\dlt_{i_1,j_2}\dlt_{j_1,i_2}
\dlt_{r_1+r_2,m_1+m_2}r_1!r_2!\left(\!\!\begin{array}{c}m_1\\
r_1+r_2+1\end{array}
\!\!\right)\kappa\hspace{5.3cm}(3.4)\end{eqnarray*} for
$i,j\in\ol{1,n},\; m_1,m_2\in\mbb{Z},\;r_1,r_2\in\mbb{N}$ and
$\mu_1,\mu_2\in\mbb{C}$. The Lie algebras
$$\hat{\cal A}\cong \widehat{gl}(n,\mbb{A}).\eqno(3.5)$$
Fixed an element $\vec\ell\in\mbb{N}^{\:n}$. We have the Lie
subalgebra
$$\widehat{gl}(\vec\ell,\mbb{A})=\sum_{i,j=1}^n\mbb{A}\ptl_t^{\ell_j}E_{i,j}+\mbb{C}\kappa
\eqno(3.6)$$ of $\widehat{gl}(n,\mbb{A})$ of type $\hat{\cal
A}_{\vec\ell}$ (cf. (2.48)).

Let $H$ be a toral Cartan subalgebra of a Lie algebra ${\cal G}$.
We always denote
$$H^\ast=\mbox{the space of linear functions
on}\;\;H.\eqno(3.7)$$ A $\cal G$-module ${\cal M}$ is called {\it
weighted} if
$${\cal M}=\bigoplus_{\nu\in\T^\ast}{\cal M}_\nu,\qquad {\cal
M}_\nu=\{u\in{\cal M}\mid h(u)=\nu(h)u,\;h\in H\}.\eqno(3.8)$$

Recall the algebra
$$\td{gl}(\infty)=\ol{gl}(\infty)\oplus\mbb{C}\kappa_0\eqno(3.9)$$
with the Lie bracket:\begin{eqnarray*}\hspace{1cm}&
&[\E_{l_1,l_2}+\mu_1\kappa_0,\E_{k_1,k_2}+\mu_2\kappa_0]
=\E_{l_1,l_2}\E_{k_1,k_2}- \E_{k_1,k_2}\E_{l_1,l_2}
\\ & &+\dlt_{l_1+k_2,0}\dlt_{l_2+k_1,0}[H(l_1)H(l_2)-H(k_1)H(k_2)]\kappa_0\hspace{5.6cm}
(3.10)\end{eqnarray*} for $l_1,l_2,k_1,k_2\in\Z$ and
$\mu_1,\mu_2\in\mbb{C}$ (cf. (2.23)), where $\kappa_0$ is a base
element. By (2.21)-(2.23) and (3.2), the linear map
$$\kappa_0\leftrightarrow\kappa_0,\;\;\E_{i,j}(l+1/2,k-1/2)\leftrightarrow
\E_{ln+i-1/2,kn-j+1/2}
\qquad\for\;\;i,j\in\ol{1,n},\;l,k\in\mbb{Z},\eqno(3.11)$$ gives a
Lie algebra isomorphism between the Lie algebra $\check{\cal A}$
and $ \td{gl}(\infty)$ under our assumption (3.3). Moreover, the
subspace
$$\T=\sum_{l\in\Z}\mbb{C}\E_{l,-l}+\mbb{C}\kappa_0\eqno(3.12)$$
forms a toral Cartan subalgebra of $\td{gl}(\infty)$. In fact, the
root structures of $\ol{gl}(\infty)$ and $\td{gl}(\infty)$ are the
same. Thus they have exactly the same representation theory. We
can take
$$\{\E_{l+1,-l}\mid l\in\Z\}\;\;\mbox{as positive simple root
vectors}\eqno(3.13)$$ and
$$\{\E_{l,-l-1}\mid l\in\Z\}\;\;\mbox{as negative simple root
vectors}.\eqno(3.14)$$ Furthermore,
$$
[\E_{l+1,-l},\E_{l,-l-1}]=\E_{l+1,-l-1}-\E_{l,-l}\qquad\for\;\;-\frac{1}{2}\neq
l\in\Z\eqno(3.15)$$ and
$$[\E_{1/2,1/2},\E_{-1/2,-1/2}]=\E_{1/2,-1/2}-\E_{-1/2,1/2}+\kappa_0.\eqno(3.16)$$

Set
$$\td{gl}^m(\infty)=\mbox{Span}\:\{\E_{l_1,-l_2},\E_{-l_1,l_2}\mid
m< l_1,l_2\in\Z\}\qquad\for\;\;m\in\mbb{N}.\eqno(3.17)$$ Then
(3.1) and (3.10) show that $\{\td{gl}^m(\infty)\mid m\in\mbb{N}\}$
are Lie subalgebras of $\td{gl}(\infty)$. Let ${\cal M}$ be a
weighted $\td{gl}(\infty)$-module
$$\mbox{generated by a subspace}\;\;{\cal M}_0\;\;\mbox{such
that}\;\;\td{gl}^m(\infty)({\cal M}_0)=\{0\}\;\;\mbox{for
some}\;\;m\in\mbb{N}.\eqno(3.18)$$ Fixed a constant
$\iota\in\mbb{C}$. By Theorem 2.2 and (3.2), we have the following
representation $\sgm_{\cal M}^\iota$ of
$\widehat{gl}(\vec\ell,\mbb{A})$: $\sgm_{\cal
M}^\iota(\kappa)=\kappa_0$ and $$\sgm_{\cal
M}^\iota(t^m\ptl_t^rE_{i,j})
=\sum_{l\in\mbb{Z}}\la\iota-l\ra_r\E_{(m-r-l)n+i-1/2,ln-j+1/2}+r!\dlt_{i,j}\dlt_{r,m}\Im_{0,r}
\kappa_0\eqno(3.19)$$
 for $ m\in\mbb{Z},\;r\in\mbb{N}$ and
$i,j\in\ol{1,n}$. \psp

{\bf Theorem 3.1} {\it Suppose $\iota\not\in\mbb{Z}$. Then the
representation $\sgm_{\cal M}^\iota$ of
$\widehat{gl}(\vec\ell,\mbb{A})$ is irreducible if and only if
$\cal M$ is an irreducible $\td{gl}(\infty)$-module.}

{\it Proof}. Denote $$h_{i,r}=\sgm_{\cal
M}^\iota(t^r\ptl_t^rE_{i,i})=\sum_{l\in\mbb{Z}}\la l+\iota\ra_
r\E_{ln+i-1/2,-ln-i+1/2}+\Im_{0,r}\kappa_0\eqno(3.20)$$ for
$i\in\ol{1,n}$ and $r\in\mbb{N}+\ell_i$. Set
$$H=\sum_{i=1}^n\:\sum_{r=0}^\infty\mbb{C}h_{i,r}\subset\mbox{End}\:{\cal M},\eqno(3.21)$$
the space of linear transformations on $\cal M$. As operators on
$\cal M$,
$$[h_{i,r},\E_{ln+j_1-1/2,kn-j_2+1/2}]=[\dlt_{i,j_1}\la
l+\iota\ra_r-\dlt_{i,j_2}\la-k+\iota\ra_r\E_{ln+j_1-1/2,kn-j_2+1/2}\eqno(3.22)$$
for $j_1,j_2\in\ol{1,n}$ and $l,k\in\mbb{Z}$. Using generating
functions, we have:
\begin{eqnarray*}& &[\sum_{r=0}^{\infty}\frac{x^r}{r!}h_{i,r+\ell_i},\E_{ln+j_1-1/2,kn-j_2+1/2}]
\\ &=&\left[\dlt_{i,j_1}\sum_{r=0}^{\infty}\frac{\la
l+\iota\ra_{r+\ell_i}x^r}{r!}-\dlt_{i,j_2}\sum_{r=0}^{\infty}\frac{\la
-k+\iota\ra_{r+\ell_i}x^r}{r!}\right]\E_{ln+j_1-1/2,kn-j_2+1/2}\\
&=&\left[\dlt_{i,j_1}\la
l+\iota\ra_{\ell_i}\sum_{r=0}^{\infty}\frac{\la
l+\iota-\ell_i\ra_rx^r}{r!}-\dlt_{i,j_2}\la
-k+\iota\ra_{\ell_i}\sum_{r=0}^{\infty}\frac{\la
-k+\iota-\ell_i\ra_rx^r}{r!}\right]\\ & &\times\E_{ln+j_1-1/2,kn-j_2+1/2}\\
&=& [\dlt_{i,j_1}\la l+\iota\ra_{\ell_i}(x+1)^{
l+\iota-\ell_i}-\dlt_{i,j_2}\la-k+\iota\ra_{\ell_i}(x+1)^{-k+\iota-\ell_i}]
\E_{ln+j_1-1/2,kn-j_2+1/2}
\\ &=&\frac{d^{\ell_i}}{dx^{\ell_i}}[\dlt_{i,j_1}(x+1)^{
l+\iota}-\dlt_{i,j_2}(x+1)^{-k+\iota}]
\E_{ln+j_1-1/2,kn-j_2+1/2}.\hspace{4.3cm}(3.23)\end{eqnarray*}
Note that when $\iota\not\in\mbb{Z}$,
$$\frac{d^{\ell_i}}{dx^{\ell_i}}(x+1)^{m_1+\iota}=
\frac{d^{\ell_i}}{dx^{\ell_i}}(x+1)^{m_2+\iota}\;\for\;m_1,m_2\in\mbb{Z}\llra
m_1=m_2\eqno(3.24)$$ and
$$\frac{d^{\ell_i}}{dx^{\ell_i}}[(x+1)^{
l_1+\iota}-(x+1)^{-k_1+\iota}]=\frac{d^{\ell_i}}{dx^{\ell_i}}[(x+1)^{
l_2+\iota}-(x+1)^{-k_2+\iota}]\eqno(3.25)$$
$$\for\;l_1,l_2,k_1,k_2\in\mbb{Z},\;l_1\neq -k_1\llra
l_1=l_2,\;k_1=k_2.\eqno(3.26)$$

 Denote by $H^\ast$ the space of linear functions on $H$. Given
$\rho\in H^\ast$, we set
$$\td{gl}(\infty)_{(\rho)}=\{\xi\in \td{gl}(\infty)\mid
[h,\xi]=\rho(h)\xi\;\for\;h\in H\}\eqno(3.27)$$ and $${\cal
M}_{(\rho)}=\{w\in{\cal M}\mid h(w)=\rho(h)w\;\for\;h\in
H\}.\eqno(3.28)$$ Then
$$\td{gl}(\infty)_{(0)}=\T,\;\;\dim\td{gl}(\infty)_{(\rho)}=1\qquad\for\;\;0\neq\rho
\in H^\ast\eqno(3.29)$$ by (3.23)-(3.26). Moreover,
$$\td{gl}(\infty)=\bigoplus_{\rho\in H^\ast}\td{gl}(\infty)_{(\rho)}.\eqno(3.30)$$
Since ${\cal M}$ is a weighted $\td{gl}(\infty)$-module, we have
$${\cal M}=\bigoplus_{\rho\in H^\ast}{\cal
M}_{(\rho)}\eqno(3.31)$$ by (3.12), (3.20) and (3.21). If $V$ is a
$\td{gl}(\infty)$-submodule of $\cal M$, then $V$ is a
$\widehat{gl}(\vec\ell,\mbb{A})$-submodule of $\cal M$ by (3.19).
Suppose that $U$ is a $\widehat{gl}(\vec\ell,\mbb{A})$-submodule
of $\cal M$. Then
$$U=\bigoplus_{\rho\in H^\ast}U_{(\rho)},\qquad
U_{(\rho)}=U\bigcap{\cal M}_{(\rho)}.\eqno(3.32)$$ Since
\begin{eqnarray*}\sgm_{\cal
M}^\iota(t^m\ptl_t^rE_{i,j})(U_{(\rho)})
&=&\sum_{l\in\mbb{Z}}\la\iota-l\ra_r\E_{(m-r-l)n+i-1/2,ln-j+1/2}(U_{(\rho)})\\
& &+r!\dlt_{i,j}\dlt_{r,m}\Im_{0,r}\kappa_0(U_{(\rho)})\subset U,
\hspace{6cm}(3.33)\end{eqnarray*}we have
$$\E_{(m-r-l)n+i-1/2,ln-j+1/2}(U_{(\rho)})\subset U\eqno(3.34)$$
for $m\in\mbb{Z},\;r\in\mbb{N},\; i,j\in\ol{1,n}$ such that $
(m,i)\neq (0,j)$ by (3.29). Observe that $\{\E_{l_1,l_2}\mid
l_1,l_2\in\Z;\;l_1\neq -l_2\}$ generates the Lie algebra
$\td{gl}(\infty)$. Thus (3.33) implies
$$\td{gl}(\infty)U_{(\rho)})\subset U\qquad\for\;\;\rho\in
H^\ast.\eqno(3.35)$$ Therefore, $U$ is a
$\td{gl}(\infty)$-submodule.$\qquad\Box$ \psp

The above theorem shows that we construct a family of irreducible
representations $\{\sgm_{\cal M}^\iota\mid
\iota\in\mbb{C}\setminus\mbb{Z}\}$ from any irreducible weighted
$\td{gl}(\infty)$-module ${\cal M}$ satisfying (3.18). Set
$$E_{i,j}(r,z)=\sum_{m\in\mbb{Z}}t^m\ptl_t^rE_{i,j}z^{-m-1}\qquad\for\;\;i,j\in\ol{1,n},
\;r\in\mbb{N}.\eqno(3.36)$$ By a similar proof as that of Theorem
3.1, we have the following result, which was proved by Ma [M2] in
more complicated form when $\iota=0$ . \psp

 {\bf Theorem 3.3}.
{\it Suppose $\iota\in\mbb{Z}$. Let ${\cal M}$ be a weighted
$\td{gl}(\infty)$-module satisfying (3.18). We have the following
representation $\sgm_{\vec\ell,{\cal M}}$ of
$\widehat{gl}(\vec\ell,\mbb{A})$: $\sgm_{\vec\ell,{\cal
M}}(\kappa)=\kappa_0$ and
\begin{eqnarray*} & & \sgm_{\vec\ell,{\cal M}}(E_{i,j}(r,z))\\&=&
\sum_{l,k=0}^n[\la
-k-1\ra_r\E_{(l-\iota)n+i-1/2,(k+\iota+1)n-j+1/2}z^{-l-k-\ell_i-r-2}\\
& &+ \la k+\ell_j\ra_
r\E_{(l-\iota)n+i-1/2,(\iota-k)n-j+/2}z^{-l+k+\ell_j-\ell_i-r-1}\\
& &+\la k+\ell_j \ra_r\E_{-(l+\iota+1)n+i-1/2,(\iota-k)n-j+1/2}z^{l+k+\ell_j-r}\\
& & +\la-k-1\ra_r
\E_{-(l+\iota+1)n+i-1/2,(k+\iota+1)n-j+1/2}z^{l-k-r-1}]+\dlt_{i,j}r!\Im_{0,r}\kappa_0z^{-r-1}.
\hspace{1.8cm}(3.37)\end{eqnarray*} for $i,j\in\ol{1,n}$ and
$r\in\mbb{N}+\ell_j$. Moreover, $\sgm_{\vec\ell,{\cal M}}$ if and
only if $\cal M$ is irreducible. } \psp

Denote
$$i^\ast=n+1-i\qquad\for\;\;i\in\ol{1,n}.\eqno(3.38)$$
We fix $\es\in\{0,1\}$ and take
$$\vec\ell=(\ell_1,\ell_2,...,\ell_n)\in\mbb{N}^{\:n}\;\;\mbox{such
that}\;\;\{\ell_1,\ell_2,...,\ell_n\}\subset
2\mbb{N}+\es\eqno(3.39)$$ and
$$\ell_i=\ell_{i^\ast}\qquad\for\;\;i\in\ol{1,n}.\eqno(3.40)$$
\newpage

 For any $$\sum_{i,j=1}^n\mu_{i,j}E_{i,j}\in{\cal A},\eqno(3.41)$$
we define
$$(\sum_{i,j=1}^n\mu_{i,j}E_{i,j})^\ast=\sum_{i,j=1}^n\mu_{i,j}E_{j^\ast,i^\ast}.\eqno(3.42)$$
Then $\ast$ is an involution of ${\cal A}$ (cf. (2.53), (2.54)).
Now the Lie algebra $\hat{\cal A}_{\vec\ell}^\ast$ (cf. (2.59))
becomes
$$\hat{o}(\vec\ell,\mbb{A})=\check{A}_{\vec\ell}^\ast=\sum_{i,j=1}^n\:\sum_{r=0}^{\infty}
\sum_{m\in\mbb{Z}}
\mbb{C}(t^m\ptl_t^{r+\ell_j}E_{i,j}-(-1)^\es(-\ptl)^rt^m\ptl_t^{\ell_i}E_{j^\ast,i^\ast})
+ \mbb{C}\kappa.\eqno(3.43)$$ Again let $\cal M$ be a weighted
$\td{gl}(\infty)$-module satisfying (3.18). Then we have the
restrict representation of $\sgm_{\cal M}^\iota$ on
$\hat{o}(\vec\ell,\mbb{A})$ with $\sgm_{\cal
M}^\iota(\kappa)=\kappa_0$ and
\begin{eqnarray*}& &\sgm_{\cal
M}^\iota(t^{m+r}\ptl_t^{r+\ell_j}E_{i,j}-(-1)^\es(-\ptl)^rt^{m+r}
\ptl_t^{\ell_i}E_{j^\ast,i^\ast})\\
&=&\sum_{l\in\mbb{Z}}[\la
l+\ell_j+\iota\ra_{r+\ell_j}\E_{(m+l)n+i-1/2,(-l-\ell_j)n-j+1/2}-(-1)^\es\la
l-\iota\ra_r\\\ & &\times
\la-m-l+\ell_i+\iota-1\ra_{\ell_i}\E_{-ln-j+1/2,(m+l-\ell_i)n+i-1/2}]\\
&&+[(r+\ell_i)!\Im_{0,r+\ell_i}-(-1)^\es
r!\ell_i!\Im_{r,\ell_i}]\dlt_{m,\ell_i}\dlt_{i,j}\kappa_0\hspace{6.2cm}(3.44)\end{eqnarray*}
for $i,j\in\ol{1,n},\;m\in\mbb{Z}$ and $r\in\mbb{N}$.\psp

{\bf Theorem 3.3} {\it Suppose $\iota\not\in\mbb{Z}/2$. Then the
representation $\sgm_{\cal M}^\iota$ of
$\hat{o}(\vec\ell,\mbb{A})$ is irreducible if and $\cal M$ is an
irreducible $\td{gl}(\infty)$-module.}

{\it Proof}. For $i\in\ol{1,n}$ and $r\in\mbb{N}$, we define
\begin{eqnarray*}& &\eta_{i,r}=\sgm_{\cal M}^\iota(t^{r+\ell_i}\ptl_t^{r+\ell_i}E_{i,i}-
(-1)^\es(-\ptl)^rt^{r+\ell_i} \ptl_t^{\ell_i}E_{i^\ast,i^\ast})\\
&=&\sum_{l\in\mbb{Z}}(\la
l+\iota\ra_{r+\ell_i}\E_{ln+i-1/2,-ln-i+1/2}- \la
\ell_i-l-\iota\ra_{r+\ell_i}\E_{ln-i+1/2,-ln+i-1/2})
\\ & &+[(r+\ell_i)!\Im_{0,r+\ell_i}-(-1)^\es
r!\ell_i!\Im_{r,\ell_i}]\kappa_0.\hspace{7.6cm}(3.45)\end{eqnarray*}

Set
$$H_o=\sum_{i=1}^n\:\sum_{r=0}^\infty\mbb{C}\eta_{i,r}\subset\mbox{End}\:{\cal M},\eqno(3.46)$$
the space of linear transformations on $\cal M$. As operators on
$\cal M$,
\begin{eqnarray*} &
&[\eta_{i,r},\E_{ln+j_1-1/2,kn-j_2+1/2}]= [\dlt_{i,j_1}\la
l+\iota\ra_{r+\ell_i}- \dlt_{i^\ast,j_1}\la
\ell_i-l-1-\iota\ra_{r+\ell_i}
\\
& &-\dlt_{i,j_2}\la \iota-k\ra_{r+\ell_i}+\dlt_{i^\ast,j_2} \la
k+\ell_i-1-\iota\ra_{r+\ell_i}]\E_{ln+j_1-1/2,kn-j_2+1/2}\hspace{3.9cm}(3.47)\end{eqnarray*}
for $j_1,j_2\in\ol{1,n}$ and $l,k\in\mbb{Z}$. Using generating
functions, we have:
\begin{eqnarray*}& &[\sum_{r=0}^{\infty}\frac{x^r}{r!}\eta_{i,r},\E_{ln+j_1-1/2,kn-j_2+1/2}]
\\ &=&\frac{d^{\ell_i}}{dx^{\ell_i}}[\dlt_{i,j_1}(x+1)^{
l+\iota}-\dlt_{i^\ast,j_1}(x+1)^{\ell_i-l-1-\iota})-\dlt_{i,j_2}(x+1)^{\iota-k}\\
& & +\dlt_{i^\ast,j_2} (x+1)^{k+\ell_i-1-\iota})]
\E_{ln+j_1-1/2,kn-j_2+1/2}.\hspace{6.9cm}(3.48)\end{eqnarray*}

Since $\iota\not\in\mbb{Z}/2$, we have
$$l+\iota\neq
k-\iota\qquad\for\;\;l,k\in\mbb{Z},\;i\in\ol{1,n}.\eqno(3.49)$$
Thus
$$(x+1)^{l_1+\iota}-(x+1)^{k_1-\iota}=(x+1)^{l_2+\iota}-(x+1)^{k_2-\iota}
\;\;\for\;\;l_1,l_2,k_1,k_2\in\mbb{Z}\llra
l_1=l_2,\;k_1=k_2\eqno(3.50)$$ and
\begin{eqnarray*} \qquad& &(x+1)^{
l_1+\iota}-(x+1)^{\ell_i-l_1-1-\iota}-(x+1)^{\iota-k_1}+(x+1)^{k_1+\ell_i-1-\iota}
\\ &=&(x+1)^{
l_2+\iota}-(x+1)^{\ell_i-l_2-1-\iota}-(x+1)^{\iota-k_2}+(x+1)^{k_2+\ell_i-1-\iota}
\hspace{2.8cm}(3.51)\end{eqnarray*} for
$l_1,k_1,l_2,k_2\in\mbb{Z}$ with $l_1\neq -k_1,\;l_2\neq -k_2$ if
and only if $l_1=l_2$ and $k_1=k_2$. Denote by $H_o^\ast$ the
space of linear functions on $H_o$. Given $\rho\in H_o^\ast$, we
set
$$\td{gl}(\infty)_{[\rho]}=\{\xi\in \td{gl}(\infty)\mid
[h,\xi]=\rho(h)\xi\;\for\;h\in H_o\}.\eqno(3.52)$$ Then
$$\td{gl}(\infty)_{[0]}=\T,\;\;\dim\td{gl}(\infty)_{[\rho]}=1\qquad\for\;\;0\neq\rho
\in H_o^\ast\eqno(3.52)$$ by the above arguments. Therefore, the
conclusion follows the same arguments as those in
(3.29)-(3.35).$\qquad\Box$ \psp

Next we suppose that
$$n=2n_0\;\; \mbox{is an even positive integer}. \eqno(3.54)$$
Moreover, we define the parity of indices:
$$p(i)=0,\qquad
p(n_0+i)=1\qquad\for\;\;i\in\ol{1,n_0}.\eqno(3.55)$$

For any element in (3.41), we define
$$(\sum_{i,j=1}^n\mu_{i,j}E_{i,j})^{\dg}
=\sum_{i,j=1}^n(-1)^{p(i)+p(j)}\mu_{i,j}E_{j^\ast,i^\ast}
\eqno(3.56)$$ Then $\dg$ is an involution of $\cal A$ (cf. (2.53),
(2.54)). Now the Lie algebra $\hat{\cal A}_{\vec\ell}^\dg$ (cf.
(2.59)) becomes
$$\widehat{sp}(\vec\ell,\mbb{A})=\sum_{i,j=1}^n\:\sum_{r=0}^{\infty}
\sum_{m\in\mbb{Z}}\mbb{C}(t^m\ptl_t^{r+\ell_j}E_{i,j}-(-1)^{p(i)+p(j)+\es}(-\ptl_t)^rt^m
\ptl_t^{\ell_i}E_{j^\ast,i^\ast}) + \mbb{C}\kappa.\eqno(3.57)$$
 Let $\cal M$ be  a weighted
$\td{gl}(\infty)$-module satisfying (3.18). Then we have the
restricted representation of $\sgm_{\cal M}^\iota$ on
$\widehat{sp}(\vec\ell,\mbb{A})$ with $\sgm_{\cal
M}^\iota(\kappa)=\kappa_0$ and
\begin{eqnarray*}& &\sgm_{\cal M}^\iota(t^{m+r}\ptl_t^{r+\ell_j}E_{i,j}-(-1)^{p(i)+p(j)
+\es}(-\ptl)^rt^{m+r}
\ptl_t^{\ell_i}E_{j^\ast,i^\ast})\\
&=&\sum_{l\in\mbb{Z}}[\la
l+\ell_j+\iota\ra_{r+\ell_j}\E_{(m+l)n+i-1/2,(-l-\ell_j)n-j+1/2}-(-1)^{p(i)+p(j)+\es}\la
l-\iota\ra_r\\\ & &\times
\la-m-l+\ell_i+\iota-1\ra_{\ell_i}\E_{-ln-j+1/2,(m+l-\ell_i)n+i-1/2}]\\
&&+[(r+\ell_i)!\Im_{0,r+\ell_i}-(-1)^{p(i)+p(j)+\es}
r!\ell_i!\Im_{r,\ell_i}]\dlt_{i,j}\dlt_{m,\ell_i}\kappa_0\hspace{4.7cm}(3.58)\end{eqnarray*}
for $i,j\in\ol{1,n},\;m\in\mbb{Z}$ and $r\in\mbb{N}$.\psp

{\bf Theorem 3.4}. {\it Suppose $\iota\not\in\mbb{Z}/2$. Then the
representation $\sgm_{\cal M}^\iota$ of
$\widehat{sp}(\vec\ell,\mbb{A})$ is irreducible if and $\cal M$ is
an irreducible $\td{gl}(\infty)$-module.} \psp

 {\bf Example 3.1}.
Let $\lmd\in\T^\ast$ (cf. (3.12)) such that there exists a
positive integer $m_0$ for which
$$\lmd(\E_{l,-l})=\lmd(\E_{-l,l})=0\qquad\for\;\;m_0<l\in\Z.\eqno(3.59)$$
Moreover, we set
$$\td{gl}(\infty)_+=\mbox{Span}\:\{\E_{l,k}\mid
l,k\in\Z;\;l+k>0\},\;\;\td{gl}(\infty)_-=\mbox{Span}\:\{\E_{l,k}\mid
l,k\in\Z;\;l+k<0\}\eqno(3.60)$$ and
$$\td{gl}(\infty)_0=\td{gl}(\infty)_++\T\eqno(3.61)$$
(cf. (3.12)). Then $\td{gl}(\infty)_{\pm}$ and $\td{gl}(\infty)_0$
are Lie subalgebras of $\td{gl}(\infty)$. Define a one-dimensional
$\td{gl}(\infty)_0$-module $\mbb{C}v_\lmd$ by
$$\td{gl}(\infty)_+(v_\lmd)=\{0\},\;\;h(v_\lmd)=\lmd(h)v_\lmd\qquad\for\;\;h\in\T.\eqno(3.62)$$

Form an induced $\td{gl}(\infty)$-module
$$M_\lmd=U(\td{gl}(\infty))\otimes_{U(\td{gl}(\infty)_0)}\mbb{C}v_\lmd\cong
U(\td{gl}(\infty)_-)\otimes_{\mbb{C}}\mbb{C}v_\lmd ,\eqno(3.63)$$
that is, a Verma module. It is known that $M_\lmd$ has a unique
maximal proper submodule $N_\lmd$. Thus
$${\cal M}=M_\lmd/N_\lmd\eqno(3.64)$$
is a weighted irreducible  $\td{gl}(\infty)$-module  satisfying
(3.18). Suppose
$$\lmd_l=\lmd(\E_{l+1,-l-1}-\E_{l,-l}),\;\lmd_{-1/2}=\lmd(\E_{1/2,-1/2}-\E_{-1/2,1/2}+\kappa_0)
\in\mbb{N}\eqno(3.65)$$ for $-1/2\neq l\in \Z$. Then
$$N_\lmd=\sum_{r\in\Z}U(\td{gl}(\infty)_-)\E_{r,-r-1}^{\lmd_r+1}\otimes
v_\lmd.\eqno(3.66)$$ \psp

{\bf Example 3.2}. Let $m$ be a fixed positive integer. Set
$$I=\{i-1/2,-i+1/2\mid i\in\ol{1,m}\},\qquad J=\Z\setminus
I.\eqno(2.67)$$ The subspaces
$${\cal G}=\mbox{Span}\:\{\E_{l,k},\kappa_0\mid l,k\in
I\}\;\;\mbox{and}\;\;\ol{\cal
G}=\mbox{Span}\:\{\E_{l,k},\kappa_0\mid l,k\in J\} \eqno(3.68)$$
forms  Lie subalgebras of $\td{gl}(\infty)$. The algebra ${\cal
G}$  is isomorphic a one-dimensional central extension of
$gl(2m,\mbb{C})$ with
$$H=\sum_{l\in I}\mbb{C}\E_{l,-l}+\mbb{C}\kappa_0\eqno(3.69)$$
as a Cartan subalgebra. Moreover, the algebra $\ol{\cal G}$ is
isomorphic $\td{gl}(\infty)$ with
$$\ol{H}=\sum_{l\in J}\mbb{C}\E_{l,-l}+\mbb{C}\kappa_0\eqno(3.70)$$
as a Cartan subalgebra.

Denote
$${\cal L}_0={\cal G}+\ol{\cal G},\;\;{\cal L}_-=\mbox{Span}\:\{\E_{l,r}\mid l\in J,\;r\in I\},
\;\; {\cal L}_+=\mbox{Span}\:\{\E_{r,l}\mid r\in I,\;l\in
J\}.\eqno(3.71)$$ Then ${\cal L}_0$ and ${\cal L}_{\pm}$ are Lie
subalgebras of $\td{gl}(\infty)$. In fact,
$$\td{gl}(\infty)={\cal L}_-\oplus {\cal L}_0\oplus {\cal
L}_+.\eqno(3.72)$$ Note $$\ol{\cal G}_\pm=\ol{\cal G}\bigcap
\td{gl}(\infty)_{\pm}\eqno(3.73)$$ (cf. (3.60)) are Lie
subalgebras of $\ol{\cal G}$ and
$$\ol{\cal G}=\ol{\cal G}_-\oplus\ol{H}\oplus \ol{\cal
G}_-.\eqno(3.74)$$ In particular, we have a Borel subalgebra
$$\ol{\cal G}_0=\ol{H}+\ol{\cal G}_+.\eqno(3.75)$$

Take any weighted irreducible $\cal G$-module $M_0$, which may not
necessarily be of highest weight type. We extend the action of
$\kappa_0$ to that of $\ol{\cal G}_0$ by
$$\ol{\cal G}_+(M_0)=\E_{k,-k}(M_0)=\{0\}\qquad\for\;\;k\in
J.\eqno(3.76)$$ Form an induced $\ol{\cal G}$-module
$$M_1=U(\ol{\cal G})\otimes_{U(\ol{\cal G}_0)}M_0\cong
U(\ol{\cal G}_-)\otimes_{\mbb{C}}M_0.\eqno(3.77)$$ Since
$$[{\cal G},\ol{\cal G}]=\{0\},\eqno(3.78)$$
$M_1$ becomes an ${\cal L}_0$-module by letting ${\cal G}$ act on
the second factor. Moreover, $M_1$ has a unique maximal proper
${\cal L}_0$-submodule $M_2$.  Form a quotient ${\cal L}_0$-module
$${\cal M}_0=M_1/M_2.\eqno(3.79)$$
In fact, (3.76) yields
$$\E_{l,-k}(M_0+M_2)=\E_{-l,k}(M_0+M_2)\subset
M_2\qquad\for\;\;m<l,k\in\Z.\eqno(3.80)$$

Note that
$$[{\cal L}_0,{\cal L}_\pm]\subset {\cal L}_\pm.\eqno(3.81)$$
So
$$ {\cal L}'={\cal L}_0+{\cal L}_+\eqno(3.82)$$
form a Lie subalgebra of $\td{gl}(\infty)$. We extend an action of
${\cal L}'$ on ${\cal M}_0$ from that of ${\cal L}_0$ by
$${\cal L}_+({\cal M}_0)=\{0\}.\eqno(3.83)$$
The expression (3.81) implies that ${\cal M}_0$ becomes an ${\cal
L}'$-module. Form an induced $\td{gl}(\infty)$-module:
$${\cal M}_1=U(\td{gl}(\infty))\otimes_{U({\cal L}')}{\cal
M}_0\cong U({\cal L}_-)\otimes_{\mbb{C}}{\cal M}_0.\eqno(3.84)$$
It can be verified that ${\cal M}_1$ has a unique maximal proper
$\td{gl}(\infty)$-submodule ${\cal M}_2$. The  quotient
$${\cal M}={\cal M}_1/{\cal M}_2\eqno(3.85)$$
is a weighted irreducible $\td{gl}(\infty)$-module satisfying
(3.18) by (3.80).

\section{Modules with $\iota\in\mbb{Z}+1/2$ Related to Skew Elements}

In this section, we give detailed constructions of irreducible
modules of the Lie algebras $\hat{\cal A}_{\vec\ell}^\tau$ in
(2.59) with $\iota\in\mbb{Z}+1/2$ when ${\cal A}$ is the $n\times
n$ matrix algebra, from weighted irreducible modules of central
extensions of  the Lie algebras  of infinite skew matrices with
finite number of nonzero entries.

Recall the Lie algebra $\ol{gl}(\infty)$ defined in (3.1). The
subspaces
$$\bar{o}_d(\infty)=\sum_{k,l\in\Z}\mbb{C}(\E_{l,k}-\E_{k,l}),\eqno(4.1)$$
$$\bar{o}_b(\infty)=\sum_{k,l\in\Z}\mbb{C}(\E_{l,k}-\E_{k-1,l+1}),\eqno(4.2)$$
and
$$\ol{sp}(\infty)=\sum_{k,l\in\Z;\;kl<0}\mbb{C}(\E_{k,l}-\E_{l,k})+\sum_{k,l\in\Z;\;kl>0}
\mbb{C}(\E_{k,l}+\E_{l,k})\eqno(4.3)$$ forms Lie subalgebras skew
elements in $\ol{gl}(\infty)$. Take
$$\{\E_{l+1,-l}-\E_{-l,l+1},\E_{3/2,1/2}-\E_{1/2,3/2}\mid
0<l\in\Z\}\eqno(4.4)$$ as positive simple root vectors of
$\bar{o}_d(\infty)$ and
$$\{\E_{l,-l-1}-\E_{-l-1,l},\E_{-1/2,-3/2}-\E_{-3/2,-1/2}\mid
0<l\in\Z\}\eqno(4.5)$$ as negative simple root vectors of
$\bar{o}_d(\infty)$. Choose
$$\{\E_{l,1-l}-\E_{-l,l+1}\mid
0<l\in\Z\}\eqno(4.6)$$ as positive simple root vectors of
$\bar{o}_b(\infty)$ and
$$\{\E_{l-1,-l}-\E_{-l-1,l}\mid
0<l\in\Z\}\eqno(4.7)$$ as negative simple root vectors of
$\bar{o}_b(\infty)$. Pick
$$\{\E_{l+1,-l}-\E_{-l,l+1},\E_{1/2,1/2}\mid
0<l\in\Z\}\eqno(4.8)$$ as positive simple root vectors of
$\ol{sp}(\infty)$ and
$$\{\E_{l,-l-1}-\E_{-l-1,l},\E_{-1/2,-1/2}\mid
0<l\in\Z\}\eqno(4.9)$$ as negative simple root vectors of
$\ol{sp}(\infty)$.

Again we assume ${\cal A}=M_{n\times n}(\mbb{C})$ in the
   settings of Section 2.  Recall the Lie algebra $\check{\cal A}^{\tau,\iota}_{\vec\ell}$
defined (2.69), the assumption (3.39), the involution $\ast$
defined in (3.42) and the involution $\dg$ defined in (3.56). \psp

{\bf Theorem 4.1}.  {\it We have the following Lie algebra
isomorphisms:
$$\check{\cal A}^{\ast,\iota}_{\vec\ell}/\mbb{C}\kappa_0\cong
\left\{\begin{array}{ll}\bar{o}_b(\infty)&\mbox{\it
if}\;\;\es=0,\;n\in 2\mbb{N}+1,\\
\bar{o}_d(\infty)&\mbox{if}\;\;\es=0,\;n\in 2\mbb{N},\\
\ol{sp}(\infty)&\mbox{if}\;\;\es=1
\end{array}\right.\eqno(4.10)$$
and
$$\check{\cal A}^{\dg,\iota}_{\vec\ell}/\mbb{C}\kappa_0\cong
\left\{\begin{array}{ll}\ol{sp}(\infty)&\mbox{\it
if}\;\;\es=0,\\
\bar{o}_d(\infty)&\mbox{if}\;\;\es=1
\end{array}\right.\eqno(4.11)$$
if $n$ is even.}

{\it Proof}. We write
$$\iota=\iota_0+1/2,\qquad\iota_0\in\mbb{Z}.\eqno(4.12)$$ By assumption (3.39) and (3.40),
$$\ell_i=\ell_{i^\ast}=2m_i+\es\;\;\mbox{with}\;\;m_i\in\mbb{N}\qquad\for\;\;i\in\ol{1,n}.
\eqno(4.13)$$ Thus (2.64) and (2.69) give \begin{eqnarray*}
\check{\cal
A}^{\ast,\iota}_{\vec\ell}&=&\sum_{i,j=1}^n\:\sum_{l,k\in\Z}\mbb{C}((-1)^\es\la
k+m_j+\es-1\ra_{\ell_j}E_{i,j}(l+m_i-\iota_0,k+\iota_0-m_j)\\
& &-\la
l+m_i\ra_{\ell_i}E_{j^\ast,i^\ast}(k+m_j-\iota_0+\es-1,l-m_i+\iota_0+1-\es))+\mbb{C}\kappa_0,
\hspace{2.3cm}(4.14)
\end{eqnarray*}
and if $n=2n_0$ is even,  \begin{eqnarray*} \check{\cal
A}^{\dg,\iota}_{\vec\ell}&=&\sum_{i,j=1}^n\:\sum_{l,k\in\Z}\mbb{C}((-1)^\es\la
k+m_j+\es-1\ra_{\ell_j}E_{i,j}(l+m_i-\iota_0,k+\iota_0-m_j)-(-1)^{p(i)+p(j)}\\
& &\times\la l+m_i\ra_{\ell_i}
E_{j^\ast,i^\ast}(k+m_j-\iota_0+\es-1,l+\iota_0-m_i+1-\es))+\mbb{C}\kappa_0
\hspace{2.5cm}(4.15)
\end{eqnarray*}
(cf. (3.55)). According to (3.11), we set \begin{eqnarray*}&
&{\cal G}^{\ast,\vec m
}_\es=\sum_{i,j=1}^n\:\sum_{l,k\in\mbb{Z}}\mbb{C}((-1)^\es\la
k+m_j+\es-3/2\ra_{\ell_j}\E_{(l+m_i-\iota_0)n+i-1/2,(k-m_j+\iota_0)n-j+1/2}\\
& &-\la
l+m_i+1/2\ra_{\ell_i}\E_{(k+m_j-\iota_0+\es-1)n-j+1/2,(l-m_i+\iota_0+1-\es)n+i-1/2})
+\mbb{C}\kappa_0, \hspace{2.7cm}(4.16)
\end{eqnarray*}
and
\begin{eqnarray*}& &{\cal
G}^{\dg,\vec
m}_\es=\sum_{i,j=1}^n\:\sum_{l,k\in\mbb{Z}}\mbb{C}((-1)^\es\la
k+m_j+\es-3/2\ra_{\ell_j}\E_{(l+m_i-\iota_0)n+i-1/2,(k-m_j+\iota_0)n-j+1/2}\\
& &-(-1)^{p(i)+p(j)}\la
l+m_i+1/2\ra_{\ell_i}\E_{(k+m_j-\iota_0+\es-1)n-j+1/2,(l-m_i+\iota_0+1-\es)n+i-1/2})
+\mbb{C}\kappa_0 \hspace{0.9cm}(4.17)
\end{eqnarray*}
if $n$ is even.
 Then ${\cal
G}^{\ast,\vec m}_\es$ and ${\cal G}^{\dg,\vec m}_\es$ are Lie
subalgebras of $\td{gl}(\infty)$. Moreover,
$$\check{\cal A}^{\ast,\iota}_{\vec\ell}\cong {\cal
G}^{\ast,\vec m}_\es,\;\;\check{\cal
A}^{\dg,\iota}_{\vec\ell}\cong {\cal G}^{\dg,\vec
m}_\es\eqno(4.18)$$ by (3.11).

 Define
\begin{eqnarray*}\bar{\cal G}^{\ast,\vec m}_\es&=&\sum_{i,j=1}^n\:\sum_{l,k\in\mbb{Z}}\mbb{C}(
(-1)^\es\la k+m_j+\es-3/2\ra_{\ell_j}\E_{ln+i-1/2,kn-j+1/2}\\ &
&-\la
l+m_i+1/2\ra_{\ell_i}\E_{(k+\es-1)n-j+1/2,(l+1-\es)n+i-1/2})\hspace{5.2cm}(4.19)
\end{eqnarray*}
and
\begin{eqnarray*}\bar{\cal G}^{\dg,\vec m}_\es&=&\sum_{i,j=1}^n\:\sum_{l,k\in\mbb{Z}}\mbb{C}
((-1)^\es\la
k+m_j+\es-3/2\ra_{\ell_j}\E_{ln+i-1/2,kn-j+1/2}\\ &
&-(-1)^{p(i)+p(j)}\la
l+m_i+1/2\ra_{\ell_i}\E_{(k+\es-1)n-j+1/2,(l+1-\es)n+i-1/2})\hspace{3.2cm}(4.20)
\end{eqnarray*}
if $n$ is even. Then $\bar{\cal G}^{\ast,\vec m}_\es$ and
$\bar{\cal G}^{\dg,\vec m}_\es$ are Lie subalgebras of
$\ol{gl}(\infty)$. The conclusion follows from (4.18) and the
following facts:
$${\cal G}^{\ast,\vec m}_0/\mbb{C}\kappa_0\cong\bar{\cal G}^{\ast,\vec m}_0\cong
\bar{o}_b(\infty)\qquad\mbox{if}\;n\;\mbox{is odd},\eqno(4.21)$$
$${\cal G}^{\ast,\vec m}_0/\mbb{C}\kappa_0\cong\bar{\cal G}^{\ast,\vec m}_0\cong
\bar{o}_d(\infty)\qquad\mbox{if}\;n\;\mbox{is even},\eqno(4.22)$$
$${\cal G}^{\ast,\vec m}_1/\mbb{C}\kappa_0\cong\bar{\cal G}^{\ast,\vec m}_1
\cong\ol{sp}(\infty),\qquad{\cal G}^{\dg,\vec m
}_0/\mbb{C}\kappa_0\cong\bar{\cal G}^{\dg,\vec m
}_0\cong\ol{sp}(\infty),\eqno(4.23)$$
$${\cal G}^{\dg,\vec m}_1/\mbb{C}\kappa_0\cong\bar{\cal G}^{\dg,\vec m}_1\cong\bar{o}_d(\infty).
\qquad\Box\eqno(4.24)$$ \vspace{0.1cm}

Next we want to study highest weight irreducible  modules.
Isomorphisms in (4.21)-(4.24) motivate us to adjust the definition
of the Lie algebra $\td{gl}(\infty)$ (cf. (3.9) and (3.10)) by
modifying the coefficient of $\kappa_0$ in (3.10). Moreover, for
any $l\in\mbb{Z}$, we write
$$l=l_Qn+l_R,\qquad l_Q\in\mbb{Z},\;l_R\in\ol{1,n}.\eqno(4.25)$$
Let $\vec m=(m_1,m_2,...,m_n)\in\mbb{N}^{\:n}$ and define a map
$\al_{\vec m}^{\iota_0}:\Z^{\:4}\rta\mbb{C}$ by:
\begin{eqnarray*}& &\al^{\iota_0}_{\vec
m}(l_1,l_2;k_1,k_2)=(H(l_1+(m_{(l_1+1/2)_R}-\iota_0)n)
H(l_2+(\iota_0-m_{(-l_2+1/2)_R})n)\\
&&-H(k_1+(m_{(k_1+1/2)_R}-\iota_0)n)H(k_2+(\iota_0-m_{(-k_2+1/2)_R})n)\dlt_{l_1+k_2,0}
\dlt_{l_2+k_1,0}\hspace{2.2cm} (4.26)\end{eqnarray*} for
$l_1,l_2,k_1,k_2\in\Z$ (cf. (1.15)). Set
$$\td{gl}^{\iota_0}_{\vec m}(\infty)=\ol{gl}(\infty)\oplus\mbb{C}\kappa_0\eqno(4.27)$$ (cf.
(3.1)), where $\kappa_0$ is a base element. We have the following
Lie bracket on $\td{gl}^{\iota_0}_{\vec m}(\infty)$:
$$[\E_{l_1,l_2}+\mu_1\kappa_0,\E_{k_1,k_2}+\mu_2\kappa_0]=\E_{l_1,l_2}\E_{k_1,k_2}-
\E_{k_1,k_2}\E_{l_1,l_2}+\al^{\iota_0}_{\vec
m}(l_1,l_2;k_1,k_2)\kappa_0\eqno(4.28)$$ for
$l_1,l_2,k_1,k_2\in\Z$. In particular, the Lie algebras
$$\td{gl}^0_{\vec 0}(\infty)\cong \td{gl}(\infty).\eqno(4.29)$$

Assume that (4.13) holds. Now we define
\begin{eqnarray*}{\cal L}^{\ast,\vec
m}_{\iota_0,\es}&=&\sum_{i,j=1}^n\:\sum_{l,k\in\mbb{Z}}\mbb{C}((-1)^\es\la
k+m_j+\es-3/2\ra_{\ell_j}\E_{ln+i-1/2,kn-j+1/2}\\ & &-\la
l+m_i+1/2\ra_{\ell_i}\E_{(k+\es-1)n-j+1/2,(l+1-\es)n+i-1/2})+\mbb{C}\kappa_0,
\hspace{3.9cm}(4.30)\end{eqnarray*} and \begin{eqnarray*}{\cal
L}^{\dg,\vec m
}_{\iota_0,\es}&=&\sum_{i,j=1}^n\:\sum_{l,k\in\mbb{Z}}\mbb{C}((-1)^\es\la
k+m_j+\es-3/2\ra_{\ell_j}\E_{ln+i-1/2,kn-j+1/2}-(-1)^{p(i)+p(j)}\\
& &\times\la
l+m_i+1/2\ra_{\ell_i}\E_{(k+\es-1)n-j+1/2,(l+1-\es)n+i-1/2})+\mbb{C}\kappa_0.
\hspace{4cm}(4.31)\end{eqnarray*} if $n$ is even. Then ${\cal
L}^{\ast,\vec m}_{\iota_0,\es}$ and ${\cal L}^{\dg,\vec m
}_{\iota_0,\es}$ are Lie subalgebras of $\td{gl}^{\iota_0}_{\vec
m}(\infty)$. Moreover,
$$\check{\cal A}^{\ast,\iota}_{\vec \ell}\cong {\cal
G}^{\ast,\vec m}_\es\cong {\cal L}^{\ast,\vec
m}_{\iota_0,\es},\;\;\check{\cal A}^{\dg,\iota}_{\vec\ell}\cong
{\cal G}^{\dg,\vec m}_\es\cong {\cal L}^{\dg,\vec
m}_{\iota_0,\es}.\eqno(4.32)$$

Write
$$n=2n_0+\ves,\qquad n_0\in\mbb{N},\;\ves\in\{0,1\}.\eqno(4.33)$$
 We denote
\begin{eqnarray*} e^\ast_{\es,l}&=&\frac{1}{\la
m_{l_R}-l_Q+\es-3/2\ra_{\ell_{l_R}}\la
m_{(l+1)_R}+(l+1)_Q+1/2\ra_{\ell_{(l+1)_R}}}\\ & &\times[(-1)^\es
\la
m_{l_R}-l_Q+\es-3/2\ra_{\ell_{l_R}}\E_{l+1/2,-l+1/2}\\
& & -\la
m_{(l+1)_R}+(l+1)_Q+1/2\ra_{\ell_{(l+1)_R}}\E_{-l+(\es-1)n+1/2,l+(1-\es)n+1/2}]
\hspace{3.3cm}(4.34)\end{eqnarray*} and
\begin{eqnarray*} f^\ast_{\es,l}&=&(-1)^\es\la
m_{(l+1)_R}-(l+1)_Q+\es-3/2\ra_{\ell_{(l+1)_R}}\E_{l-1/2,-l-1/2}\\
& & -\la
m_{l_R}+l_Q+1/2\ra_{\ell_{l_R}}\E_{-l+(\es-1)n-1/2,l+(1-\es)n-1/2}
\hspace{5.5cm}(4.35)\end{eqnarray*} for
$l\in\mbb{N}+\dlt_{\es,0}(\dlt_{\ves,0}-n_0)+\dlt_{\es,1}$. When
$\es=\ves=0$, we define
\begin{eqnarray*}e^\ast_{0,-n_0}&=&\frac{1}{\la m_{n_0}-1/2\ra_{\ell_{n_0}}\la
m_{(n_0+2)_R}+(n_0+2)_Q-1/2\ra_{\ell_{(n_0+2)_R}}}\\ & &\times[
\la m_{n_0}-1/2\ra_{\ell_{n_0}}\E_{-n_0+3/2, n_0+1/2}
\\ & &-\la
m_{(n_0+2)_R}+(n_0+2)_Q-1/2\ra_{\ell_{(n_0+2)_R}}\E_{-n_0+1/2,n_0+3/2}]
\hspace{3.9cm}(4.36)\end{eqnarray*} and
\begin{eqnarray*}f^\ast_{0,-n_0}&=&
\la m_{(n_0+2)_R}-(n_0+2)_Q-1/2\ra_{\ell_{(n_0+2)_R}}\E_{-n_0-1/2,
n_0-3/2} \\ & &-\la
m_{n_0}-1/2\ra_{\ell_{n_0}}\E_{-n_0-3/2,n_0-1/2}.
\hspace{7.5cm}(4.37)\end{eqnarray*} Furthermore, we let
$$ e^\ast_{1,0}=\E_{1/2,1/2},\qquad
f^\ast_{1,0}=\E_{-1/2,-1/2}.\eqno(4.38)$$

Suppose $\ves=0$. We define
\begin{eqnarray*} e^\dg_{\es,l}&=&\frac{1}{\la
m_{l_R}-l_Q+\es-3/2\ra_{\ell_{l_R}}\la
m_{(l+1)_R}+(l+1)_Q+1/2\ra_{\ell_{(l+1)_R}}}\\ & &\times[(-1)^\es
\la m_{l_R}-l_Q+\es-3/2\ra_{\ell_{l_R}}\E_{l+1/2,-l+1/2}
-(-1)^{p((l+1)_R)+p(l_R)}\\
& &\times\la
m_{(l+1)_R}+(l+1)_Q+1/2\ra_{\ell_{(l+1)_R}}\E_{-l+(\es-1)n+1/2,l+(1-\es)n+1/2}],
\hspace{3.1cm}(4.39)\end{eqnarray*} and
\begin{eqnarray*} f^\dg_{\es,l}&=&(-1)^\es\la
m_{(l+1)_R}-(l+1)_Q+\es-3/2\ra_{\ell_{(l+1)_R}}\E_{l-1/2,-l-1/2}\\
& & -(-1)^{p((l+1)_R)+p(l_R)}\la
m_{l_R}+l_Q+1/2\ra_{\ell_{l_R}}\E_{-l+(\es-1)n-1/2,l+(1-\es)n-1/2}
\hspace{2.4cm}(4.40)\end{eqnarray*} for
$l\in\mbb{N}+1-n_0\dlt_{\es,0}$. Moreover, we set
$$e^\dg_{0,-n_0}=\E_{-n_0+1/2,n_0+1/2},\qquad
f^\dg_{0,0}=\E_{-n_0-1/2,n_0-1/2}.\eqno(4.41)$$ Furthermore, we
let
\begin{eqnarray*}\qquad\qquad e^\dg_{1,,0}&=&\frac{1}{\la m_1+1/2\ra_{\ell_1}\la
m_{2_R}+2_Q+1/2\ra_{\ell_{2_R}}}\\ & &\times[ \la
m_1+1/2\ra_{\ell_1}\E_{3/2, 1/2} +(-1)^{p(n)+p(2)}\\ & &\times\la
m_{2_R}+2_Q+1/2\ra_{\ell_{2_R}}\E_{1/2,3/2}]
\hspace{6.8cm}(4.42)\end{eqnarray*} and
\begin{eqnarray*}\qquad \qquad f^\dg_{1,,0}&=&\la
m_{2_R}-2_Q-1/2\ra_{\ell_{2_R}}\E_{-1/2, -3/2}
\\ & &+(-1)^{p(n)+p(2)}\la
m_1+3/2\ra_{\ell_1}\E_{-3/2,-1/2}.\hspace{5.5cm}(4.43)\end{eqnarray*}

For convenience, we always assume
$$\tau\in\{\ast,\dg\}.\eqno(4.44)$$
Under the above settings,
$$\{ e^\tau_{\es,l}\mid l\in \mbb{N}-\dlt_{0,\es}n_0\}\;\;\mbox{is set of
 positive simple root
vectors of}\;\;{\cal L}^{\tau,\vec\ell}_{\iota_0,\es}\eqno(4.45)$$
 and
$$\{ f^\tau_{\es,l}\mid l\in \mbb{N}-\dlt_{0,\es}n_0\}\;\;
\mbox{is set of negative simple root vectors of}\;\;{\cal L
}^{\tau,\vec m}_{\iota,\es}.\eqno(4.46)$$

Set
$$\vt^\es_l=\E_{l-1/2,-l+1/2}-\E_{-l+(\es-1)n+1/2,l+(1-\es)n-1/2}\eqno(4.47)$$
for $l\in\mbb{N}-n_0\dlt_{\es,0}$. Define
\begin{eqnarray*}\omega^\tau_{\es,l}&=&\al^{\iota_0}_{\vec
m}(l+1/2,-l+1/2;l-1/2,-l-1/2)+\al^{\iota_0}_{\vec m}(-l+(\es-1)n+1/2,\\
&
&l+(1-\es)n+1/2;-l+(\es-1)n-1/2,l+(1-\es)n-1/2)\hspace{3.7cm}(4.48)\end{eqnarray*}
for $l\in\mbb{N}+\dlt_{\es,0}(\dlt_{\ves,0}-n_0)+\dlt_{\es,1}$ if
$\tau=\ast$, and $l\in\mbb{N}+1-n_0\dlt_{\es,0}$ when $\tau=\dg$.
Moreover, we let
\begin{eqnarray*}\qquad\omega_{0,-n_0}^\ast&=&\al^{\iota_0}_{\vec m}(-n_0+3/2,
n_0+1/2;-n_0-1/2,n_0-3/2)\\ & &+\al^{\iota_0}_{\vec m
}(-n_0+1/2,n_0+3/2;-n_0-3/2,n_0-1/2)\hspace{4.2cm}(4.49)\end{eqnarray*}
when $\es=\ves=0$,
$$\omega_{1,0}^\ast=\al^{\iota_0}_{\vec
m}(1/2,1/2;-1/2,-1/2),\eqno(2.50)$$
$$\omega_{0,-n_0}^\dg=\al^{\iota_0}_{\vec
m}(-n_0+1/2,n_0+1/2;-n_0-1/2,n_0-1/2)\eqno(4.51)$$ and
$$\omega^\dg_{1,0}=\al^{\iota_0}_{\vec m}(3/2,1/2;-1/2,-3/2)+\al^{\iota_0}_{\vec m
}(1/2,3/2;-3/2,-1/2).\eqno(4.52)$$

Set
$$T^\tau_{\es,l}=\vt_{l+1}^\es-\vt_l^\es+\omega^\tau_{\es,l}\kappa_0\eqno(4.53)$$
for $l\in\mbb{N}+\dlt_{\es,0}(\dlt_{\ves,0}-n_0)+\dlt_{\es,1}$ if
$\tau=\ast$, and $l\in\mbb{N}+1-n_0\dlt_{\es,0}$ when $\tau=\dg$.
Moreover, we set
$$T_{0,-n_0}^\ast=\vt^0_{-n_0}+\vt^0_{-n_0+2}+\omega_{0,-n_0}^\ast\kappa_0\eqno(4.54)$$
when $\es=\ves=0$,
$$T_{1,0}^\ast=\vt^1_1+\omega_{1,0}^\ast\kappa_0,\;\;T^\dg_{0,-n_0}=\vt^0_{-n_0+1}
+\omega_{0,-n_0}^\dg\kappa_0\eqno(4.55)$$ and
$$T^\dg_{1,0}=\vt^1_1+\vt^1_2+\omega^\dg_{1,0}\kappa_0.\eqno(4.56)$$

It can be verified that
$$[e^\tau_{\es,l},f^\tau_{\es,l}]=T^\tau_{\es,l}\qquad\for\;\;l\in\mbb{N}-\dlt_{0,\es}n_0.
\eqno(4.57)$$ Obviously
$$[T^\tau_{\es,l},e^\tau_{\es,l}]=2e^\tau_{\es,l},\;\;[T^\tau_{\es,l},f^\tau_{\es,l}]=
-2f^\tau_{\es,l}\qquad\for\;\; l\in\mbb{N}-\dlt_{0,\es}n_0.
\eqno(4.58)$$
\newpage

 Set
$$\T^\es=\sum_{l=-n_0\dlt_{\es,0}}^{\infty}\mbb{C}\vt_l^\es+\mbb{C}\kappa_0.\eqno(4.59)$$
 Moreover, we let
\begin{eqnarray*}{\cal L}^{\ast,\vec m}_{\iota_0,\es,\pm}&=&\sum_{i,j=1}^n\:
\sum_{l,k\in\mbb{Z};\pm(l+k)>0}\mbb{C}((-1)^\es\la
k+m_j+\es-3/2\ra_{\ell_j}\E_{ln+i-1/2,kn-j+1/2}\\ & &-\la
l+m_i+1/2\ra_{\ell_i}\E_{(k+\es-1)n-j+1/2,(l+1-\es)n+i-1/2})\hspace{4.9cm}(4.60)\end{eqnarray*}
and
\begin{eqnarray*}{\cal L}^{\dg,\vec m}_{\iota_0,\es,\pm}&=&\sum_{i,j=1}^n\:\sum_{l,k\in\mbb{Z};
\pm(l+k)>0}\mbb{C}((-1)^\es\la
k+m_j+\es-3/2\ra_{\ell_j}\E_{ln+i-1/2,kn-j+1/2}-(-1)^{p(i)+p(j)}\\
& &\times\la
l+m_i+1/2\ra_{\ell_i}\E_{(k+\es-1)n-j+1/2,(l+1-\es)n+i-1/2}).\hspace{4.8cm}(4.61)\end{eqnarray*}
Then ${\cal L}^{\tau,\vec m}_{\iota_0,\es,\pm}$ are Lie
subalgebras of ${\cal L}^{\tau,\vec m}_{\iota_0,\es}$, $\T^\es$ is
a toral Cartan subalgebra of ${\cal L}^{\tau,\vec
m}_{\iota_0,\es}$ and
$${\cal L}^{\tau,\vec m}_{\iota_0,\es}={\cal L}^{\tau,\vec
m}_{\iota_0,\es,-}\oplus\T^\es\oplus{\cal L}^{\tau,\vec
m}_{\iota_0,\es,+}.\eqno(4.62)$$ Furthermore, we have Borel
subalgebras:
$${\cal L}^{\tau,\vec m}_{\iota_0,\es,0}=\T^\es+{\cal L}^{\tau,\vec
m}_{\iota_0,\es,+}\eqno(4.63)$$ of ${\cal L}^{\tau,\vec
m}_{\iota_0,\es}.$

Denote by $(\T^\es)^\ast$ the space of linear functions on
$\T^\es$. Fix an element
$$k_0\in\mbb{N}-n_0\dlt_{0,\es}.\eqno(4.64)$$
Take $\lmd\in (\T^\es)^\ast$ such that
$$\lmd(\vt^\es_l)=0\qquad\for\;\;k_0\leq
l\in\mbb{N}-n_0\dlt_{\es,0}.\eqno(4.65)$$ Define a one-dimensional
${\cal L}^{\tau,\vec m}_{\iota_0,\es,0}$-module
$\mbb{C}v^{\tau,\es}_\lmd$ by
$${\cal L}^{\tau,\vec
m}_{\iota_0,\es,+}(v^{\tau,\es}_\lmd)=\{0\},\;\;h(v^{\tau,\es}_\lmd)=\lmd(h)v^{\tau,\es}_\lmd
\qquad\for\;\;h\in\T^\es.\eqno(4.66)$$ Form an induced ${\cal
L}^{\tau,\vec m}_{\iota,\es}$-module:
$$M^{\tau,\es}_\lmd=U({\cal L}^{\tau,\vec m}_{\iota_0,\es})\otimes_{U({\cal L}^{\tau,\vec
m}_{\iota_0,\es,0})}\mbb{C}v^{\tau,\es}_\lmd\cong U({\cal
L}^{\tau,\vec
m}_{\iota_0,\es,-})\otimes_{\mbb{C}}\mbb{C}v^{\tau,\es}_\lmd.\eqno(4.67)$$
There exists a unique maximal proper submodule $N^{\tau,\es}_\lmd$
of $M^{\tau,\es}_\lmd$ , and the quotient module
$${\cal M}^{\tau,\es}_\lmd=M^{\tau,\es}_\lmd/N^{\tau,\es}_\lmd\eqno(4.68)$$
is a weighted irreducible ${\cal L}^{\tau,\vec
m}_{\iota_0,\es}$-module. Identify $1\otimes v^{\tau,\es}_\lmd$
with $v^{\tau,\es}_\lmd$. If
$$\lmd^\es_l=\lmd(T^\tau_{\es,l})\in\mbb{N}\qquad\for\;\;l\in\mbb{N}-n_0\dlt_{\es,0},\eqno(4.69)$$
then
$$N^{\tau,\es}_\lmd=\sum_{l=-n_0\dlt_{0,\es}}^{\infty}U({\cal L}^{\tau,\vec
m}_{\iota_0,\es,-})(f^\tau_{\es,l})^{\lmd_l^\es+1}v^{\tau,\es}_\lmd.\eqno(4.70)$$

Set
$$\hat\ell=|\iota_0|+1+\mbox{max}\:\{m_1,m_2,...,m_n\}.\eqno(4.71)$$
For $\hat{\ell}<s\in\mbb{N}$, we define
\begin{eqnarray*}{\cal L}^{\ast,\vec
m,s}_{\iota_0,\es}&=&\sum_{i,j=1}^n\:\sum_{l,-k\in \pm(\mbb{N}+s)}
\mbb{C}((-1)^\es\la k+m_j+\es-3/2\ra_{\ell_j}\E_{ln+i-1/2,kn-j+1/2}\\
& &-\la
l+m_i+1/2\ra_{\ell_i}\E_{(k+\es-1)n-j+1/2,(l+1-\es)n+i-1/2})\hspace{5cm}(4.72)\end{eqnarray*}
and
\begin{eqnarray*}{\cal L}^{\dg,\vec
m,s}_{\iota_0,\es}&=&\sum_{i,j=1}^n\:\sum_{l,-k\in
\pm(\mbb{N}+s)}\mbb{C}((-1)^\es\la
k+m_j+\es-3/2\ra_{\ell_j}\E_{ln+i-1/2,kn-j+1/2}\\ &
&-(-1)^{p(i)+p(j)}\la
l+m_i+1/2\ra_{\ell_i}\E_{(k+\es-1)n-j+1/2,(l+1-\es)n+i-1/2}).
\hspace{2.9cm}(4.73)\end{eqnarray*} Then ${\cal L}^{\tau,\vec
m,s}_{\iota_0,\es}$ is a Lie subalgebra of ${\cal L}^{\tau,\vec
m}_{\iota,\es}$ by (4.16)-(4.26).

Suppose that ${\cal M}$ is an ${\cal L}^{\tau,\vec
m}_{\iota_0,\es}$-module
$$\mbox{generated by a subspace}\;{\cal M}_0\;\mbox{such
that}\;{\cal L}^{\tau,\vec m,s}_{\iota_0,\es}({\cal
M}_0)=\{0\}\;\;\mbox{for
some}\;\hat{\ell}<s\in\mbb{N}.\eqno(4.74)$$ For instance, the
above module ${\cal M}^{\tau,\es}_\lmd$ is such an module with
${\cal M}_0=\mbb{C}v^{\tau,\es}_\lmd$ and
$s=\max\{\hat\ell,k_0+2\}$. We can also construct $\cal M$ as
Example 3.2. According to (2.60), (2.64) and (2.73), we define a
representation $\sgm_\ast$ of $\hat{o}(\vec\ell,\mbb{A})$ and a
representation $\sgm_\dg$ of $\widehat{sp}(\vec\ell,\mbb{A})$ on
${\cal M}$ as follows:
$\sgm_\ast(\kappa)=\kappa_0,\;\sgm_\dg(\kappa)=\kappa_0$,
\begin{eqnarray*}&
&\sgm_\ast(t^{k+m_i+m_j+r+\es}\ptl_t^{r+\ell_j}E_{i,j}-
(-\ptl_t)^rt^{k+m_i+m_j+r+\es}\ptl_t^{\ell_i}E_{j^\ast,i^\ast})\\
&=&\sum_{l\in\mbb{Z}}\la l-m_j-\es+1/2\ra_r((-1)^\es\la
-l+m_j+\es-3/2\ra_{\ell_j}\E_{(k+l)n+i-1/2,-ln-j+1/2}\\ & &-\la
k+l+m_i+1/2\ra_{\ell_i}\E_{(-l+\es-1)n-j+1/2,(k+l+1-\es)n+i-1/2})\\
&&+((r+\ell_i)!\Im_{0,r+\ell_i}-r!\ell_i!\Im_{r,\ell_i})\dlt_{k+\es,0}
\dlt_{i,j}\kappa_0\hspace{7.2cm}(4.75)\end{eqnarray*} and
\begin{eqnarray*}&
&\sgm_\dg(t^{k+m_i+m_j+r+\es}\ptl_t^{r+\ell_j}E_{i,j}-
(-1)^{p(i)+p(j)}(-\ptl_t)^rt^{k+m_i+m_j+r+\es}\ptl_t^{\ell_i}E_{j^\ast,i^\ast})\\
&=&\sum_{l\in\mbb{Z}}\la l-m_j-\es+1/2\ra_r((-1)^\es\la
-l+m_j+\es-3/2\ra_{\ell_j}\E_{(k+l)n+i-1/2,-ln-j+1/2}\\ & &
-(-1)^{p(i)+p(j)}\la
k+l+m_i+1/2\ra_{\ell_i}\E_{(-l+\es-1)n-j+1/2,(k+l+1-\es)n+i-1/2})\\
&&+((r+\ell_i)!\Im_{0,r+\ell_i}-r!\ell_i!\Im_{r,\ell_i})\dlt_{k+\es,0}
\dlt_{i,j}\kappa_0.\hspace{7cm}(4.76)\end{eqnarray*}
\vspace{0.1cm}

{\bf Theorem 4.2}. {\it Suppose that $\cal M$ is a weighted ${\cal
L}^{\tau,\vec m}_{\iota_0,\es}$-module satisfying (4.74). Then the
representation $\sgm_\tau$ is irreducible if and only if $\cal M$
is irreducible.}
\newpage

{\it Proof}. We only prove the statement for $\sgm_\ast$ when
$\{\ell_s=2m_s\mid s\in\ol{1,n}\}$ and $n=2n_0$  are even. The
other cases can be proved similarly. Define
\begin{eqnarray*}h_{i,r}&=&\sgm_\ast(t^{\ell_i+r}\ptl_t^{r+\ell_i}E_{i,i}
-(-\ptl_t)^rt^{\ell_i+r}\ptl_t^{\ell_i}E_{i^\ast,i^\ast})\\
&=&\sum_{l\in\mbb{Z}}\la
l+m_i+1/2\ra_{2m_i+r}(\E_{ln+i-1/2,-ln-i+1/2}
-\E_{-(l+1)n-i+1/2,(l+1)n+i-1/2})\\
&&+((r+\ell_i)!\Im_{0,r+\ell_i}-r!\ell_i!\Im_{r,\ell_i})\kappa_0\hspace{8cm}(4.77)\end{eqnarray*}
for $i\in\ol{1,n} $ and $r\in\mbb{N}$ by (4.75). Set
$$H=\sum_{i=1}^n\:\sum_{r=0}^\infty\mbb{C}h_{i,r}\subset\mbox{End}\:{\cal M},\eqno(4.78)$$
the space of linear transformations on $\cal M$. As operators on
$\cal M$,
\begin{eqnarray*}& &[h_{i,r},\la
k+m_{j_2}-3/2\ra_{\ell_{j_2}}\E_{ln+j_1-1/2,kn-j_2+1/2}\\ & &-\la
l+m_{j_1}+1/2\ra_{\ell_{j_1}}\E_{(k-1)n-j_2+1/2,(l+1)n+j_1-1/2}]
\\ &=&[\dlt_{i,j_1}\la
l+m_i+1/2\ra_{\ell_i+r}-\dlt_{i^\ast,j_1}\la-l+m_i-3/2\ra_{\ell_i+r}+\dlt_{i^\ast,j_2}\la
k+m_i-3/2\ra_{\ell_i+r}\\ & &-\dlt_{i,j_2}\la
-k+m_i+1/2\ra_{\ell_i+r}] \la
k+m_{j_2}-3/2\ra_{\ell_{j_2}}\E_{(l+j_1-1/2,kn-j_2+1/2}\\ & &-\la
l+m_{j_1}+1/2\ra_{\ell_{j_1}}\E_{(k-1)n-j_2+1/2,(l+1)n+j_1-1/2}.
\hspace{6.1cm}(4.79)\end{eqnarray*} Using generating functions, we
get
\begin{eqnarray*}& &[\sum_{r=0}^\infty h_{i,r}x^r,\la
k+m_{j_2}-3/2\ra_{\ell_{j_2}}\E_{ln+j_1-1/2,kn-j_2+1/2}\\ &&-\la
l+m_{j_1}+1/2\ra_{\ell_{j_1}}\E_{(k-1)n-j_2+1/2,(l+1)n+j_1-1/2}]
\\ &=&\frac{d^{\ell_i}}{dx^{\ell_i}}[\dlt_{i,j_1}(1+x)^{l+m_i+1/2}
-\dlt_{i^\ast,j_1}(1+x)^{-l+m_i
-3/2}+\dlt_{i^\ast,j_2}(1+x)^{k+m_i-3/2}\\ & &-\dlt_{i,j_2}(1+x)^{
-k+m_i+1/2}]\la
k+m_{j_2}-3/2\ra_{\ell_{j_2}}\E_{ln+j_1-1/2,kn-j_2+1/2}\\ & &-\la
l+m_{j_1}+1/2\ra_{\ell_{j_1}}\E_{(k-1)n-j_2+1/2,(l+1)n+j_1-1/2}.
\hspace{6.1cm}(4.80)\end{eqnarray*}

Denote by $H^\ast$ the space of linear functions on $H$. Given
$\rho\in H^\ast$, we set
$$({\cal L}^{\ast,\vec m}_{\iota_0,0})_{(\rho)}=\{\xi\in{\cal L}^{\ast,\vec
m}_{\iota_0,0}\mid [h,\xi]=\rho(h)\xi\;\for\;h\in
H\}.\eqno(4.81)$$ Then
$$({\cal L}^{\ast,\vec m}_{\iota_0,0})_{(0)}=\T,\;\;\dim({\cal L}^{\ast,\vec
m}_{\iota_0,0})_{(\rho)}=1\qquad\for\;\;0\neq\rho\in H^\ast
\eqno(4.82)$$ by (4.80). Moreover,
$${\cal L}^{\ast,\vec m}_{\iota_0,0}=\bigoplus_{\rho\in H^\ast}({\cal L}^{\ast,\vec m}_{1,0})_{(\rho)}.
\eqno(4.83)$$ The conclusion can be proved exactly as the proof of
Theorem 3.1.$\qquad\Box$

\section{Modules with $\iota\in\mbb{Z}$ Related to Skew Elements}

In this section, we give detailed constructions of irreducible
modules of the Lie algebras $\hat{\cal A}_{\vec\ell}^\tau$ in
(2.59) with $\iota\in\mbb{Z}$ when ${\cal A}$ is the $n\times n$
matrix algebra, from weighted irreducible modules of central
extensions of  the Lie algebras  of infinite skew matrices with
finite number of nonzero entries.

 Recall the Lie algebra $\check{\cal D}_{\vec\ell}^{\tau,\iota}$ defined
 in (2.76), and the Lie algebra $\td{gl}(\infty)$ defined in (3.9) and (3.10). Note
\begin{eqnarray*}& &\check{\cal
D}_{\vec\ell}^{\ast,\iota}
=\sum_{i,j=1}^n\:\sum_{0<l,k\in\Z}[\mbb{C}((-1)^\es\la-k-1/2\ra_{\ell_j}
E_{i,j}(-l-\iota,-k+\iota-\ell_j)-\la -l-1/2\ra_{\ell_i}\\ &
&\times E_{j^\ast,i^\ast}
(-k-\iota,-l+\iota-\ell_i))+\mbb{C}((-1)^\es\la-k-1/2\ra_{\ell_j}E_{i,j}(l-\iota+\ell_i,-k+\iota-\ell_j)
\\ & &-\la
l+\ell_i-1/2\ra_{\ell_i}E_{j^\ast,i^\ast}(-k-\iota,l+\iota))+\mbb{C}((-1)^\es\la
k+\ell_j-1/2\ra_{\ell_j}E_{i,j}(l-\iota+\ell_i,k+\iota)\\ & &-\la
l+\ell_i-1/2\ra_{\ell_i}E_{j^\ast,i^\ast}(k-\iota+\ell_j,l+\iota))]+\mbb{C}\kappa_0
\hspace{6.7cm}(5.1)
\end{eqnarray*}
and
\begin{eqnarray*}& &\check{\cal
D}_{\vec\ell}^{\dg,\iota}
=\sum_{i,j=1}^n\:\sum_{0<l,k\in\Z}[\mbb{C}((-1)^\es\la-k-1/2\ra_{\ell_j}
E_{i,j}(-l-\iota,-k+\iota-\ell_j)\\
& &-(-1)^{p(i)+p(j)}\la
-l-1/2\ra_{\ell_i}E_{j^\ast,i^\ast}(-k-\iota,-l+\iota-\ell_i))\\ &
& +\mbb{C}((-1)^\es\la-k-1/2\ra_{\ell_j}E_{i,j}
(l-\iota+\ell_i,-k+\iota-\ell_j) \\ & &-(-1)^{p(i)+p(j)}\la
l+\ell_i-1/2\ra_{\ell_i}E_{j^\ast,i^\ast}(-k-\iota,l+\iota))\\ & &
+\mbb{C}((-1)^\es\la k+\ell_j-1/2\ra_{\ell_j}
E_{i,j}(l-\iota+\ell_i,k+\iota) \\ & &-(-1)^{p(i)+p(j)}\la
l+\ell_i-1/2\ra_{\ell_i}E_{j^\ast,i^\ast}(k-\iota+\ell_j,l+\iota))]+\mbb{C}\kappa_0
\hspace{4.2cm}(5.2)\end{eqnarray*} if $n$ is even. By (3.11), we
set
\begin{eqnarray*}{\cal
G}_{\vec\ell}^{\ast,\iota}
&=&\sum_{i,j=1}^n\:\sum_{l,k=0}^\infty[\mbb{C}((-1)^\es\la-k-1\ra_{\ell_j}
\E_{-(l+\iota+1)n+i-1/2,(-k+\iota-\ell_j)n-j+1/2}\\ & &-\la
-l-1\ra_{\ell_i}\E_{
-(k+\iota)n-j+1/2,(-l+\iota-\ell_i-1)n+i-1/2})
\\ &&+\mbb{C}((-1)^\es\la-k-1\ra_{\ell_j}\E_{(l-\iota+\ell_i)n+i-1/2,(-k+\iota-\ell_j)n-j+1/2}\\
& &-\la l+\ell_i\ra_{\ell_i}\E_{-(k+\iota)n-j+1/2,(l+\iota)n+i-1/2})\\
& &+\mbb{C}((-1)^\es\la
k+\ell_j\ra_{\ell_j}E_{(l-\iota+\ell_i)n+i-1/2,(k+\iota+1)n-j+1/2}\\
& & -\la
l+\ell_i\ra_{\ell_i}E_{(k-\iota+\ell_j+1)n-j+1/2,(l+\iota)n+i-1/2})]+\mbb{C}\kappa_0
\hspace{5cm}(5.3)
\end{eqnarray*}
and
\begin{eqnarray*}{\cal
G}_{\vec\ell}^{\dg,\iota}
&=&\sum_{i,j=1}^n\:\sum_{l,k=0}^\infty[\mbb{C}((-1)^\es\la-k-1\ra_{\ell_j}
\E_{-(l+\iota+1)n+i-1/2,(-k+\iota-\ell_j)n-j+1/2}\\ &
&-(-1)^{p(i)+p(j)}\la -l-1\ra_{\ell_i}\E_{
-(k+\iota)n-j+1/2,(-l+\iota-\ell_i-1)n+i-1/2})
\\ &&+\mbb{C}((-1)^\es\la-k-1\ra_{\ell_j}\E_{(l-\iota+\ell_i)n+i-1/2,(-k+\iota-\ell_j)n-j+1/2}\hspace{8cm}
\end{eqnarray*}
\begin{eqnarray*}\hspace{1cm}
& &-(-1)^{p(i)+p(j)}\la
l+\ell_i\ra_{\ell_i}\E_{-(k+\iota)n-j+1/2,(l+\iota)n+i-1/2})\\&
&+\mbb{C}((-1)^\es\la
k+\ell_j\ra_{\ell_j}E_{(l-\iota+\ell_i)n+i-1/2,(k+\iota+1)n-j+1/2}\\
& & -(-1)^{p(i)+p(j)}\la
l+\ell_i\ra_{\ell_i}E_{(k-\iota+\ell_j+1)n-j+1/2,(l+\iota)n+i-1/2})]+\mbb{C}\kappa_0
\hspace{3.2cm}(5.4)
\end{eqnarray*}
Then ${\cal G}_{\vec\ell}^{\ast,\iota}$ and ${\cal
G}_{\vec\ell}^{\dg,\iota}$ are Lie subalgebras of
$\td{gl}(\infty)$, and the map in (3.11) induces
$$\check{\cal D}_{\vec\ell}^{\tau,\iota}\cong{\cal
G}_{\vec\ell}^{\tau,\iota}\eqno(5.5)$$ (cf. (4.44)).

We define two functions $\hat{H}_1,\:\hat{H}_2:\Z\rta\mbb{C}$ by
$$\hat{H}_1(l)=\left\{\begin{array}{ll}1&\mbox{if}\:\:n\iota<l<0\;\mbox{or}\;l>0,
 (\iota-\ell_{(l+1/2)_R})n,\\
0&\mbox{otherwsie}\end{array}\right.\eqno(5.6)$$ and
$$\hat{H}_2(l)=\left\{\begin{array}{ll}1&\mbox{if}\:\:-\iota
n,0<l\;\mbox{or}\;(\ell_{(-l+1/2)_R}-\iota)n<l<0,\\
0&\mbox{otherwsie}\end{array}\right.\eqno(5.7)$$ (cf. (4.25)).
Moreover, we define a map $\be^{\iota}_{\vec
\ell}:\Z^{\:4}\rta\mbb{C}$ by:
 $$\be^{\iota}_{\vec
\ell}(l_1,l_2;k_1,k_2)=(H_1(l_1)H_2(l_2)-H_1(k_1)H_2(k_2))
\dlt_{l_1+k_2,0} \dlt_{l_2+k_1,0}.\eqno(5.8)$$ Set
$$\td{gl}^{(\iota)}_{\vec\ell}(\infty)=\ol{gl}(\infty)\oplus\mbb{C}\kappa_0\eqno(5.9)$$ (cf.
(3.1)), where $\kappa_0$ is a base element. We have the following
Lie bracket on $\td{gl}^{(\iota)}_{\vec\ell}(\infty)$:
$$[\E_{l_1,l_2}+\mu_1\kappa_0,\E_{k_1,k_2}+\mu_2\kappa_0]=\E_{l_1,l_2}\E_{k_1,k_2}-
\E_{k_1,k_2}\E_{l_1,l_2}+\be^{\iota}_{\vec\ell}(l_1,l_2;k_1,k_2)\kappa_0\eqno(5.10)$$
for $l_1,l_2,k_1,k_2\in\Z$. In particular, the Lie algebras
$$\td{gl}^{(0)}_{\vec\ell}(\infty)\cong \td{gl}(\infty).\eqno(5.11)$$

Next we set
\begin{eqnarray*}{\cal L}^{\ast,\vec\ell}_{\iota}&=&\sum_{i,j=1}^n\:\sum_{l,k=0}^{\infty}[\mbb{C}
((-1)^\es\la-k-1\ra_{\ell_j} \E_{(-l-1)n+i-1/2,-kn-j+1/2}-\la
-l-1\ra_{\ell_i}\\ &
&\times\E_{-kn-j+1/2,(-l-1)n+i-1/2})+\mbb{C}((-1)^\es\la-k-1\ra_{\ell_j}\E_{ln+i-1/2,-kn-j+1/2}\\
& & -\la
l+\ell_i\ra_{\ell_i}\E_{-kn-j+1/2,ln+i-1/2})+\mbb{C}((-1)^\es\la
k+\ell_j\ra_{\ell_j}\E_{ln+i-1/2,(k+1)n-j+1/2}\\ & &-\la
l+\ell_i\ra_{\ell_i}\E_{(k+1)n-j+1/2,ln+i-1/2})]+\mbb{C}\kappa_0\hspace{6.6cm}(5.12)
\end{eqnarray*}
and
\begin{eqnarray*}{\cal L}^{\dg,\vec\ell}_\iota&=&\sum_{i,j=1}^n\:\sum_{l,k=0}^{\infty}[\mbb{C}
((-1)^\es\la-k-1\ra_{\ell_j}
\E_{(-l-1)n+i-1/2,-kn-j+1/2}-(-1)^{p(i)+p(j)}\la
-l-1\ra_{\ell_i}\\ & &\times
\E_{-kn-j+1/2,(-l-1)n+i-1/2})+\mbb{C}((-1)^\es\la-k-1\ra_{\ell_j}\E_{ln+i-1/2,-kn-j+1/2}
-(-1)^{p(i)+p(j)}\\
& &\times\la
l+\ell_i\ra_{\ell_i}\E_{-kn-j+1/2,ln+i-1/2})+\mbb{C}((-1)^\es\la
k+\ell_j\ra_{\ell_j}\E_{ln+i-1/2,(k+1)n-j+1/2}\\ &
&-(-1)^{p(i)+p(j)}\la
l+\ell_i\ra_{\ell_i}\E_{(k+1)n-j+1/2,ln+i-1/2})]+\mbb{C}\kappa_0\hspace{4.5cm}(5.13)
\end{eqnarray*}
if $n$ is even. Then ${\cal L}^{\tau,\vec\ell}_{\iota}$ are Lie
subalgebras of  $\td{gl}^{(\iota)}_{\vec\ell}(\infty)$, and
$$\check{\cal D}_{\vec\ell}^{\tau,\iota}\cong {\cal
G}_{\vec\ell}^{\tau,\iota}\cong{\cal
L}^{\tau,\vec\ell}_\iota\eqno(5.14)$$ (cf. (4.44)) by (5.3) and
(5.4). Thus we have: \psp

{\bf Theorem 5.1}. {\it The Lie algebras:
$$\check{\cal D}_{\vec\ell}^{\ast,\iota}/\mbb{C}\kappa_0\cong
\left\{\begin{array}{ll} \bar{o}_d(\infty)&\mbox{\it if}\;\;\es=0,\\
\ol{sp}(\infty)&\mbox{\it
if}\;\;\es=1\end{array}\right.\eqno(5.15)$$ and
$$\check{\cal D}_{\vec\ell}^{\dg,0}/\mbb{C}\kappa_0\cong
\left\{\begin{array}{ll}\ol{sp}(\infty) &\mbox{\it if}\;\;\es=0,\\
\bar{o}_d(\infty)&\mbox{\it
if}\;\;\es=1\end{array}\right.\eqno(5.16)$$
 (cf. (4.1), (4.3)).}
\psp

Next we want to study highest weight irreducible modules of ${\cal
L}^{\tau,\vec\ell}_\iota$. Recall the notions in (4.25). Set
\begin{eqnarray*}{\cal L}^{\ast,\vec\ell}_{\iota,+}&=&\sum_{0<k<l}\mbb{C}
((-1)^\es\la-k_Q-1\ra_{\ell_{k_R}}\E_{l-1/2,-k+1/2}-\la
l_Q+\ell_{l_R}\ra_{\ell_{l_R}}\E_{-k+1/2,l-1/2})\\ &
&+\sum_{l,k=1}^{\infty}\mbb{C}((-1)^\es\la
k_Q+\ell_{k_R}\ra_{\ell_{k_R}}\E_{l-1/2,k-1/2}-\la
l_Q+\ell_R\ra_{\ell_{l_R}}\E_{k-1/2,l-1/2}),\hspace{1.7cm}(5.17)
\end{eqnarray*}
\begin{eqnarray*}{\cal L}^{\ast,\vec\ell}_{\iota,-}&=&\sum_{l,k=1}^{\infty}\mbb{C}
((-1)^\es\la-k_Q-1\ra_{\ell_{k_R}}
\E_{-l+1/2,-k+1/2}-\la -l_Q-1\ra_{\ell_{l_R}}\E_{-k+1/2,-l+1/2})\\
&
&+\sum_{0<l<k}\mbb{C}((-1)^\es\la-k_Q-1\ra_{\ell_{k_R}}\E_{l-1/2,-k+1/2}-\la
l_Q+\ell_{l_R}\ra_{\ell_{l_R}}\E_{-k+1/2,l-1/2})\hspace{1.2cm}(5.18)
\end{eqnarray*}
and
\begin{eqnarray*}{\cal L}^{\dg,\vec\ell}_{\iota,+}&=&\sum_{0<k<l}\mbb{C}
((-1)^\es\la-k_Q-1\ra_{\ell_{k_R}}\E_{l-1/2,-k+1/2}-(-1)^{p(l_R)+p(k_R)}\\
& &\times \la
l_Q+\ell_{l_R}\ra_{\ell_{l_R}}\E_{-k+1/2,l-1/2})+\sum_{l,k=1}^{\infty}\mbb{C}((-1)^\es\la
k_Q+\ell_{k_R}\ra_{\ell_{k_R}}\E_{l-1/2,k-1/2}\\ &
&-(-1)^{p(l_R)+p((k_R)^\ast)}\la
l_Q+\ell_R\ra_{\ell_{l_R}}\E_{k-1/2,l-1/2},\hspace{5.9cm}(5.19)
\end{eqnarray*}
\begin{eqnarray*}{\cal L}^{\dg,\vec\ell}_{\iota,-}&=&\sum_{l,k=1}^{\infty}\mbb{C}
((-1)^\es\la-k_Q-1\ra_{\ell_{k_R}}
\E_{-l+1/2,-k+1/2}-(-1)^{p((l_R)^\ast)+p(k_Q)}\\ & &\times\la
-l_Q-1\ra_{\ell_{l_R}}
\E_{-k+1/2,-l+1/2})+\sum_{0<l<k}\mbb{C}((-1)^\es\la-k_Q-1\ra_{\ell_{k_R}}\E_{l-1/2,-k+1/2}\\
& &-(-1)^{p(l_R)+p(k_R)}\la
l_Q+\ell_{l_R}\ra_{\ell_{l_R}}\E_{-k+1/2,l-1/2}).\hspace{6cm}(5.20)
\end{eqnarray*}
Then ${\cal L}^{\tau,\vec\ell}_{\iota,\pm}$ are Lie subalgebras of
${\cal L}^{\tau,\vec\ell}_\iota$ (cf. (4.44)).

Denote
$$\vt_l=\E_{l+1/2,-l-1/2}-\E_{-l-1/2,l+1/2}\qquad\for\;\;l\in\mbb{N}.\eqno(5.21)$$
Set
$$\T=\sum_{l=0}^{\infty}\vt_l+\mbb{C}\kappa_0.\eqno(5.22)$$
Then $\T$ is a toral Cartan subalgebra of ${\cal
L}^{\tau,\vec\ell}_\iota$. Moreover,
$${\cal L}^{\tau,\vec\ell}_\iota={\cal L}^{\tau,\vec
\ell}_{\iota,-}\oplus \T\oplus {\cal
L}^{\tau,\vec\ell}_{\iota,+}.\eqno(5.23)$$

Define
$$f^\ast_{\es,l}=\la-(l+1)_Q-1\ra_{\ell_{(l+1)_R}}\E_{l-1/2,-l-1/2}-(-1)^\es\la
l_Q+\ell_{l_R}\ra_{\ell_{l_R}}\E_{-l-1/2,l-1/2}\eqno(5.24)$$ and
\begin{eqnarray*}\hspace{1.9cm}f^\dg_{\es,l}&=&\la-(l+1)_Q-1\ra_{\ell_{(l+1)_R}}\E_{l-1/2,-l-1/2}
\\ & &-(-1)^{\es+p(l_R)+p((l+1)_R)}\la
l_Q+\ell_{l_R}\ra_{\ell_{l_R}}\E_{-l-1/2,l-1/2}\hspace{3.8cm}(5.25)
\end{eqnarray*} for
$l\in\mbb{N}+1$. Moreover, we define
$$f^\ast_{0,0}=\la-1\ra_{\ell_1}
\E_{-3/2,-1/2}-\la-2_Q-1\ra_{\ell_{2_R}}\E_{-1/2,-3/2},\eqno(5.26)$$
$$f^\dg_{1,0}=\la-1\ra_{\ell_1}
\E_{-3/2,-1/2}+(-1)^{p(1)+p((2_R)^\ast)}\la-2_Q-1\ra_{\ell_{2_R}}\E_{-1/2,-3/2},\eqno(5.27)$$
$$f^\ast_{1,0}=f^\dg_{0,0}=\E_{-1/2,-1/2}.\eqno(5.28)$$ Then $\{f^\tau_{\es,l}\mid
l\in\mbb{N}\}$ is a set of negative simple root vectors of ${\cal
L}^{\tau,\vec\ell}_\iota$.

For $l\in\mbb{N}+1$, we define
$$\omega_l=\be^{\iota}_{\vec
\ell}(l+1/2,-l+1/2;l-1/2,-l-1/2)+\be^{\iota}_{\vec
\ell}(-l+1/2,l+1/2;-l-1/2,l-1/2).\eqno(5.29)$$ Moreover, we let
$$\omega^\ast_{0,0}=\be^{\iota}_{\vec
\ell}(1/2,3/2;-3/2,-1/2)+\be^{\iota}_{\vec
\ell}(3/2,1/2;-1/2,-3/2)\eqno(5.30)$$ and
$$\omega^\ast_{1,0}=\be^{\iota}_{\vec
\ell}(1/2,1/2;-1/2,-1/2).\eqno(5.31)$$ Furthermore, we set
$$T^\tau_{\es,l}=\vt_l-\vt_{l-1}+\omega_l\kappa_0\qquad\for\;\;l\in\mbb{N}+1,
\eqno(5.32)$$
$$T^\ast_{0,0}=T^\dg_{1,0}=\vt_1+\vt_0+\omega^\ast_{0,0}\kappa_0\eqno(5.33)$$
and
$$T^\ast_{1,0}=T^\dg_{0,0}=\vt_0+\omega^\ast_{1,0}\kappa_0.\eqno(5.34)$$

Denote
$${\cal L}^{\tau,\vec\ell}_{\iota,0}=\T+ {\cal L}^{\tau,\vec\ell}_{\iota,+}.\eqno(5.35)$$ Let $\lmd^{\tau,\es}$ a linear function
on $\T$ such that there exists $k_0\in\mbb{N}$ for which
$$\lmd^{\tau,\es}(\vt_l)=0\qquad\for\;\;k_0\leq l\in\mbb{N}.\eqno(5.36)$$
Define a one-dimensional ${\cal
L}^{\tau,\vec\ell}_{\iota,0}$-module $\mbb{C}v^{\tau,\es}$ by:
$${\cal L}^{\tau,\vec\ell}_{\iota,+}(v^{\tau,\es})=\{0\},\;\;h
(v^{\tau,\es})=\lmd^{\tau,\es}(h)v^{\tau,\es}\qquad\for\;\;h\in\T.\eqno(5.37)$$
Form an induced ${\cal L}^{\tau,\vec\ell}_\iota$-module:
$$M_{\lmd^{\tau,\es}}=U({\cal L}^{\tau,\vec\ell}_\iota)\otimes_{U({\cal L}^{\tau,\vec\ell}_{\iota,0})}
\mbb{C}v^{\tau,\es}\cong U({\cal L}^{\tau,\vec
\ell}_{\iota,-})\otimes_{\mbb{C}}\mbb{C}v^{\tau,\es}.\eqno(5.38)$$
There exists a unique maximal proper submodule
$N_{\lmd^{\tau,\es}}$ of  $M_{\lmd^{\tau,\es}}$,
 and the quotient
$${\cal
M}_{\lmd^{\tau,\es}}=M_{\lmd^{\tau,\es}}/N_{\lmd^{\tau,\es}}\eqno(5.39)$$
is a weighted irreducible ${\cal L}^{\tau,\vec\ell}_\iota$-module.
Identify $1\otimes v^{\tau,\es}$ with $v^{\tau,\es}$. When
$$\lmd^{\tau,\es}_k=\lmd^{\tau,\es}(T^\tau_{\es,k})
\in\mbb{N}\qquad\for\;\; k\in\mbb{N},\eqno(5.40)$$ the submodule
$$N_{\lmd^{\tau,\es}}=\sum_{l=0}^{\infty}U({\cal L}^{\tau,\vec
\ell}_{\iota,-})(f^\tau_{\es,l})^{\lmd^{\tau,\es}_l+1}v^{\tau,\es.}\eqno(5.41)$$

Based on (5.3) and (5.4), we denote
$$\bar{\ell}=(|\iota|+\mbox{max}\:\{\ell_1,\ell_2,...,\ell_n\})n.\eqno(5.42)$$
For any $\bar\ell<s\in\mbb{N}$, we define:
$${\cal L}^{\ast,\vec\ell}_{\iota,s}=\sum_{l,k=s}^{\infty}\mbb{C}
((-1)^\es\la-k_Q-1\ra_{\ell_{k_R}}\E_{l-1/2,-k+1/2}-\la
l_Q+\ell_{l_R}\ra_{\ell_{l_R}}\E_{-k+1/2,l-1/2})\eqno(5.43)$$ and
\begin{eqnarray*}\qquad{\cal L}^{\dg,\vec\ell}_{\iota,s}&=&\sum_{l,k=s}^{\infty}\mbb{C}
((-1)^\es\la-k_Q-1\ra_{\ell_{k_R}}\E_{l-1/2,-k+1/2}\\ & &
-(-1)^{p(l_R)+p(k_R)}\la
l_Q+\ell_{l_R}\ra_{\ell_{l_R}}\E_{-k+1/2,l-1/2}).\hspace{5.4cm}(5.44)\end{eqnarray*}
Suppose that ${\cal M}$ is an ${\cal L}^{\tau,\vec
\ell}_\iota$-module
$$\mbox{generated by a subspace}\;{\cal M}_0\;\mbox{such
that}\;{\cal L}^{\tau,\vec\ell}_{\iota,s}({\cal
M}_0)=\{0\}\;\;\mbox{for some}\;\bar\ell<s\in\mbb{N}.\eqno(5.45)$$
For instance, the above module ${\cal M}^{\tau,\es}$ is such an
module with ${\cal M}_0=\mbb{C}v^{\tau,\es}$ and
$s=\max\{k_0+2,\bar\ell+1\}$. We can also construct $\cal M$ as
Example 3.2.

 Recall the notion defined in (2.60).
For $i,j\in\ol{1,n}$ and $r\in\mbb{N}$, we have
$$(E_{i,j})^\ast_{\vec\ell}(r,z)=\sum_{l\in\mbb{Z}}(E_{i,j}\otimes
t^l\ptl_t^{r+\ell_j}-(-1)^\es E_{j^\ast,i^\ast}\otimes
(-\ptl_t)^rt^l\ptl_t^{\ell_i})z^{-l-1}\eqno(5.46)$$ and
$$(E_{i,j})^\dg_{\vec\ell}(r,z)=\sum_{l\in\mbb{Z}}(E_{i,j}\otimes
t^l\ptl_t^{r+\ell_j}-(-1)^{\es+p(i)+p(j)} E_{j^\ast,i^\ast}\otimes
(-\ptl_t)^rt^l\ptl_t^{\ell_i})z^{-l-1}.\eqno(5.47)$$ By Theorem
2.5, we have a representation $\sgm_\ast$ of
$\hat{o}(\vec\ell,\mbb{A})$ and a representation $\sgm_\dg$ of
$\widehat{sp}(\vec\ell,\mbb{A})$ on ${\cal M}$ with
$\sgm_\ast(\kappa)=\kappa_0,\;\sgm_\dg(\kappa)=\kappa_0$ and
\begin{eqnarray*}& &\sgm_\ast((E_{i,j})^\ast_{\vec\ell}(r,z))\\
&=&\sum_{l,k=0}^{\infty}[\la k\ra_r ((-1)^\es\la-k-1\ra_{\ell_j}
\E_{(-l-1)n+i-1/2,-kn-j+1/2}\\ & &-\la
-l-1\ra_{\ell_i}\E_{-kn-j+1/2,(-l-1)n+i-1/2})z^{l+k-r} \\ & &+\la
k\ra_r ((-1)^\es\la-k-1\ra_{\ell_j}\E_{ln+i-1/2,-kn-j+1/2}\\
&&-\la
l+\ell_i\ra_{\ell_i}\E_{-kn-j+1/2,ln+i-1/2})z^{-l+k-\ell_i-r-1}\\
&&+\la-k-\ell_j-1\ra_r((-1)^\es\la
k+\ell_j\ra_{\ell_j}\E_{-(l+1)n+i-1/2,(k+1)n-j+1/2}\\ & &-\la
-l-1\ra_{\ell_i}\E_{(k+1)n-j+1/2,-(l+1)n+i-1/2})z^{l-k-\ell_j-r-1}\\
& &+\la -k-\ell_j-1\ra_r((-1)^\es\la
k+\ell_j\ra_{\ell_j}\E_{ln+i-1/2,(k+1)n-j+1/2}\\ & &-\la
l+\ell_i\ra_{\ell_i}\E_{(k+1)n-j+1/2,ln+i-1/2})z^{-l-k-\ell_i-\ell_j-r-2}\\
& & + [(r+\ell_i)!\Im_{0,r+\ell_i} -(-1)^\es
r!\ell_i!\Im_{r,\ell_i}]\dlt_{i,j}\kappa_0z^{-r-\ell_i-1}\hspace{5.8cm}(5.48)\end{eqnarray*}
and
\begin{eqnarray*}& &\sgm_\dg((E_{i,j})^\dg_{\vec\ell}(r,z))\\
&=&\sum_{l,k=0}^{\infty}[\la k\ra_r ((-1)^\es\la-k-1\ra_{\ell_j}
\E_{(-l-1)n+i-1/2,-kn-j+1/2}\\ & &-(-1)^{p(i)+p(j)}\la
-l-1\ra_{\ell_i}\E_{-kn-j+1/2,(-l-1)n+i-1/2})z^{l+k-r}\\ && +\la
k\ra_r ((-1)^\es\la-k-1\ra_{\ell_j}\E_{ln+i-1/2,-kn-j+1/2}\\ & &
-(-1)^{p(i)+p(j)}\la
l+\ell_i\ra_{\ell_i}\E_{-kn-j+1/2,ln+i-1/2})z^{-l+k-\ell_i-r-1}\\
&&+\la-k-\ell_j-1\ra_r((-1)^\es\la
k+\ell_j\ra_{\ell_j}\E_{-(l+1)n+i-1/2,(k+1)n-j+1/2}\\ & &
-(-1)^{p(i)+p(j)}\la
-l-1\ra_{\ell_i}\E_{(k+1)n-j+1/2,-(l+1)n+i-1/2})z^{l-k-\ell_j-r-1}\\
& &+\la -k-\ell_j-1\ra_r((-1)^\es\la
k+\ell_j\ra_{\ell_j}\E_{ln+i-1/2,(k+1)n-j+1/2}\\ & &
-(-1)^{p(i)+p(j)}\la
l+\ell_i\ra_{\ell_i}\E_{(k+1)n-j+1/2,ln+i-1/2})z^{-l-k-\ell_i-\ell_j-r-2}\\
& & +[(r+\ell_i)!\Im_{0,r+\ell_i} -(-1)^\es
r!\ell_i!\Im_{r,\ell_i}]\dlt_{i,j}\kappa_0z^{-r-\ell_i-1}\hspace{5.8cm}(5.49)\end{eqnarray*}
By a similar proof as that of Theorem 4.2, we obtain the following
theorem, which was proved by Ma [M2] in a different form when
$\iota=0$.\psp

{\bf Theorem 5.2}. {\it Suppose that $\cal M$ is a weighted ${\cal
L}^{\tau,\vec\ell}_\iota$-module satisfying (5.45). Then the
representation $\sgm_\tau$ is irreducible if and only if $\cal M$
is irreducible.}

\section{Vacuum Representation of
$\widehat{gl}(\vec\ell,\mbb{A})$}

In this section, we study the  vacuum representation of the Lie
algebra $\widehat{gl}(\vec\ell,\mbb{A})$ in (3.6) and its vertex
algebra structure. Its vertex algebra irreducible representations
are investigated.

Recall the algebra $\mbb{A}$ of differential operators given in
(2.2) and (2.3). Observe that
$$\mbb{A}_-=\sum_{i=0}^{\infty}\mbb{C}[t^{-1}]t^{-1}\ptl_t^i\;\;\mbox{and}\;\;
\mbb{A}_+=\sum_{i=0}^{\infty}\mbb{C}[t]\ptl_t^i\eqno(6.1)$$ forms
associative subalgebras of $\mbb{A}$. Moreover,
$$\mbb{A}=\mbb{A}_-+\mbb{A}_+.\eqno(6.2)$$
 Recall the Lie
algebra $\widehat{gl}(n,\mbb{A})$ in (3.3) and (3.4). Note that in
(3.4),
$$\dlt_{r_1+r_2,m_1+m_2}r_1!r_2!\left(\!\!\begin{array}{c}m_1\\
r_1+r_2+1\end{array} \!\!\right)=0\qquad
\for\;\;m_1,m_2,r_1,r_2\in\mbb{N},\eqno(6.3)$$ because if
$m_1+m_2=r_1+r_2$, then $m_1<r_1+r_2+1$. Thus we have the
following Lie subalgebras of $\widehat{gl}(n,\mbb{A})$:
$$\widehat{gl}(n,\mbb{A})_\pm=M_{n\times n}(\mbb{A}_\pm).\eqno(6.4)$$
Moreover,
$$\widehat{gl}(n,\mbb{A})=\widehat{gl}(n,\mbb{A})_-+\widehat{gl}(n,\mbb{A})_+
+\mbb{C}\kappa.\eqno(6.5)$$

Suppose that ${\cal G}$ is Lie subalgebra of
$\widehat{gl}(n,\mbb{A})$ such that
$${\cal G}={\cal G}_-+{\cal G}_++\mbb{C}\kappa,\qquad {\cal
G}_\pm={\cal G}\bigcap \widehat{gl}(n,\mbb{A})_\pm.\eqno(6.6)$$
Then ${\cal G}_\pm$ and
$${\cal G}_0={\cal G}_++\mbb{C}\kappa\eqno(6.7)$$
are Lie subalgebras of ${\cal G}$.  Take a nonzero constant
$\chi\in\mbb{C}$. Form a one-dimensional ${\cal G}_0$-module
$\mbb{C}\vcm$ by:
$$\kappa(\vcm)=\chi\vcm,\;\;{\cal G}_+(\vcm)=\{0\}.\eqno(6.8)$$
The induced ${\cal G}$-module
$${\cal V}_\chi({\cal G})=U({\cal G})\otimes_{U({\cal
G}_0)}\mbb{C}\vcm\cong U({\cal
G}_-)\otimes_{\mbb{C}}\mbb{C}\vcm\eqno(6.9)$$ is called the {\it
vacuum module} of ${\cal G}$ and the corresponding representation
is called the {\it vacuum representation of ${\cal G}$ at level
$\chi$}. The main objective in the rest of paper is to study
${\cal V}_\chi({\cal G})$ and the related vertex algebra structure
when ${\cal G}$ is one of the Lie algebras
$\widehat{gl}(\vec\ell,\mbb{A})$, $\hat{o}(\vec\ell,\mbb{A})$ or
$\widehat{sp}(\vec\ell,\mbb{A})$. For convenience, we simply
denote
$$u\vcm=u\otimes\vcm\qquad\for\;\;u\in U({\cal G}).\eqno(6.10)$$
In the rest of this section, we will deal only with
$\widehat{gl}(\vec\ell,\mbb{A})$.\psp

{\bf Theorem 6.1}. {\it The module ${\cal
V}_\chi(\widehat{gl}(\vec\ell,\mbb{A}))$ is irreducible if
$\chi\not\in\mbb{Z}$. When $\chi\in\mbb{Z}$,  the module ${\cal
V}_\chi(\widehat{gl}(\vec\ell,\mbb{A}))$ has a unique maximal
proper submodule $\bar{\cal
V}_\chi(\widehat{gl}(\vec\ell,\mbb{A}))$, and the quotient
$$V_\chi(\widehat{gl}(\vec\ell,\mbb{A}))={\cal
V}(\widehat{gl}(\vec\ell,\mbb{A}))/\bar{\cal
V}_\chi(\widehat{gl}(\vec\ell,\mbb{A}))\eqno(6.11)$$ is an
irreducible $\widehat{gl}(\vec\ell,\mbb{A})$-module. If $n>1$ and
$\chi\in\mbb{N}$, the submodule}
$$\bar{\cal V}_\chi(\widehat{gl}(\vec\ell,\mbb{A}))=
U(\widehat{gl}(\vec\ell,\mbb{A}))(t^{-1}\ptl^{\ell_1}E_{n,1})^{\chi+1}\vcm.\eqno(6.12)$$

{\it Proof}. We define
$$\widehat{gl}(\vec\ell,\mbb{A}))_{(k)}=\sum_{i,j=1}^n\sum_{r=0}^{\infty}\mbb{C}
t^{r+\ell_j-k}\ptl_t^{r+\ell_j}E_{i,j}+\mbb{C}\dlt_{k,0}\kappa\eqno(6.13)$$
for $k\in\mbb{Z}$. Then
$$\widehat{gl}(\vec\ell,\mbb{A}))=\bigoplus_{k\in\mbb{Z}}\widehat{gl}(\vec\ell,\mbb{A}))_{(k)}
\eqno(6.14)$$ is a $\mbb{Z}$-graded Lie algebra. Moreover, we
define a $\mbb{Z}$-grading on ${\cal
V}_{\chi}(\widehat{gl}(\vec\ell,\mbb{A}))$ by
$${\cal V}_\chi(\widehat{gl}(\vec\ell,\mbb{A}))_{(0)}=\mbb{C}\vcm,
\;\;{\cal
V}_\chi(\widehat{gl}(\vec\ell,\mbb{A}))_{(-m)}=\{0\}\qquad\for\;\;m\in\mbb{N}+1\eqno(6.15)$$
and
$${\cal
V}_\chi(\widehat{gl}(\vec\ell,\mbb{A}))_{(m)}=\mbox{Span}\{u_1u_2\cdots
u_s\vcm\mid u_i\in \widehat{gl}(\vec\ell,\mbb{A})_-\bigcap
\widehat{gl}(\vec\ell,\mbb{A})_{(k_i)};\;\sum_{i=1}^sk_i=m\}\eqno(6.16)$$
for $m\in\mbb{N}+1$. Then
$${\cal
V}_\chi(\widehat{gl}(\vec\ell,\mbb{A}))=\bigoplus_{k\in\mbb{Z}}{\cal
V}_\chi(\widehat{gl}(\vec\ell,\mbb{A}))_{(k)}\eqno(6.17)$$ is a
$\mbb{Z}$-graded $\widehat{gl}(\vec\ell,\mbb{A})$-module. Since
$$(\mbb{A}\ptl_t^{\ell_j}E_{i,j})\bigcap\widehat{gl}(\vec\ell,\mbb{A}))_{(k+\ell_j)}
\bigcap\widehat{gl}(\vec\ell,\mbb{A}))_-=\sum_{r=0}^{k-1}\mbb{C}
t^{r-k}\ptl_t^{r+\ell_j}E_{i,j}\eqno(6.18)$$ (cf. (6.4) and
(6.6)), the character
$$d({\cal
V}_\chi(\widehat{gl}(\vec\ell,\mbb{A})),q)=\sum_{k=0}^{\infty}(\dim{\cal
V}_\chi(\widehat{gl}(\vec\ell,\mbb{A}))_{(k)})q^k=
\prod_{i=1}^n\prod_{r=1}^{\infty}\frac{1}{(1-q^{\ell_i+r})^{rn}}.\eqno(6.19)$$

Recall the Lie algebra $\td{gl}(\infty)$ defined in (3.9) and
(3.10). Set
$$\td{gl}(\infty)_{(-)}=\sum_{0>j,k\in\Z}\mbb{C}\E_{j,k},\;\;
\td{gl}(\infty)_{(+)}=
\sum_{0<j,k\in\Z}(\mbb{C}\E_{j,k}+\mbb{C}\E_{j,-k}+\mbb{C}\E_{-j,k}).\eqno(6.20)$$
By (3.10), $\td{gl}(\infty)_{(\pm)}$ are Lie subalgebras of
$\td{gl}(\infty)$ and
$$\td{gl}(\infty)=\td{gl}(\infty)_{(-)}+\td{gl}(\infty)_{(+)}+\mbb{C}\kappa_0.\eqno(6.21)$$
Hence we have the Lie subalgebra
$$\td{gl}(\infty)_{(0)}=\td{gl}(\infty)_{(+)}+\mbb{C}\kappa_0.\eqno(6.22)$$
Define a one-dimensional $\td{gl}(\infty)_{(0)}$-module
$\mbb{C}{\bf 1}$ by
$$\kappa_0({\bf 1})=\chi
{\bf 1},\;\;\td{gl}(\infty)_{(+)}({\bf 1})=\{0\}.\eqno(6.23)$$
Form an induced $\td{gl}(\infty)$-module
$$U_\chi=U(\td{gl}(\infty))\otimes_{U(\td{gl}(\infty)_{(-)})}\mbb{C}{\bf 1}\cong
U(\td{gl}(\infty)_{(-)})\otimes_{\mbb{C}}\mbb{C}{\bf
1},\eqno(6.24)$$ which satisfies the condition (3.18) with $m=0$
and ${\cal M}_0=\mbb{C}1\otimes{\bf 1}$. For convenience, we
denote
$$ v{\bf 1}=v\otimes {\bf 1}\qquad\for\;\;v\in
U(\td{gl}(\infty)).\eqno(6.25)$$

Note  our notion
$$
E_{i,j}(r,z)=\sum_{m\in\mbb{Z}}t^m\ptl_t^rz^{-m-1}E_{i,j}
\qquad\for\;\;i,j\in\ol{1,n},\;r\in\mbb{N}+\ell_j.\eqno(6.26)$$
According to (3.37) (also cf. (2.52) and (3.11)), we obtain a
$\widehat{gl}(\vec\ell,\mbb{A})$-module structure  on $U_\chi$
defined by $\kappa=\chi\mbox{Id}_{U_\chi}$ and
\begin{eqnarray*}\hspace{2cm}E_{i,j}(r,z) &=&\sum_{l,k=0}^\infty[\la-k-1\ra_r(\E_{ln+i-1/2,(k+1)n-j+1/2}
z^{-l-k-\ell_i-r-2}\\ & &+\E_{-(l+1)n+i-1/2,(k+1)n-j+1/2} z^{l-k-r-1})\\
& &+\la
k+\ell_j\ra_r(\E_{ln+i-1/2,-kn-j+1/2}z^{-l+k+\ell_j-\ell_i-r-1}\\
& & +
\E_{-(l+1)n+i-1/2,-kn-j+1/2}z^{l+k+\ell_j-r})]\hspace{4.3cm}(6.27)\end{eqnarray*}
for $i,j\in\ol{1,n}$ and $r\in\mbb{N}+\ell_j$. In particular,
$$E_{i,j}(r,z)({\bf 1})=\sum_{l,k=0}^\infty \E_{-(l+1)n+i-1/2,-kn-j+1/2}{\bf
1}z^{l+k+\ell_j-r}.\eqno(6.28)$$ Thus we have
$$\widehat{gl}(\vec\ell,\mbb{A})_+({\bf 1})=\{0\}.\eqno(6.29)$$
By a similar proof as that of Theorem 3.1,
$$U_\chi=U(\widehat{gl}(\vec\ell,\mbb{A})_-){\bf 1}.\eqno(6.30)$$
Therefore, we have a Lie algebra module epimorphism $\nu:{\cal
V}_\chi(\widehat{gl}(\vec\ell,\mbb{A}))\rta U_\chi$ defined by
$$\nu(u\vcm)=u{\bf 1}\qquad\for\;\;u\in U(\widehat{gl}(\vec\ell,\mbb{A})_-).\eqno(6.31)$$

Set
$$\ell=\mbox{min}\:\{\ell_1,\ell_2,...,\ell_n\}.\eqno(6.32)$$
For $m\in\mbb{N}+\ell+1$, we let
$$\td{gl}(\infty)_{(-)}^{(m)}=\mbox{Span}\:\{\E_{-(l+1)n+i-1/2,-k-j+1/2}\mid
i,j\in\ol{1,n},\;l,k\in\mbb{N};\;l+k+\ell_j+1=m\}.\eqno(6.33)$$
Then
$$\td{gl}(\infty)_{(-)}=\bigoplus_{m=\ell+1}^{\infty}\td{gl}(\infty)_{(-)}^{(m)}.\eqno(6.34)$$
Moreover, we define
$$U_\chi^{(0)}=\mbb{C}{\bf
1},\;\;U_\chi^{(m)}=\{0\}\qquad\for\;\;m\in(-\mbb{N}-1)\bigcup\ol{1,\ell}\eqno(6.35)$$
and
$$U_\chi^{(m)}=\mbox{Span}\:\{u_1u_2\cdots u_s{\bf 1}\mid u_i\in
\td{gl}(\infty)_{(-)}^{(m_i)};\;\sum_{i=1}^sm_i=m\}.\eqno(6.36)$$
By (6.30),
$$ U_\chi=\bigoplus_{m\in\mbb{Z}}U_\chi^{(m)}\eqno(6.37)$$
is a $\mbb{Z}$-graded $\widehat{gl}(\vec\ell,\mbb{A})$-module.
Furthermore, (6.33) implies the character
$$d(U_\chi,q)=\sum_{k=0}^{\infty}(\dim U_\chi^{(k)})q^k=
\prod_{i=1}^n\prod_{r=1}^{\infty}\frac{1}{(1-q^{\ell_i+r})^{rn}}.\eqno(6.38)$$
Therefore, (6.19) and (6.38) imply
$${\cal V}_\chi(\widehat{gl}(\vec\ell,\mbb{A}))\cong U_\chi.\eqno(6.39)$$

Let $\lmd$ be a linear function on $\T$ (cf. (3.12)) such that
$$\lmd(\kappa_0)=\chi,\;\;\lmd(\E_{l,-l})=0\qquad\for\;\;l\in\Z.\eqno(6.40)$$
Recall the Verma module $M_\lmd$ defined in (3.63). Note
$$U_\chi\cong M_\lmd/(\sum_{-1/2\neq
l\in\Z}U(\td{gl}(\infty)_-)(\E_{l+1,l}\otimes
v_\lmd)),\eqno(6.41)$$ which is irreducible if
$\chi\not\in\mbb{Z}$ by [J1-J3]. When $\chi\in\mbb{N}$,
\begin{eqnarray*} & &\bar{U}_\chi=U(\td{gl}(\infty)_-)\E_{-1/2,-1/2}^{\chi+1}{\bf
1}\\ &&\cong
(\sum_{l\in\Z}U(\td{gl}(\infty)_-)(\E_{l+1,l}^{\lmd_l+1}\otimes
v_\lmd))/(\sum_{-1/2\neq
l\in\Z}U(\td{gl}(\infty)_-)(\E_{l+1,l}^{\lmd_l+1}\otimes
v_\lmd))\hspace{3cm}(6.42)\end{eqnarray*} is the unique maximal
proper submodule of $U_\chi$ (cf. (3.65)). Thus
$$U(\widehat{gl}(\vec\ell,\mbb{A}))\nu^{-1}(\E_{-1/2,-1/2}^{\chi+1}{\bf
1})\eqno(6.43)$$ is the unique maximal proper submodule of ${\cal
V}_\chi(\widehat{gl}(\vec\ell,\mbb{A}))$ (cf. (6.31)). When $n>1$,
(6.23) and (6.27) imply
$$\nu^{-1}(\E_{-1/2,-1/2}^{\chi+1}{\bf
1})=(t^{-1}\ptl_t^{\ell_1}E_{n,1})^{\chi+1}\vcm.\qquad\Box\eqno(6.44)$$
\vspace{0.1cm}

Next we want to present the definitions of vertex algebra and its
module. For any two vector spaces $U$ and $W$, we denote by
$\mbox{LM}(U,W)$ the set of all linear maps from $U$ to $W$. Let
$z_1$ and $z_2$ be two formal variables. We have the following
convention of binomial expansions:
$$(z_1-z_2)^r=\sum_{l=0}^{\infty}(-1)^l(^{\;r}_{\;l})z_1^{r-l}z_2^l\qquad
\mbox{for}\;\;r\in \Bbb{C}.\eqno(6.45)$$ For a vector space $V$,
we denote
$$V[z^{-1},z]]=\{\sum_{i=m}^\infty v_iz^i\mid v_i\in
V,\;m\in\mbb{Z}\},\eqno(6.46)$$ the space of formal  Laurent
series with coefficients in $V$.

A {\it vertex  algebra} is a vector space $V$  with a linear map
$Y(\cdot, z):\;V\rightarrow \mbox{ LM}(V, V[z^{-1}:z]])$, an
element $\ptl\in \mbox{End}\:V$ and an element $\vcm \in V$,
satisfying the following conditions: given $u,v\in V$,
$$Y(\vcm,z)=\mbox{Id}_V,\eqno(6.47)$$
$$[\ptl,Y(v,z)]=\frac{d}{dz}Y(v,z),\;\;\; Y(v,z)\vcm=e^{z\ptl}v,\eqno(6.48)$$
$$ (z_1-z_2)^mY(u,z_1)Y(v,z_2)=(z_1-z_2)^mY(v,z_2)Y(u,z_1)\eqno(6.49)$$
for some positive integer $m$.

An ideal $U$ of a vertex algebra $V$ is a vector space of $V$ such
that
$$\ptl(U)\subset U,\;\;Y(v,z)U\subset U[z^{-1},z]]\;\;\for\;\;v\in
V.\eqno(6.50)$$ A vertex algebra without proper nonzero ideals is
called {\it simple}.

A module $W$ of a vertex algebra $(V,Y(\cdot,z),\vcm,\ptl)$ is a
vector space with a linear map $Y_W(\cdot, z):\;V\rightarrow
\mbox{ LM}(W, W[z^{-1}:z]])$ such that given $u,v\in V$ and $w\in
W$,
$$Y_W(\vcm,z)=\mbox{Id}_W,\;\;Y_W(\ptl
v,z)=\frac{d}{dz}Y(v,z),\eqno(6.51)$$
$$(z_1-z_2)^mY_W(u,z_1)Y_W(v,z_2)=(z_1-z_2)^mY_W(v,z_2)Y_W(u,z_1),\eqno(6.52)$$
$$(z_0+z_2)^mY_W(u,z_0+z_2)Y_W(v,z_2)w=(z_2+z_0)^mY_W(Y(u,z_0)v,z_2)w\eqno(6.53)$$
for some positive
integer $m$.

A {\it submodule} $W_1$ of $W$ is a subspace such that
$$ Y_W(v,z)W_1\subset W_1[z^{-1},z]]\qquad\for\;\;v\in
V.\eqno(6.54)$$ A module without proper nonzero submodule is
called {\it irreducible}.

For a vector space $U$ and any formal power series
$$f(z)=\sum_{r\in\mbb{Z}}u_rz^{-r-1}\;\;\mbox{with}\;\;u_r\in
U,\eqno(6.55)$$ we define
$$f(z)^+=\sum_{r=0}^{\infty}u_rz^{-r-1},\qquad
f(z)^-=\sum_{r=0}^\infty u_{-r-1}z^r.\eqno(6.56)$$ Now we define a
linear transformation $\ptl$ on $\widehat{gl}(\vec\ell,\mbb{A})$
by
$$\ptl(\kappa)=0,\;\;\ptl(t^m\ptl_t^rE_{i,j})=-mt^{m-1}\ptl_t^rE_{i,j}\eqno(6.57)$$
for $i,j\in\ol{1,n},\;r\in\mbb{N}$ and $m\in\mbb{Z}$. Moreover, we
define a linear transformation $\ptl$ on ${\cal
V}_\chi(\widehat{gl}(\vec\ell,\mbb{A}))$ by
$$\ptl(\vcm)=0,\;\;\ptl(u_1u_2\cdots u_s\vcm)=\sum_{i=1}^su_1\cdots u_{i-1}
\ptl(u)_iu_{i+1}\cdots u_s\vcm\eqno(6.58)$$ for
$u_i\in\widehat{gl}(\vec\ell,\mbb{A})_-$. For $i,j\in\ol{1,n}$ and
$r\in\mbb{N}+\ell_j$, we denote
$$E_{i,j}(r,z)=\sum_{l\in\mbb{Z}}^{\infty}t^l\ptl_t^rE_{i,j}z^{-l-1}.\eqno(6.59)$$
Furthermore, we define linear maps
$$Y^\pm(\cdot,z):\widehat{gl}(\vec\ell,\mbb{A})_-\rta LM({\cal
V}_\chi(\widehat{gl}(\vec\ell,\mbb{A})),{\cal
V}_\chi(\widehat{gl}(\vec\ell,\mbb{A}))[z^{-1},z]])\eqno(6.60)$$
by
$$Y^\pm(t^{-m-1}\ptl_t^rE_{i,j},z)=\frac{1}{m!}\frac{d^m}{dz^m}E_{i,j}(r,z)^{\pm}.\eqno(6.61)$$
Now we define a linear map
$$Y(\cdot,z):{\cal V}_\chi(\widehat{gl}(\vec\ell,\mbb{A}))\rta LM({\cal
V}_\chi(\widehat{gl}(\vec\ell,\mbb{A})),{\cal
V}_\chi(\widehat{gl}(\vec\ell,\mbb{A}))[z^{-1},z]])\eqno(6.62)$$
by induction:
$$ Y(\vcm,z)=\mbox{Id}_{{\cal
V}_\chi(\widehat{gl}(\vec\ell,\mbb{A}))},\;\;Y(uv,z)=Y^-(u,z)Y(v,z)+Y(v,z)Y^+(u,z)\eqno(6.63)$$
for $u\in \widehat{gl}(\vec\ell,\mbb{A})_-$ and $v\in {\cal
V}_\chi(\widehat{gl}(\vec\ell,\mbb{A}))$. In particular,
$$Y((t^{-m-1}\ptl_t^rE_{i,j})\vcm,z)=\frac{1}{m!}\frac{d^m}{dz^m}E_{i,j}(r,z)
\qquad\for\;\;i,j\in\ol{1,n},\;r\in\mbb{N}+\ell_j,\;m\in\mbb{N}.\eqno(6.64)$$
 By Lemma 2.1 and the general theory of conformal algebras (cf.
[K] and [X2]), we have:\psp

{\bf Theorem 6.2}. {\it The family $({\cal
V}_\chi(\widehat{gl}(\vec\ell,\mbb{A})),Y(\cdot,z),\ptl,\vcm)$
forms a vertex algebra . If $\chi\not\in \mbb{Z}$, the vertex
algebra $({\cal
V}_\chi(\widehat{gl}(\vec\ell,\mbb{A})),Y(\cdot,z),\ptl,\vcm)$ is
simple. When $\chi\in\mbb{N}$, the quotient space
$V_\chi(\widehat{gl}(\vec\ell,\mbb{A}))$ forms a simple vertex
algebra. If $\ell_i\in\{0,1\}$ for $i\in\ol{1,n}$, ${\cal
V}_\chi(\widehat{gl}(\vec\ell,\mbb{A}))$ with
$\chi\in\mbb{C}\setminus\mbb{Z}$ and
$V_\chi(\widehat{gl}(\vec\ell,\mbb{A}))$ with $\chi\in\mbb{Z}$ are
simple vertex operator algebras with the Virasoro element
$-t^{-1}\ptl_t(\sum_{i=1}^nE_{i,i})\vcm$.} \psp

Take constants $\chi,\iota\in\mbb{C}$ such that $\chi\neq 0$. Let
${\cal M}$ be a weighted irreducible $\td{gl}(\infty)$-module
satisfying (3.18) (also cf. (3.17)) and $\kappa_0|_{\cal
M}=\chi\mbox{Id}_{\cal M}$. For $i,j\in\ol{1,n}$ and
$r\in\mbb{N}+\ell_j$, we denote
$$E_{i,j}^\iota(r,z)=\sum_{l,k\in\mbb{Z}}\la
\iota-k\ra_r\E_{ln+i-1/2,kn-j+1/2}z^{-l-k-r-1}+r!\Im_{0,r}\dlt_{i,j}\kappa_0
z^{-r-1}\eqno(6.65)$$ if $\iota\not\in\mbb{Z}$, and
\begin{eqnarray*} & & E_{i,j}^\iota(r,z)\\&=&
\sum_{l,k=0}^n[\la
-k-1\ra_r\E_{(l-\iota)n+i-1/2,(k+\iota+1)n-j+1/2}z^{-l-k-\ell_i-r-2}\\
& &+\la \la k+\ell_j\ra_
r\E_{(l-\iota)n+i-1/2,(\iota-k)n-j+/2}z^{-l+k+\ell_j-\ell_i-r-1}\\
& &+\la k+\ell_j \ra_r\E_{-(l+\iota+1)n+i-1/2,(\iota-k)n-j+1/2}z^{l+k+\ell_j-r}\\
& & +\la-k-1\ra_r
\E_{-(l+\iota+1)n+i-1/2,(k+\iota+1)n-j+1/2}z^{l-k-r-1}]+\dlt_{i,j}r!\Im_{0,r}\kappa_0z^{-r-1}
\hspace{2.1cm}(6.66)\end{eqnarray*} (cf. (2.26)). Furthermore, we
define linear maps
$$Y^{\iota,\pm}_{\cal M}(\cdot,z):\widehat{gl}(\vec\ell,\mbb{A})_-\rta LM({\cal
M},{\cal M}[z^{-1},z]])\eqno(6.67)$$ by
$$Y^{\iota,\pm}_{\cal M}(t^{-m-1}\ptl_t^rE_{i,j},z)=\frac{1}{m!}\frac{d^m}{dz^m}
E^\iota_{i,j}(r,z)^{\pm}.\eqno(6.68)$$ Now we define a linear map
$$Y^\iota_{\cal M}(\cdot,z):{\cal V}_\chi(\widehat{gl}(\vec\ell,\mbb{A}))\rta
LM({\cal M},{\cal M} [z^{-1},z]])\eqno(6.69)$$ by induction:
$$ Y^\iota_{\cal M}(\vcm,z)=\mbox{Id}_{{\cal
M}},\;\;Y(uv,z)=Y^{\iota,-}_{\cal M}(u,z)Y_{\cal M}
^\iota(v,z)+Y^\iota_{\cal M}(v,z)Y^{\iota,+}_{\cal
M}(u,z)\eqno(6.70)$$ for $u\in \widehat{gl}(\vec\ell,\mbb{A})_-$
and $v\in {\cal V}_\chi(\widehat{gl}(\vec\ell,\mbb{A})).$ By
Theorem 3.1, Theorem 3.2 and the general theory for vertex
algebras (e.g. cf. Section 4.1 in [X2]), we have:\psp

{\bf Theorem 6.3}. {\it The family $({\cal M} ,Y^\iota_{\cal
M}(\cdot,z))$ forms a irreducible module of the vertex algebra
$({\cal
V}_\chi(\widehat{gl}(\vec\ell,\mbb{A})),Y(\cdot,z),\ptl,\vcm)$.}
\psp

In the rest of this section, we want to show certain familiar
unitary highest weight irreducible $\td{gl}(\infty)$-modules
induce irreducible modules of the quotient simple vertex algebra
$(V(\widehat{gl}(\vec\ell,\mbb{A})),Y(\cdot,z),\\ \ptl,\vcm)$ when
$\chi\in\mbb{Z}$. To this end, we need to use charged free fields.

First we use charged free fermionic field realization. Denote
$$\Z_+=\mbb{N}+1/2,\qquad \Z_-=-\Z_+.\eqno(6.71)$$
Then
$$\Z=\Z_+\bigcup\Z_-.\eqno(6.72)$$
Let $\{\cta_l,\bar{\cta}_l\mid l\in\Z_-\}$ be a set of odd
variables, that is,
$$\cta_l\cta_k=-\cta_k\cta_l,\;\;\cta_l\bar\cta_k=-\bar\cta_k\cta_l,\;\;
\bar\cta_l\bar\cta_k=-\bar\cta_k\bar\cta_l.\eqno(6.73)$$ Set
$$V_f=\mbb{C}[\cta_l,\bar{\cta}_l\mid l\in\Z_-],\eqno(6.74)$$
the polynomial algebra in odd variables $\{\cta_l,\bar{\cta}_l\mid
l\in\Z_-\}$, where the subindex ``f" stands for ``fermionic".
Moreover, we denote
$$\cta_l=\ptl_{\bar\cta_{-l}},\qquad
\bar\cta_l=\ptl_{\cta_{-l}}\qquad\for\;\;l\in\Z_+.\eqno(6.75)$$ It
can be verified that we the following representation of
$\td{gl}(\infty)$ on $V_f$:
$$\kappa_0|_{V_f}=\mbox{Id}_{V_f},\;\;\E_{l,k}|_{V_f}=\left\{\begin{array}{ll}
\bar\cta_l\cta_k&\mbox{if}
\;l\in\Z_-\;\mbox{or}\;-k\neq l\in\Z_+,\\
-\cta_k\bar\cta_l&\mbox{if}\;\;-k=l\in\Z_+\end{array}\right.\eqno(6.76)$$
for $l,k\in\Z$. Set
$$\ptl=\sum_{l\in\Z}l\E_{-1-l,l}.\eqno(6.77)$$
Given $i\in\ol{1,n}$, we define
$$\cta(\iota,i,z)=\sum_{l\in\mbb{Z}}\cta_{ln-i+1/2}z^{\iota-l},\;\;\bar\cta(\iota,i,z)
=\sum_{l\in\mbb{Z}}\bar\cta_{ln+i-1/2}z^{-\iota-l-1}\eqno(6.78)$$
for $\iota\in\mbb{C}\setminus\mbb{Z}$, and
$$\cta(\iota,i,z)=\sum_{l=0}^\infty[\cta_{(-l+\iota)n-i+1/2}z^{\ell_i+l}
+\cta_{(l+\iota+1)n-i+1/2}z^{-l-1}], \eqno(6.79)$$
$$\bar\cta(\iota,i,z)=\sum_{l=0}^\infty[\bar\cta_{(l-\iota)n+i-1/2}z^{-\ell_i-l-1}
+\bar\cta_{-(l+\iota+1)+i-1/2}z^l]\eqno(6.80)$$
 for $\iota\in\mbb{Z}$. Then $\{\cta(\iota,i,z),\bar\cta(\iota,i,z)\mid
 i\in\ol{1,n}\}$ are {\it charged free fermionic fields}.

  Set
$$\Theta=\sum_{l\in\Z_-}\mbb{C}\cta_l,\qquad
\bar\Theta=\sum_{l\in\Z_-}\mbb{C}\bar\cta_l,\eqno(6.81)$$ and
$$V_{f,0}=\mbb{C}+\sum_{s=1}^{\infty}\bar\Theta^s\Theta^s,\;\;V_{f,r}=V_{f,0}\Theta^r,
\;\;V_{f,-r}=\bar\Theta^rV_{f,0}\eqno(6.82)$$ for $r\in\mbb{N}+1$.
Then
$$V_f=\bigoplus_{k\in\mbb{Z}}V_{f,k}.\eqno(6.83)$$
We define a linear map $Y^\iota(\cdot,z):V_f\rta LM(V_f,V_f\{z\})$
by induction:
$$Y^\iota(1,z)=\mbox{Id}_{V_f},\eqno(6.84)$$
\begin{eqnarray*}Y^\iota(\cta_{-rn-i+1/2}u,z)&=&\rd_{z_1}\sum_{s=0}^\infty(-1)^s\left(\!\!\begin{array}{c}-\iota\\
s\end{array}\!\!\right)z_1^{-\iota-s}z^\iota[(z_1-z)^{s-r-\ell_i-1}\cta(\iota,i,z_1)Y^\iota(u,z)\\
&
&-(-1)^k(-z+z_1)^{s-r-\ell_i-1}Y^\iota(u,z)\cta(\iota,i,z_1)]\hspace{3.7cm}(6.85)\end{eqnarray*}
and
\begin{eqnarray*}Y^\iota(\bar\cta_{-(r+1)n+i-1/2}u,z)&=&\rd_{z_1}\sum_{s=0}^\infty(-1)^s
\left(\!\!\begin{array}{c}\iota\\
s\end{array}\!\!\right)z_1^{\iota-s}z^{-\iota}[(z_1-z)^{s-r-1}\bar\cta(\iota,i,z_1)Y^\iota(u,z)
\\ &&-(-1)^k(-z+z_1)^{s-r-1}Y^\iota(u,z)\bar\cta(\iota,i,z_1)]\hspace{3.6cm}(6.86)\end{eqnarray*}
for $i\in\ol{1,n},\;r\in\mbb{N}$ and $u\in V_{f,k}$. By Section
4.1 in [X2], we have
\begin{eqnarray*}& & z_0^{-1}\dlt\left({z_1-z_2\over z_0}\right)Y^\iota(u,z_1)Y^\iota(v,z_2)
-(-1)^{k_1k_2}z_0^{-1}\dlt\left({z_2-z_1\over
-z_0}\right)Y^\iota(v,z_2)Y^\iota(u,z_1)
\\&=& z_2^{-1}\left({z_1-z_0\over z_2}\right)^{k_1\iota}\dlt\left({z_1-z_0\over z_2}\right)
Y^\iota(Y^0(u,z_0)v,z_2)\hspace{6cm}(6.87)\end{eqnarray*} for
$u\in V_{f,k_1}$ and $v\in V_{f,k_2}$. In fact,
$$E^\iota_{i,j}(r,z)|_{V_f}=Y^\iota(\bar\cta_{-n+i-1/2}\cta_{-rn-j+1/2},z)\eqno(6.88)$$
(cf. (6.65), (6.66)),
$$Y^\iota(\cta_{-rn+i-1/2},z)=\frac{1}{(r+\ell_i)!}\frac{d^{r+\ell_i}}{z^{r+\ell_i}}
\cta(\iota,i.z)\eqno(6.89)$$ and
$$Y^\iota(\bar\cta_{-(r+1)n-i+1/2},z)=\frac{1}{r!}\frac{d^r}{z^r}\bar\cta(\iota,i.z)\eqno(6.90)$$
for $i,j\in\ol{1,n}$.

Let $k\in\mbb{N}+1$. Denote
$$v_0=1,\;\;v_k=\cta_{-1/2}\cta_{-3/2}\cdots\cta_{1/2-k},\;\;
v_{-k}=\bar\cta_{-1/2}\bar\cta_{-3/2}\cdots\bar\cta_{1/2-k}.\eqno(6.91)$$
Moreover, we define linear functions on $\T$ in (3.10):
$$\lmd^0(\kappa_0)=1,\qquad\lmd^0(\E_{l,-l})=0\qquad\for\;\;l\in\Z,\eqno(6.92)$$
$$\lmd^{-k}(\kappa_0)=\lmd^{-k}(\E_{-r-1/2,r+1/2})=1,\;\;\lmd^{-k}(E_{-l-1/2,l+1/2})=0\eqno(6.93)$$ for
for $r\in\ol{0,k-1}$ and $l\in\mbb{Z}\setminus\ol{0,k-1}$, and
$$\lmd^k(\kappa_0)=-\lmd^k(\E_{r+1/2,-r-1/2})=1,\;\;\lmd^k(E_{l+1/2,-l-1/2})=0\eqno(6.94)$$
for $r\in\ol{0,k-1}$ and $l\in\mbb{Z}\setminus\ol{0,k-1}$. \psp

{\bf Lemma 6.4}. {\it Each space $V_{f,k}$ forms a unitary
irreducible highest weight module of $\td{gl}(\infty)$ with the
highest weight $\lmd^k$. Moreover, a linear function $\lmd$ on
$\T$ satisfies (3.59), (3.65) and $\lmd(\kappa_0)=1$, if and only
if $\lmd$ is of the form $\lmd^k$ for some $k\in\mbb{Z}$. The
family $(V_{f,0},Y^0(\cdot,z),1,\ptl)$ is a simple vertex operator
algebra isomorphic to
$(V_1(\widehat{gl}(\vec\ell,\mbb{A})),Y(\cdot,z),\vcm,\ptl)$ in
Theorem 6.2. Moreover, each
$(V_{f,k},Y^\iota(|_{V_{f,0}},z)|_{V_{f,k}})$ is an irreducible
$V_{f,0}$-module. The map $Y^\iota(|_{V_{f,l}},z)|_{V_{f,k}}$ is
an intertwining operator of type
$[_{V_{f,l}}^{\;\;\;\;\;\;V_{f,l+k}}\:_{V_{f,k}}]$.} \psp

Now we want to deal with higher level case. Assume
$1<\chi\in\mbb{N}$. Form $\chi$th $\td{gl}(\infty)$-module tensor:
$$V_f^{\la \chi\ra}=V_f\otimes V_f\otimes\cdots\otimes
V_f\qquad(\chi\;\mbox{copies}).\eqno(6.95)$$ Then
$$\kappa_0|_{V_f^{\la \chi\ra}}=\chi\mbox{Id}_{V_f^{\la
\chi\ra}}.\eqno(6.96)$$ Since $V_f$ is an unitary
$\td{gl}(\infty)$-module that is a direct sum of the highest
weight irreducible modules $V_{f,k}$ with weights $\lmd^k$
satisfying (3.59), we can apply the tensor theory of
finite-dimensional general linear Lie algebras to $V_f^{\la
\chi\ra}$. Thus $V_f^{\la \chi\ra}$ is a completely reducible
$\td{gl}(\infty)$-module. Moreover, for each weight $\lmd$
satisfying (3.59), (3.65) and $\lmd(\kappa_0)=\chi$, there exists
a component of $V_f^{\la \chi\ra}$ that is a highest weight
irreducible $\td{gl}(\infty)$-module with weight $\lmd$.

Denote
$$1_\chi=1\otimes 1\otimes\cdots\otimes 1\qquad(\chi\;\mbox{copies}).\eqno(6.97)$$
Since $V_f$ is also a polynomial algebra in odd variables,
$V_f^{\la \chi\ra}$ has an extended associative tenor algebra
structure. Note the space
$$\bar\Theta\Theta=\sum_{l,k\in\Z_-}\mbb{C}\bar\cta_l\cta_k.\eqno(6.98)$$
We define the diagonal linear map $\varrho:\bar\Theta\Theta\rta
V_f^{\la \chi\ra}$ by
$$\varrho(\bar\cta_l\cta_k)=\sum_{i=1}^\chi 1\otimes\cdots\otimes 1\otimes\stl{i}
{\bar\cta_l\cta_k}\otimes 1\otimes\cdots\otimes 1\eqno(6.99)$$ for
$l,k\in\Z_-$. Set
$$V_{f,0}^{[\chi]}=\mbb{C}1_\chi+\sum_{r=1}^{\infty}[\varrho(\bar\Theta\Theta)]^r.\eqno(6.100)$$
Moreover, define the map
$$Y^\iota_\chi(\cdot,z)=Y^\iota(\cdot,z)\otimes
Y^\iota(\cdot,z)\otimes\cdots\otimes
Y^\iota(\cdot,z)\qquad(\chi\;\mbox{copies})\eqno(6.101)$$ and
$$\ptl^{\la\chi\ra}=\sum_{i=1}^\chi 1\otimes\cdots\otimes 1\otimes\stl{i}
{\ptl}\otimes 1\otimes\cdots\otimes 1\eqno(6.102)$$ (cf. (6.77)).
 Then the family
$(V_{f,0}^{[\chi]},Y^0_\chi(\cdot,z),1_\chi,\ptl^{\la\chi\ra})$
forms a simple vertex algebra isomorphic to
$(V_\chi(\widehat{gl}(\vec\ell,\mbb{A})), Y(\cdot,z),\vcm,\ptl)$
in Theorem 6.2. For each irreducible $\td{gl}(\infty)$-module
component $U$ of $V_f^{\la \chi\ra}$, the family
$(U,Y^\iota_\chi(|_{V_{f,0}^{[\chi]}},z)|_U)$ forms an irreducible
module of
$(V_{f,0}^{[\chi]},Y^0_\chi(\cdot,z),1_\chi,\ptl^{\la\chi\ra})$.
 \psp

{\bf Lemma 6.5}. {\it  Suppose that $\chi$ is a positive integer
and $\lmd$ is a weight of $\td{gl}(\infty)$ satisfying (3.59),
(3.65) and $\lmd(\kappa_0)=\chi$. Let $\cal M$ be the irreducible
highest weight $\td{gl}(\infty)$-module with highest weight
$\lmd$. Then the family $({\cal M},Y^\iota_{\cal M})$ defined in
(6.65)-(6.70) forms an irreducible module of the simple vertex
algebra $(V_\chi(\widehat{gl}(\vec\ell,\mbb{A})),Y(\cdot,z),
\vcm,\ptl)$ in Theorem 6.2, equivalently,
$$Y^\iota_{\cal
M}((t^{-1}\ptl^{\ell_1}E_{n,1})^{\chi+1}\vcm,z)=0\eqno(6.103)$$
when $n>1$.}\psp

We remark that (6.103) can be proved easily by using the affine
Lie algebra $\widehat{sl}(2,\mbb{C})$ when $\ell_1=\ell_n=0$. Let
$\cal M$ be the highest weight irreducible module in Example 3.1.
By the locality of $Y^\iota_{\cal M}(\cdot,z)$, (6.103) holds if
and only if
$$Y^\iota_{\cal
M}((t^{-1}\ptl^{\ell_1}E_{n,1})^{\chi+1}\vcm,z)v_\lmd=0.\eqno(6.104)$$
According to  the above lemma with $n=2$, (4.104) holds if
$$
\lmd(\E_{rn+1/2,-rn-1/2}-\E_{rn-1/2,-rn+1/2}+\dlt_{r,0}\kappa_0)\in\mbb{N}
\qquad\for\;\;r\in\mbb{Z}\eqno(4.105)$$ and when $n>2$
$$\lmd(\E_{rn+n-1/2,-rn-n+1/2}-\E_{rn+1/2,-rn-1/2}+\dlt_{r,0})\in\mbb{N}
\qquad\for\;\;r\in\mbb{Z},\eqno(4.106)$$ Note that the condition
of (6.105) and (6.106) is weaker than (3.65) when $n>3$. Under the
condition of (6.105) and (6.106), $\cal M$ may not be a unitary
$\td{gl}(\infty)$-module. Therefore, we have proved the following
main theorem: \psp

{\bf Theorem 6.6}. {\it  Suppose that $\chi$ is a positive integer
and $\lmd$ is a weight of $\td{gl}(\infty)$ satisfying (3.59),
(6.105), (6.106) and $\lmd(\kappa_0)=\chi$. Let $\cal M$ be the
irreducible highest weight $\td{gl}(\infty)$-module with highest
weight $\lmd$. Then the family $({\cal M},Y^\iota_{\cal M})$
defined in (6.65)-(6.70) forms an irreducible module of the simple
vertex algebra
$(V_\chi(\widehat{gl}(\vec\ell,\mbb{A})),Y(\cdot,z),\vcm,\ptl)$ in
Theorem 6.2.} \psp

Next we use charged free bosonic field realization to study the
case of negative integral $\chi$.
 Set
$$V_b=\mbb{C}[x_l,\bar x_l\mid l\in\Z_-],\eqno(6.107)$$
the polynomial algebra in  a set of ordinary commute variables
$\{x_l,\bar x_l\mid l\in\Z_-\}$ (cf. (6.71)),  where the subindex
``b" stands for ``bosonic". Moreover, we denote
$$x_l=\ptl_{\bar x_{-l}},\qquad
\bar x_l=-\ptl_{x_{-l}}\qquad\for\;\;l\in\Z_+.\eqno(6.108)$$ It
can be verified that we the following representation of
$\td{gl}(\infty)$ on $V_b$:
$$\kappa_0|_{V_b}=-\mbox{Id}_{V_b},\;\;\E_{-l,-k}|_{V_b}=-\ptl_{x_l}\ptl_{\bar
x_k},\;\;\E_{l,-k}|_{V_b}=\bar x_l\ptl_{\bar x_k},\eqno(6.109)$$
$$\E_{-l,k}|_{V_b}=-x_k\ptl_{x_l},\qquad \E_{l,k}|_{V_b}=\bar
x_lx_k\eqno(6.110)$$ for $l,k\in\Z_-$. The operator $\ptl$ acts on
$V_b$ by (6.77). Given $i\in\ol{1,n}$, we define
$$x(\iota,i,z)=\sum_{l\in\mbb{Z}}x_{ln-i+1/2}z^{\iota-l},\;\;\bar x(\iota,i,z)
=\sum_{l\in\mbb{Z}}\bar x_{ln+i-1/2}z^{-\iota-l-1}\eqno(6.111)$$
for $\iota\in\mbb{C}\setminus\mbb{Z}$, and
$$x(\iota,i,z)=\sum_{l=0}^\infty[x_{(-l+\iota)n-i+1/2}z^{\ell_i+l}
+x_{(l+\iota+1)n-i+1/2}z^{-l-1}], \eqno(6.112)$$
$$\bar x(\iota,i,z)=\sum_{l=0}^\infty[\bar x_{(l-\iota)n+i-1/2}z^{-\ell_i-l-1}
+\bar x_{-(l+\iota+1)+i-1/2}z^l]\eqno(6.113)$$
 for $\iota\in\mbb{Z}$. Then $\{x(\iota,i,z),\bar x(\iota,i,z)\mid
 i\in\ol{1,n}\}$ are {\it charged free bosonic fields}.

  Set
$$X=\sum_{l\in\Z_-}\mbb{C}x_l,\qquad \bar X=\sum_{l\in\Z_-}\mbb{C}\bar x_l,\eqno(6.114)$$ and
$$V_{b,0}=\mbb{C}+\sum_{s=1}^{\infty}\bar X^sX^s,\;\;V_{b,r}=V_{b,0}X^r,
\;\;V_{b,-r}=\bar X^rV_{b,0}\eqno(6.115)$$ for $r\in\mbb{N}+1$.
Then
$$V_b=\bigoplus_{k\in\mbb{Z}}V_{b,k}.\eqno(6.116)$$
We define a linear map $Y^\iota(\cdot,z):V_b\rta LM(V_b,V_b\{z\})$
by induction:
$$Y^\iota(1,z)=\mbox{Id}_{V_b},\eqno(6.117)$$
\begin{eqnarray*}Y^\iota(x_{-rn-i+1/2}u,z)&=&\rd_{z_1}\sum_{s=0}^\infty(-1)^s\left(\!\!\begin{array}{c}-\iota\\
s\end{array}\!\!\right)z_1^{-\iota-s}z^\iota[(z_1-z)^{s-r-\ell_i-1}x(\iota,i,z_1)Y^\iota(u,z)\\
&
&-(-z+z_1)^{s-r-\ell_i-1}Y^\iota(u,z)x(\iota,i,z_1)]\hspace{4.5cm}(6.118)\end{eqnarray*}
and
\begin{eqnarray*}Y^\iota(\bar x_{-(r+1)n+i-1/2}u,z)&=&\rd_{z_1}\sum_{s=0}^\infty(-1)^s
\left(\!\!\begin{array}{c}\iota\\
s\end{array}\!\!\right)z_1^{\iota-s}z^{-\iota}[(z_1-z)^{s-r-1}\bar
x(\iota,i,z_1)Y^\iota(u,z)
\\ &&-(-z+z_1)^{s-r-1}Y^\iota(u,z)\bar x(\iota,i,z_1)]\hspace{4.4cm}(6.119)\end{eqnarray*}
for $i\in\ol{1,n},\;r\in\mbb{N}$ and $u\in V_{f,k}$. By Section
4.1 in [X2], we have
\begin{eqnarray*}& & z_0^{-1}\dlt\left({z_1-z_2\over z_0}\right)Y^\iota(u,z_1)Y^\iota(v,z_2)
-z_0^{-1}\dlt\left({z_2-z_1\over
-z_0}\right)Y^\iota(v,z_2)Y^\iota(u,z_1)
\\&=& z_2^{-1}\left({z_1-z_0\over z_2}\right)^{k\iota}\dlt\left({z_1-z_0\over z_2}\right)
Y^\iota(Y^0(u,z_0)v,z_2)\hspace{6cm}(6.120)\end{eqnarray*} for
$u\in V_{b,k}$ and $v\in V_b$. In fact,
$$E^\iota_{i,j}(r,z)|_{V_b}=Y^\iota(\bar
x_{-n+i-1/2}x_{-rn-j+1/2},z)\eqno(6.121)$$ (cf. (6.65), (6.66))
$$Y^\iota(x_{-rn+i-1/2},z)=\frac{1}{(r+\ell_i)!}\frac{d^{r+\ell_i}}{z^{r+\ell_i}}
x(\iota,i.z)\eqno(6.122)$$ and
$$Y^\iota(\bar x_{-(r+1)n-i+1/2},z)=\frac{1}{r!}\frac{d^r}{z^r}\bar x(\iota,i.z)\eqno(6.123)$$
for $i\in\ol{1,n}$.

Let $k\in\mbb{N}+1$. Denote
$$v_0=1,\;\;v_k=x_{-1/2}^k,\qquad v_{-k}=\bar x_{-1/2}^k.\eqno(6.124)$$
Moreover, we define linear functions on $\T$ in (3.12):
$$\lmd^0(\kappa_0)=-1,\qquad\lmd^0(\E_{l,-l})=0\qquad\for\;\;l\in\Z,\eqno(6.125)$$
$$\lmd^k(\kappa_0)=-1,\;\;\lmd^k(\E_{1/2,-1/2})=-k,\;\;\lmd^k(E_{l+1/2,-l-1/2})=0\eqno(6.126)$$ for
$0\neq l\in\mbb{Z}$, and
$$\lmd^{-k}(\kappa_0)=-1,\;\;\lmd^{-k}(\E_{-1/2,1/2})=k,\;\;\lmd^{-k}(E_{l-1/2,-l+1/2})=0\eqno(6.127)$$
for $0\neq l\in\mbb{Z}$.  \psp

{\bf Lemma 6.7}. {\it Each space $V_{b,k}$ forms a unitary
irreducible highest weight module of $\td{gl}(\infty)$ with the
highest weight $\lmd^k$.  The family
$(V_{b,0},Y^0(\cdot,z),1,\ptl)$ is a simple vertex operator
algebra isomorphic to
$(V_{-1}(\widehat{gl}(\vec\ell,\mbb{A})),Y(\cdot,z),\vcm,\ptl)$ in
Theorem 6.2. Moreover, each
$(V_{b,k},Y^\iota(|_{V_{b,0}},z)|_{V_{b,k}})$ is an irreducible
$V_{b,0}$-module. The map $Y^\iota(|_{V_{b,l}},z)|_{V_{b,k}}$ is
an intertwining operator of type
$[_{V_{b,l}}^{\;\;\;V_{b,l+k}}\:_{V_{b,k}}]$.} \psp

Now we want to deal with higher level case. Assume
$1<\chi\in\mbb{N}$. Form $\chi$th $\td{gl}(\infty)$-module tensor:
$$V_b^{\la \chi\ra}=V_b\otimes V_b\otimes\cdots\otimes
V_b\qquad(\chi\;\mbox{copies}).\eqno(6.128)$$ Then
$$\kappa_0|_{V_b^{\la \chi\ra}}=\chi\mbox{Id}_{V_f^{\la
\chi\ra}}.\eqno(6.129)$$

Set
$${\cal L}=\sum_{l,k\in\Z_-}\mbb{C}x_l\ptl_{x_k},\qquad
\bar{\cal L}=\sum_{l,k\in\Z_-}\mbb{C}\bar x_l\ptl_{\bar
x_k}.\eqno(6.130)$$ Then ${\cal L}$ and $\bar{\cal L}$ are Lie
subalgebras of $\td{gl}(\infty)|_{V_b}$. Moreover, they are
infinite-dimensional general Lie algebras. Denote
\newpage

$$H=\sum_{l\in\Z_-}x_l\ptl_{x_l},\qquad \bar{H}=\sum_{l\in\Z_-}\bar x_l\ptl_{\bar
x_l}.\eqno(6.131)$$ The subspace $H$ is a toral Cartan subalgebra
of ${\cal L}$ and the subspace $\bar H$ is a toral Cartan
subalgebra of $\bar{\cal L}$. Moreover, $X^k$ is an irreducible
highest weight ${\cal L}$-module with the highest weight
$\lmd^k_X$ determined by
$$\lmd_X^k(x_l\ptl_{x_l})=k\dlt_{l,-1/2}\qquad\for\;\;l\in\Z_-,\eqno(6.132)$$
and $\bar X^k$ is an irreducible highest weight $\bar{\cal
L}$-module with the highest weight $\lmd^k_X$ determined by
$$\lmd_{\bar X}^k(\bar x_l\ptl_{\bar x_l})=k\dlt_{l,-1/2}\qquad\for\;\;l\in\Z_-.\eqno(6.133)$$

For a positive integer $s$, we define:
\begin{eqnarray*}\qquad& &\G_X^s=\{\lmd_X\in H^\ast\mid
\lmd_X(x_{1/2-r}\ptl_{x_{1/2-r}})\in\mbb{N},\\ & &
\lmd_X(x_{-1/2-s-l}\ptl_{x_{-1/2-s-l}})=0
\;\for\;r\in\ol{1,s},\;l\in\mbb{N}\}\hspace{5.4cm}(6.134)\end{eqnarray*}
and \begin{eqnarray*}\qquad& &\G_{\bar X}^s=\{\lmd_{\bar X}\in
\bar H^\ast\mid \lmd_{\bar X}(\bar x_{1/2-r}\ptl_{\bar
x_{1/2-r}})\in\mbb{N},\\ & &\lmd_{\bar X}(\bar
x_{-1/2-s-l}\ptl_{\bar x_{-1/2-s-l}})=0
\;\for\;r\in\ol{1,s},\;l\in\mbb{N}\}.\hspace{5.2cm}(6.135)\end{eqnarray*}
Set
$${\cal C}=\mbb{C}[x_l\mid l\in\Z_-]=\bigoplus_{r=0}^\infty X^r
,\qquad \bar{\cal C}=\mbb{C}[\bar x_l\mid
l\in\Z_-]=\bigoplus_{r=0}^\infty \bar X^r.\eqno(6.136)$$ Then
$$ V_b=\bar CC.\eqno(6.137)$$
By the tensor theory of modules of finite-dimensional general Lie
algebras, the $s$-tensor
$${\cal C}^{\la s\ra}={\cal C}\otimes{\cal
C}\otimes\cdots\otimes{\cal C}\;\;(s\;\mbox{copies})\eqno(6.138)$$
can be decomposed as a direct sum of highest weight irreducible
$\cal L$-submodules, whose set of highest weights is exactly
$\G^s_X$. Similarly, the $s$-tensor
$$\bar{\cal C}^{\la s\ra}=\bar{\cal C}\otimes\bar{\cal
C}\otimes\cdots\otimes{\cal C}\;\;(s\;\mbox{copies})\eqno(6.139)$$
can be decomposed as a direct sum of highest weight irreducible
$\cal L$-submodules, whose set of highest weights is exactly
$\G^s_{\bar X}$.

Since $V_b$ is an unitary $\td{gl}(\infty)$-module, $V_b^{\la
\chi\ra}$ is a completely reducible $\td{gl}(\infty)$-module. Let
$s_1,s_2\in\mbb{N}$ such that $s_1+s_2=\chi$. Suppose that $v_1$
is a highest weight vector of an irreducible component of the
$\bar{\cal L}$-module $\bar{\cal C}^{\la s\ra}$ with $\lmd_{\bar
X}$, and  $v_2$ is a highest weight vector of an irreducible
component of the ${\cal L}$-module ${\cal C}^{\la s\ra}$ with
$\lmd_X$. Then $v_1\otimes v_2$ is the highest weight vector of
some an irreducible component of the $\td{gl}(\infty)$-module
$V_b^{\la \chi\ra}$ (cf. (6.128)), whose weight $\lmd$ is given
by:
$$\lmd(\kappa_0)=-\chi,\;\;
\lmd(\E_{l,-l})=\lmd_{\bar X}(\bar x_l\ptl_{\bar x_l}),\;\;
\lmd(\E_{-l,l})=-\lmd_X(x_l\ptl_{x_l}),\qquad\for\;\;l\in\Z_-.\eqno(6.140)$$
Set
$${\cal S}_\chi=\{\{3/2-r,5/2-r,...,(2\chi+1)/2-r\}\mid
r\in\ol{1,\chi+1}\}.\eqno(6.141)$$ For $\lmd\in\T^\ast$, we define
$$\mbox{supp}\:\lmd=\{l\in\Z\mid \lmd(\E_{l,-l})\neq
0\}.\eqno(6.142)$$ Set
\begin{eqnarray*}\hspace{2cm}\G^\chi&=&\{\lmd\in\T^\ast\mid\lmd(\kappa_0)=-\chi,\;
-s^{-1}|s|\lmd(\E_{s,-s})\in\mbb{N}\;\for\;s\in\Z;\\ & &
\mbox{supp}\:\lmd\subset S\;\mbox{for some}\;S\in{\cal S}_\chi\}.
\hspace{6.5cm}(6.143)\end{eqnarray*} Expression (6.140) implies
that every element in $\G^\chi$ is the highest weight of some
irreducible component of the $\td{gl}(\infty)$-module $V_b^{\la
\chi\ra}$.

Denote
$$1_\chi=1\otimes 1\otimes\cdots\otimes 1\qquad(\chi\;\mbox{copies}).\eqno(6.144)$$
Since $V_b$ is also a polynomial algebra, $V_b^{\la \chi\ra}$ has
an extended commutative associative tenor algebra structure. Note
the space
$$\bar XX=\sum_{l,k\in\Z_-}\mbb{C}\bar x_lx_k.\eqno(6.145)$$
We define the diagonal linear map $\varrho:\bar XX\rta V_b^{\la
\chi\ra}$ by
$$\varrho(\bar x_lx_k)=\sum_{i=1}^\chi 1\otimes\cdots\otimes 1\otimes\stl{i}
{\bar x_lx_k}\otimes 1\otimes\cdots\otimes 1\eqno(6.146)$$ for
$l,k\in\Z_-$. Set
$$V_{b,0}^{[\chi]}=\mbb{C}1_\chi+\sum_{r=1}^{\infty}[\varrho(\bar XX)]^r.\eqno(6.147)$$
Moreover, define the map
$$Y^\iota_\chi(\cdot,z)=Y^\iota(\cdot,z)\otimes Y^\iota(\cdot,z)\otimes\cdots\otimes
Y^\iota(\cdot,z)\qquad(\chi\;\mbox{copies})\eqno(6.148)$$ and
$$\ptl^{\la\chi\ra}=\sum_{i=1}^\chi 1\otimes\cdots\otimes 1\otimes\stl{i}
{\ptl}\otimes 1\otimes\cdots\otimes 1\eqno(6.149)$$ (cf. (6.77)).
Then the family
$(V_{b,0}^{[\chi]},Y^0_\chi(\cdot,z),1_\chi,\ptl_{(\iota)}^{\la\chi\ra})$
forms a simple vertex algebra isomorphic to
$(V_{-\chi}(\widehat{gl}(\vec\ell,\mbb{A})),
Y(\cdot,z),\vcm,\ptl)$ in Theorem 6.2. For each irreducible
$\td{gl}(\infty)$-module component $U$ of $V_f^{\la \chi\ra}$, the
family $(U,Y^\iota_\chi(|_{V_{b,0}^{[\chi]}},z)|_U)$ forms an
irreducible module of
$(V_{b,0}^{[\chi]},Y^0_\chi(\cdot,z),1_\chi,\ptl^{\la\chi\ra})$.
Therefore, we have proved the following main theorem: \psp

{\bf Theorem 6.8}. {\it  Suppose that $\chi$ is a positive integer
and $\lmd\in\G^\chi$. Let $\cal M$ be the irreducible highest
weight $\td{gl}(\infty)$-module with highest weight $\lmd$. Then
the family $({\cal M},Y^\iota_{\cal M})$ defined in (6.65)-(6.70)
forms an irreducible module of the simple vertex algebra
$(V_{-\chi}(\widehat{gl}(\vec\ell,\mbb{A})),Y(\cdot,z),\vcm,\ptl)$
in Theorem 6.2.}

\section{Vacuum Representations of $\hat{o}(\vec\ell,\mbb{A})$ and
$\widehat{sp}(\vec\ell,\mbb{A})$}

In this section, we  study the  vacuum representations of the Lie
algebras  $\hat{o}(\vec\ell,\mbb{A})$ in (3.43) and
$\widehat{sp}(\vec\ell,\mbb{A})$ in (3.57), and their related
vertex algebra structures.  Their vertex algebra irreducible
representations are investigated.

Recall the general settings in (6.1)-(6.9). Note
$$\hat{o}(\vec\ell,\mbb{A})_-=\sum_{i,j=1}^n\:\sum_{r=0}^{\infty}
\sum_{m=1}^{\infty}
\mbb{C}(t^{-m}\ptl_t^{r+\ell_j}E_{i,j}-(-1)^\es(-\ptl)^rt^{-m}\ptl_t^{\ell_i}E_{j^\ast,i^\ast})
+ \mbb{C}\kappa\eqno(7.1)$$  and
$$\hat{o}(\vec\ell,\mbb{A})_+=\sum_{i,j=1}^n\:\sum_{r=0}^{\infty}
\sum_{m=0}^{\infty}
\mbb{C}(t^m\ptl_t^{r+\ell_j}E_{i,j}-(-1)^\es(-\ptl)^rt^m\ptl_t^{\ell_i}E_{j^\ast,i^\ast})
+ \mbb{C}\kappa.\eqno(7.2)$$ The vacuum module $${\cal
V}_\chi(\hat{o}(\vec\ell,\mbb{A}))=U(\hat{o}(\vec\ell,\mbb{A})_-)\vcm\eqno(7.3)$$
and
$$\hat{o}(\vec\ell,\mbb{A})_+\vcm=\{0\},\qquad\kappa(\vcm)=\chi\vcm,\eqno(7.4)$$
where $\chi\in\mbb{C}$.\psp

{\bf Theorem 7.1}. {\it The module ${\cal V}_\chi(\hat o
(\vec\ell,\mbb{A}))$ is irreducible if $\chi\not\in\mbb{Z}$. When
$\chi\in\mbb{Z}$,  the module ${\cal V}(\hat o(\vec\ell,\mbb{A}))$
has a unique maximal proper submodule $\bar{\cal V}_\chi(\hat o
(\vec\ell,\mbb{A}))$, and the quotient
$$V_\chi(\hat o(\vec\ell,\mbb{A}))={\cal
V}(\hat o(\vec\ell,\mbb{A}))/\bar{\cal V}_\chi(\hat o
(\vec\ell,\mbb{A}))\eqno(7.5)$$ is an irreducible $\hat o
(\vec\ell,\mbb{A})$-module. Assume $\chi\in\mbb{N}$. The submodule
$$\bar{\cal V}_\chi(\hat o(\vec\ell,\mbb{A}))=
U(\hat
o(\vec\ell,\mbb{A}))(t^{-1}\ptl^{\ell_1}E_{n,1})^{\chi+1}\vcm\eqno(7.6)$$
if $n>1$ and $\es=1$ (cf. (4.13)). When  $\es=0$ and $n>3$, the
submodule
$$\bar{\cal V}_\chi(\hat o(\vec\ell,\mbb{A}))=
U(\hat
o(\vec\ell,\mbb{A}))(t^{-1}(\ptl^{\ell_1}E_{n-1,1}-\ptl_t^{\ell_2}E_{n,2}))^{\chi+1}
\vcm.\eqno(7.7)$$}

{\it Proof}. For $k\in\mbb{Z}$, we define
$$\hat o(\vec\ell,\mbb{A}))_{(k)}=\sum_{i,j=1}^n\sum_{r=0}^{\infty}\mbb{C}
(t^{r-k}\ptl_t^{r+\ell_j}E_{i,j}-(-1)^\es(-\ptl_t)^rt^{r-k}\ptl_t^{\ell_i}E_{j^\ast,i^\ast})
+\mbb{C}\dlt_{k,0}\kappa\eqno(7.8)$$
\newpage

\noindent Then
$$\hat o(\vec\ell,\mbb{A}))=\bigoplus_{k\in\mbb{Z}}\hat o(\vec\ell,\mbb{A}))_{(k)}
\eqno(7.9)$$ is a $\mbb{Z}$-graded Lie algebra by (2.60), (2.66)
and (2.68) with $\iota=0$. We remark that this grading is not
conformal weight grading. Moreover, we define a $\mbb{Z}$-grading
on ${\cal V}_\chi(\hat{o}(\vec\ell,\mbb{A}))$ by
$${\cal V}_\chi(\hat{o}(\vec\ell,\mbb{A}))_{(0)}=\mbb{C}\vcm,
\;\;{\cal
V}_\chi(\hat{o}(\vec\ell,\mbb{A}))_{(-m)}=\{0\}\qquad\for\;\;m\in\mbb{N}+1\eqno(7.10)$$
and
$${\cal
V}_\chi(\hat{gl}(\vec\ell,\mbb{A}))_{(m)}=\mbox{Span}\{u_1u_2\cdots
u_s\mid u_i\in \hat{o}(\vec\ell,\mbb{A})_-\bigcap
\hat{o}(\vec\ell,\mbb{A})_{(k_i)};\;\sum_{i=1}^sk_i=m\}\eqno(7.11)$$
for $m\in\mbb{N}+1$. Then
$${\cal
V}_\chi(\hat{o}(\vec\ell,\mbb{A}))=\bigoplus_{k\in\mbb{Z}}{\cal
V}_\chi(\hat{o}(\vec\ell,\mbb{A}))_{(k)}\eqno(7.12)$$ is a
$\mbb{Z}$-graded $\hat{o}(\vec\ell,\mbb{A})$-module. observe
\begin{eqnarray*}\qquad& &[\sum_{r,s\in\mbb{N}}\mbb{C}(t^{-s-1}\ptl_t^{\ell_j+r}E_{i,j}
-(-1)^\es(-\ptl_t)^r t^{-s-1} \ptl_t^{\ell_i}E_{j^\ast,i^\ast})]
\bigcap \hat{o}(\vec\ell,\mbb{A}))_{(k)}\\
&=&\sum_{r=0}^{k-1}\mbb{C}
(t^{r-k}\ptl^{r+\ell_j}E_{i,j}-(-1)^\es(-\ptl_t)^rt^{r-k}\ptl_t^{\ell_j}E_{j^\ast,i^\ast}).
\hspace{5.1cm}(7.13)\end{eqnarray*} By [X3],
$$\dim(\sum_{r=0}^{k-1}\mbb{C}(t^{r-k}\ptl^r-(-1)^\es(-\ptl_t)^rt^{r-k}))=k\es+(-1)^\es[|k/2|]
.\eqno(7.14)$$ Thus the character
$$d({\cal
V}_\chi(\hat{o}(\vec\ell,\mbb{A})),q)=\sum_{k=0}^{\infty}(\dim{\cal
V}_\chi(\hat{o}(\vec\ell,\mbb{A}))_{(k)})q^k=
\prod_{r=1}^{\infty}\frac{1}{(1-q^r)^{n(r(n-1)/2+r\es+(-1)^\es[|r/2|])}}.\eqno(7.15)$$

Recall the Lie algebra $\td{gl}(\infty)$ defined in (3.9) and
(3.10), and the Lie algebra ${\cal L}^{\ast,\vec\ell}_0$ in
(5.12), which is a Lie subalgebra of $\td{gl}(\infty)$ by (5.11).
Next we set
\begin{eqnarray*}\qquad{\cal L}^{\ast,\vec\ell}_{(-)}&=&
\sum_{i,j=1}^n\:\sum_{l,k=0}^{\infty}[\mbb{C}
((-1)^\es\la-k-1\ra_{\ell_j} \E_{(-l-1)n+i-1/2,-kn-j+1/2}\\ &
&-\la
-l-1\ra_{\ell_i}\E_{-kn-j+1/2,(-l-1)n+i-1/2})\hspace{6.4cm}(7.16)\end{eqnarray*}
and
\begin{eqnarray*}{\cal L}^{\ast,\vec\ell}_{(+)}&=&\sum_{i,j=1}^n\:\sum_{l,k=0}^{\infty}
[\mbb{C}((-1)^\es\la-k-1\ra_{\ell_j}\E_{ln+i-1/2,-kn-j+1/2}\\
& & -\la
l+\ell_i\ra_{\ell_i}\E_{-kn-j+1/2,ln+i-1/2})+\mbb{C}((-1)^\es\la
k+\ell_j\ra_{\ell_j}\E_{ln+i-1/2,(k+1)n-j+1/2}\\ & &-\la
l+\ell_i\ra_{\ell_i}\E_{(k+1)n-j+1/2,ln+i-1/2})].\hspace{7.6cm}(7.17)
\end{eqnarray*}
By (3.10),  ${\cal L}^{\ast,\vec\ell}_{(\pm)}$ are Lie subalgebras
of ${\cal L}^{\ast,\vec\ell}_0$ and
$${\cal L}^{\ast,\vec\ell}_0={\cal L}^{\ast,\vec\ell}_{\iota,(-)}+
{\cal L}^{\ast,\vec\ell}_{(+)}+\mbb{C}\kappa_0.\eqno(7.18)$$
Recall the $\td{gl}(\infty)$-module $U_\chi$ defined in
(6.20)-(6.24). Note
$${\cal L}^{\ast,\vec\ell}_{(+)}{\bf 1}=\{0\} \eqno(7.19)$$
because ${\cal
L}^{\ast,\vec\ell}_{(+)}\subset\td{gl}(\infty)_{(+)}$. Thus
$$U^o_\chi=U({\cal L}^{\ast,\vec\ell}_{(-)}){\bf 1}\eqno(7.20)$$
is a ${\cal L}^{\ast,\vec\ell}_0$-module,
 which satisfies the condition (5.45) with
$s=\bar\ell+1$ (cf. (5.42)) and ${\cal M}_0=\mbb{C}1\otimes{\bf
1}$.

Note our notion (5.46). We have the
$\hat{o}(\vec\ell,\mbb{A})$-module structure  on $U^o_\chi$
defined by $\kappa=\chi\mbox{Id}_{U^o_\chi}$ and
\begin{eqnarray*}(E_{i,j})^\ast_{\vec\ell}(r,z)
&=&\sum_{l,k=0}^{\infty}[\la k\ra_r ((-1)^\es\la-k-1\ra_{\ell_j}
\E_{(-l-1)n+i-1/2,-kn-j+1/2}\\ & &-\la
-l-1\ra_{\ell_i}\E_{-kn-j+1/2,(-l-1)n+i-1/2})z^{l+k-r} \\ & &+\la
k\ra_r ((-1)^\es\la-k-1\ra_{\ell_j}\E_{ln+i-1/2,-kn-j+1/2}\\
&&-\la
l+\ell_i\ra_{\ell_i}\E_{-kn-j+1/2,ln+i-1/2})z^{-l+k-\ell_i-r-1}\\
&&+\la-k-\ell_2-1\ra_r((-1)^\es\la
k+\ell_j\ra_{\ell_j}\E_{-(l+1)n+i-1/2,(k+1)n-j+1/2}\\ & &-\la
-l-1\ra_{\ell_i}\E_{(k+1)n-j+1/2,-(l+1)n+i-1/2})z^{l-k-\ell_j-r-1}\\
& &+\la -k-\ell_j-1\ra_r((-1)^\es\la
k+\ell_j\ra_{\ell_j}\E_{ln+i-1/2,(k+1)n-j+1/2}\\ & &-\la
l+\ell_i\ra_{\ell_i}\E_{(k+1)n-j+1/2,ln+i-1/2})z^{-l-k-\ell_i-\ell_j-r-2}
\hspace{3.9cm}(7.21)\end{eqnarray*} for $i,j\in\ol{1,n}$ and
$r\in\mbb{N}$. In particular,
\begin{eqnarray*}\qquad(E_{i,j})^\ast_{\vec\ell}(r,z){\bf 1}
&=&\sum_{l,k=0}^{\infty}\la k\ra_r ((-1)^\es\la-k-1\ra_{\ell_j}
\E_{(-l-1)n+i-1/2,-kn-j+1/2}\\ & &-\la
-l-1\ra_{\ell_i}\E_{-kn-j+1/2,(-l-1)n+i-1/2})z^{l+k-r}{\bf 1}.
\hspace{3.5cm}(7.22)\end{eqnarray*}
 Thus we have
$$\hat o(\ell,\mbb{A})_+({\bf 1})=\{0\}.\eqno(7.23)$$ By a
similar proof as that of Theorem 3.1,
$$U^o_\chi=U(\hat
o(\vec\ell,\mbb{A})_-){\bf 1}.\eqno(7.24)$$ Therefore, we have a
Lie algebra module epimorphism $\nu:{\cal V}_\chi(\hat
o(\vec\ell,\mbb{A}))\rta U^o_\chi$ defined by
$$\nu(u\vcm)=u{\bf 1}\qquad\for\;\;u\in U(\hat
o(\vec\ell,\mbb{A})_-).\eqno(7.25)$$

For $m\in\mbb{N}+1$, we let \begin{eqnarray*}& &{\cal
L}^{\ast,\vec\ell}_{(-),m}=\mbox{Span}\:\{(-1)^\es\la-k-1\ra_{\ell_j}
\E_{(-l-1)n+i-1/2,-kn-j+1/2}\\ & &-\la
-l-1\ra_{\ell_i}\E_{-kn-j+1/2,(-l-1)n+i-1/2})\mid
i,j\in\ol{1,n},\;l,k\in\mbb{N};\;l+k+1=m\}.
\hspace{1.7cm}(7.26)\end{eqnarray*} Then
$${\cal L}^{\ast,\vec\ell}_{(-)}=\bigoplus_{m=1}^{\infty}{\cal
L}^{\ast,\vec\ell}_{(-),m}.\eqno(7.27)$$ Moreover, we define
$$U_\chi^{o,(0)}=\mbb{C}{\bf
1},\;\;U_\chi^{o,(m)}=\{0\}\qquad\for\;\;m\in(-\mbb{N}-1)\eqno(7.28)$$
and
$$U_\chi^{o,(m)}=\mbox{Span}\:\{u_1u_2\cdots u_s{\bf 1}\mid u_i\in
{\cal
L}^{\ast,\vec\ell}_{(-),m_i};\;\sum_{i=1}^sm_i=m\}.\eqno(7.29)$$
Expressions (7.8) (7.21) and (7.24) imply
$$ U^o_\chi=\bigoplus_{m\in\mbb{Z}}U_\chi^{o,(m)}\eqno(7.31)$$
is a $\mbb{Z}$-graded $\hat{o}(\vec\ell,\mbb{A})$-module.
Furthermore, (7.26) gives the character
$$d(U^o_\chi,q)=\sum_{m=0}^\infty(\dim U_\chi^{o,(m)})z^m=
\prod_{r=1}^{\infty}\frac{1}{(1-q^r)^{n(r(n-1)/2+r\es+(-1)^\es[|r/2|])}}.
\eqno(7.32)$$ Therefore, (7.15) and (7.32) yield
$${\cal V}_\chi(\hat{o}(\vec\ell,\mbb{A}))\cong U^o_\chi.\eqno(7.33)$$

Let $\lmd$ be a linear function on $\T$ defined in (5.22) such
that
$$\lmd(\kappa_0)=\chi,\;\;\lmd(\vt_l)=0\qquad\for\;\;l\in\mbb{N}\eqno(7.34)$$
(cf. (5.21)). Recall the Verma module $M_\lmd$ defined in (5.39)
with $\iota=0$ and $\tau=\ast$. Note
$$U^o_\chi\cong M_\lmd/(\sum_{l=1}^\infty U({\cal L}^{\ast,\vec
\ell}_{0,-})f_{\es,l}^\ast\otimes v_\lmd)),\eqno(7.35)$$ which is
irreducible if $\chi\not\in\mbb{Z}$ by [J1-J3]. When
$\chi\in\mbb{N}$,
$$\bar{U}^o_\chi=U({\cal L}^{\ast,\vec\ell}_{0,-})(f_{\es,0}^\ast)^{\chi+1}{\bf 1}\eqno(7.36)$$
 is the unique maximal
proper submodule of $U_\chi^o$. Thus
$$U(\hat{o}(\vec\ell,\mbb{A}))\nu^{-1}((f_{\es,0}^\ast)^{\chi+1}{\bf
1})\eqno(7.37)$$ is the unique maximal proper submodule of ${\cal
V}_\chi(\hat{o}(\vec\ell,\mbb{A}))$. When $n>1$, (5.28), (7.21)
and (7.25) imply
$$\nu^{-1}((f_{1,0}^\ast)^{\chi+1}{\bf
1})=(t^{-1}\ptl^{\ell_1}E_{n,1})^{\chi+1}\vcm.\eqno(7.38)$$
Moreover, if $n>3$, (5.26), (7.21) and (7.25) yield
$$\nu^{-1}((f_{0,0}^\ast)^{\chi+1}{\bf
1})=(t^{-1}(\ptl^{\ell_1}E_{n-1,1}-\ptl_t^{\ell_2}E_{n,2}))^{\chi+1}\vcm.
\qquad\Box\eqno(7.39)$$ \vspace{0.1cm}

Since $\hat o(\ell,\mbb{A})$ is a Lie subalgebra of
$\widehat{gl}(\vec\ell,\mbb{A})$, we  view ${\cal
V}_\chi(\hat{o}(\vec\ell,\mbb{A}))$ as a subspace of ${\cal
V}_\chi(\widehat{gl}(\vec\ell,\mbb{A}))$. Recall the vertex
algebra $({\cal
V}_\chi(\widehat{gl}(\vec\ell,\mbb{A})),Y(\cdot,z),\ptl,\vcm)$
defined in (6.57)-(6.63). Then
$$({\cal
V}_\chi(\hat{o}(\vec\ell,\mbb{A})),Y(|_{{\cal
V}_\chi(\hat{o}(\vec\ell,\mbb{A}))},z)|_{{\cal
V}_\chi(\hat{o}(\vec\ell,\mbb{A}))},\ptl|_{{\cal
V}_\chi(\hat{o}(\vec\ell,\mbb{A}))},\vcm)\eqno(7.40)$$ forms a
vertex subalgebra. Let ${\cal M}$ be a weighted irreducible
$\td{gl}(\infty)$-module satisfying (3.18) (also cf. (3.17)) and
$\kappa_0|_{\cal M}=\chi\mbox{Id}_{\cal M}$. Recall  the linear
map $Y^{\iota}_{\cal M}(\cdot,z)$ defined by (6.65)-(6.70). Now we
obtain:\psp

{\bf Theorem 7.2}. {\it The family (7.40) forms a vertex algebra
and $({\cal M} ,Y^\iota_{\cal M}(|_{{\cal
V}_\chi(\hat{o}(\vec\ell,\mbb{A}))},z))$ forms a irreducible
vertex algebra module of the vertex algebra (7.40). If
$\chi\not\in \mbb{Z}$, the vertex algebra $({\cal
V}_\chi(\hat{o}(\vec\ell,\mbb{A})),Y(\cdot,z),\ptl,\vcm)$ is
simple. When $\chi\in\mbb{Z}$, the quotient space
$V_\chi(\hat{o}(\vec\ell,\mbb{A}))$ forms a simple vertex algebra.
If $\ell_i=\es$ for $i\in\ol{1,n}$, ${\cal
V}_\chi(\hat{o}(\vec\ell,\mbb{A}))$ with
$\chi\in\mbb{C}\setminus\mbb{Z}$ and
$V_\chi(\hat{o}(\vec\ell,\mbb{A}))$ with $\chi\in\mbb{Z}$ are
simple vertex operator algebras with the Virasoro element}
$$\sum_{i=1}^{[|(n+1)/2|]}(-t^{-1}\ptl_tE_{i,i}+(-1)^\es
(-\ptl_t)^{\dlt_{\es,0}}t^{-1}\ptl_t^\es
E_{i^\ast,i^\ast})\vcm.\eqno(7.41)$$ \vspace{0.1cm}

Assume that $\chi$ is a positive integer. Recall the assumption
(3.39) and the  charged free fermionic field realization given in
(6.71)-(6.86). We set
$$
R^o=\mbox{Span}\;\{\bar\cta_l\cta_k-(-1)^\es\bar\cta_k\cta_l\mid
l,k\in\Z_-\}.\eqno(7.42)$$ Note the notion $V_f^{\la \chi\ra}$
defined in (6.95), the notion $1_\chi$ defined in (6.97) and the
map $\varrho:\bar\Theta\Theta\rta V_f^{\la \chi\ra}$ defined by
(6.81) and (6.99). Set
$$V^{o,f}_\chi=\mbb{C}1_\chi+\sum_{r=1}^{\infty}[\varrho(R^o)]^r.\eqno(7.43)$$
In terms of (6.101) and (6.102), the family
$(V^{o,f}_\chi,Y^0_\chi(\cdot,z),1_\chi,\ptl^{\la\chi\ra})$ forms
a simple vertex algebra isomorphic to
$(V_\chi(\hat{o}(\vec\ell,\mbb{A})), Y(\cdot,z),\vcm,\ptl)$ in
Theorem 7.2. For each irreducible $\td{gl}(\infty)$-module
component $U$ of $V_f^{\la \chi\ra}$, the family
$(U,Y^\iota_\chi(|_{V^{o,f}_\chi},z)|_U)$ forms an irreducible
module of
$(V^{o,f}_\chi,Y^0_\chi(\cdot,z),1_\chi,\ptl^{\la\chi\ra})$ when
$\iota\not\in\mbb{Z}/2$.

Recall the  charged free bosonic field realization given in
(6.107)-(6.119). We set
$$
R_o=\mbox{Span}\;\{\bar x_lx_k-(-1)^\es\bar x_kx_l\mid
l,k\in\Z_-\}.\eqno(7.44)$$ Note the notion $V_b^{\la \chi\ra}$
defined in (6.128), the notion $1_\chi$ defined in (6.144) and the
map $\varrho:\bar\Theta\Theta\rta V_b^{\la \chi\ra}$ defined by
(6.114) and (6.146). Set
$$V^{o,b}_\chi=\mbb{C}1_\chi+\sum_{r=1}^{\infty}[\varrho(R_o)]^r.\eqno(7.45)$$
According to(6.148) and (6.149), the family
$(V^{o,b}_\chi,Y^0_\chi(\cdot,z),1_\chi,\ptl^{\la\chi\ra})$ forms
a simple vertex algebra isomorphic to
$(V_{-\chi}(\hat{o}(\vec\ell,\mbb{A})), Y(\cdot,z),\vcm,\ptl)$ in
Theorem 7.2. For each irreducible $\td{gl}(\infty)$-module
component $U$ of $V_b^{\la \chi\ra}$, the family
$(U,Y^\iota_\chi(|_{V^{o,b}_\chi},z)|_U)$ forms an irreducible
module of $(V^{o,b}_\chi,Y^0_\chi(\cdot,z),
1_\chi,\ptl^{\la\chi\ra})$ when $\iota\not\in\mbb{Z}/2$. Thus we
have:\psp

{\bf Theorem 7.3}. {\it  Suppose that $\chi$ is a positive
integer. Assume $\lmd$ is a weight of $\td{gl}(\infty)$ satisfying
(3.59), (6.105), (6.106) and $\lmd(\kappa_0)=\chi$. Let $\cal M$
be the  highest weight irreducible  $\td{gl}(\infty)$-module with
highest weight $\lmd$. Then the family $({\cal M},Y^\iota_{\cal
M}(\cdot,z))$ defined in (6.65)-(6.70) forms an irreducible module
of the simple vertex algebra
$(V_\chi(\hat{o}(\vec\ell,\mbb{A})),Y(\cdot,z),\vcm,\ptl)$ in
Theorem 7.2 when $\iota\not\in\mbb{Z}/2$.

If $\lmd \in \G^\chi$ (cf. (6.143)) and $\cal M$ is the
irreducible highest weight $\td{gl}(\infty)$-module with highest
weight $\lmd$, then the family $({\cal M},Y^\iota_{\cal
M}(\cdot,z))$ defined in (6.65)-(6.70) forms an irreducible module
of the simple vertex algebra
$(V_{-\chi}(\hat{o}(\vec\ell,\mbb{A})),Y(\cdot,z),\vcm,\ptl)$ in
Theorem 7.2 when $\iota\not\in\mbb{Z}/2$. } \psp

Assume $\iota\in \mbb{Z}+1/2$. Recall the Lie algebra ${\cal
L}^{\ast,\vec m}_{\iota_0,\es}$ defined in (4.30) and the highest
weight irreducible ${\cal L}^{\ast,\vec m}_{\iota_0,\es}$-module
${\cal M}^{\ast,\es}_\lmd$ defined in (4.68) with $\tau=\ast$ and
$\lmd(\kappa_0)=\chi$. By (4.75), we define operators
\begin{eqnarray*}& &E_{i,j}^{\iota,\ast}(r,z)\\
 &=&\sum_{l,k\in\mbb{Z}}\la k-m_j-\es+1/2\ra_r((-1)^\es\la
k+m_j+\es-3/2\ra_{\ell_j}\E_{ln+i-1/2,kn-j+1/2}\\ & &-\la
l+m_i+1/2\ra_{\ell_i}\E_{(k+\es-1)n-j+1/2,(l+1-\es)n+i-1/2})z^{-l-k-m_i-m_j-r-\es-1}\\
&&+((r+\ell_i)!\Im_{0,r+\ell_i}-r!\ell_i!\Im_{r,\ell_i})
\dlt_{i,j}\kappa_0z^{-\ell_i-r-1}\hspace{6.6cm}(7.46)\end{eqnarray*}
on ${\cal M}^{\ast,\es}_\lmd$ for $i,j\in\ol{1,n}$ and
$r\in\mbb{N}$. For $l\in\mbb{Z}$, we define
$$\vf(l)=\left\{\begin{array}{ll}-1&\mbox{if}\;l-1/2+(m_{l_R}-\iota_0)n>0,\\
1&\mbox{if}\;l-1/2+(m_{l_R}-\iota_0)n<0\end{array}\right.\eqno(7.47)$$
(cf. (4.13) and (4.25)) by (4.16). Moreover, for $\lmd\in
(\T^\es)^\ast$ (cf. (4.47) and (4.59)), we define
$$\mbox{supp}\:\lmd=\{l-1/2+(m_{l_R}-\iota_0)n\mid
l\in\mbb{N}-n_0\dlt_{\es,0},\;\lmd(\vt^\es_l)\neq
0\}.\eqno(7.48)$$ Now we let
\begin{eqnarray*}\hspace{2cm}\G^\chi_{\iota,\es}&=&\{\lmd\in (\T^\es)^\ast\mid\lmd(\kappa_0)=-\chi,\;
\vf(l)\lmd(\vt^\es_l)\in\mbb{N}\;\for\;l\in\mbb{N}-n_0\dlt_{\es,0};\\
& & \mbox{supp}\:\lmd\subset S\;\mbox{for some}\;S\in{\cal
S}_\chi\} \hspace{6.6cm}(7.49)\end{eqnarray*}(cf. (6.141)).

Let $\iota\in\mbb{Z}$.  Recall the Lie algebra ${\cal
L}^{\ast,\vec\ell}_{\iota}$ defined in (5.12) and the highest
weight irreducible ${\cal L}^{\ast,\vec\ell}_{\iota}$-module
${\cal M}_{\lmd^{\ast,\es}}$ defined in (5.39) with $\tau=\ast$
and $\lmd^{\ast,\es}(\kappa_0)=\chi$. By (5.48), we define
operators
\newpage

\begin{eqnarray*}E_{i,j}^{\iota,\ast}(r,z)
 &=&\sum_{l,k=0}^{\infty}[\la k\ra_r ((-1)^\es\la-k-1\ra_{\ell_j}
\E_{(-l-1)n+i-1/2,-kn-j+1/2}\\& &-\la
-l-1\ra_{\ell_i}\E_{-kn-j+1/2,(-l-1)n+i-1/2})z^{l+k-r} \\ & &+\la
k\ra_r ((-1)^\es\la-k-1\ra_{\ell_j}\E_{ln+i-1/2,-kn-j+1/2}\\
&&-\la
l+\ell_i\ra_{\ell_i}\E_{-kn-j+1/2,ln+i-1/2})z^{-l+k-\ell_i-r-1}\\
&&+\la-k-\ell_2-1\ra_r((-1)^\es\la
k+\ell_j\ra_{\ell_j}\E_{-(l+1)n+i-1/2,(k+1)n-j+1/2}\\ & &-\la
-l-1\ra_{\ell_i}\E_{(k+1)n-j+1/2,-(l+1)n+i-1/2})z^{l-k-\ell_j-r-1}\\
& &+\la -k-\ell_j-1\ra_r((-1)^\es\la
k+\ell_j\ra_{\ell_j}\E_{ln+i-1/2,(k+1)n-j+1/2}\\ & &-\la
l+\ell_i\ra_{\ell_i}\E_{(k+1)n-j+1/2,ln+i-1/2})z^{-l-k-\ell_i-\ell_j-r-2}\\
& & + ((r+\ell_i)!\Im_{0,r+\ell_i} -(-1)^\es
r!\ell_i!\Im_{r,\ell_i})\dlt_{i,j}\kappa_0z^{-r-\ell_i-1}\hspace{4.3cm}(7.50)\end{eqnarray*}
on ${\cal M}_{\lmd^{\ast,\es}}$ for $i,j\in\ol{1,n}$ and
$r\in\mbb{N}$. For $l\in\mbb{N}$, we define
$$\bar\vf(l)=\left\{\begin{array}{ll}-1&\mbox{if}\;l-1/2-\iota>0,\\
1&\mbox{if}\;l-1/2-\iota<0\end{array}\right.\eqno(7.51)$$ by
(5.3). Moreover, for $\lmd\in \T^\ast$ (cf. (5.21) and (5.22)), we
define
$$\mbox{supp}\:\lmd=\{l-1/2-\iota\mid
l\in\mbb{N},\;\lmd(\vt_l)\neq 0\}.\eqno(7.52)$$ Now we let
\begin{eqnarray*}\hspace{2cm}\G^\chi_{\iota,\es}&=&\{\lmd\in \T^\ast\mid\lmd(\kappa_0)=-\chi,
\;
\bar\vf(l)\lmd(\vt^\es_l)\in\mbb{N}\;\for\;l\in\mbb{N};\\
& & \mbox{supp}\:\lmd\subset S\;\mbox{for some}\;S\in{\cal
S}_\chi\} \hspace{6.6cm}(7.53)\end{eqnarray*} (cf. (6.141)).

For convenience, we denote
$${\cal M}=\left\{\begin{array}{ll}{\cal M}^{\ast,\es}_\lmd\;\mbox{in
(4.68)}&\mbox{if}\;\iota\in\mbb{Z}+1/2,\\ {\cal
M}_{\lmd^{\ast,\es}}\;\mbox{in
(5.39)}&\mbox{if}\;\iota\in\mbb{Z}.\end{array}\right.\eqno(7.54)$$
We define linear maps
$$Y^{\iota,\pm}_{\cal M}(\cdot,z):\hat{o}(\vec\ell,\mbb{A})_-\rta LM({\cal
M},{\cal M}[z^{-1},z]])\eqno(7.55)$$ by
$$Y^{\iota,\pm}_{\cal M}(t^{-m-1}\ptl_t^{r+\ell_j}E_{i,j}-(-1)^\es(-\ptl_t)^r
t^{-m-1}\ptl_t^{\ell_i}E_{j^\ast,i^\ast},z)=\frac{1}{m!}\frac{d^m}{dz^m}
E^{\iota,\ast}_{i,j}(r,z)^{\pm}\eqno(7.56)$$ for $i,j\in\ol{1,n}$
and $r,m\in\mbb{N}$. Now we define a linear map
$$Y^\iota_{\cal M}(\cdot,z):{\cal V}_\chi(\hat{o}(\vec\ell,\mbb{A}))\rta LM({\cal M},{\cal M}
[z^{-1},z]])\eqno(7.57)$$ by induction:
$$ Y^\iota_{\cal M}(\vcm,z)=\mbox{Id}_{{\cal
M}},\;\;Y(uv,z)=Y^{\iota,-}_{\cal M}(u,z)Y_{\cal M}
^\iota(v,z)+Y^\iota_{\cal M}(v,z)Y^{\iota,+}_{\cal
M}(u,z)\eqno(7.58)$$ for $u\in \hat{o}(\vec\ell,\mbb{A})_-$ and
$v\in {\cal V}_\chi(\hat{o}(\vec\ell,\mbb{A})).$
\newpage

By Theorem 4.2, Theorem 5.2, the general theory for vertex
algebras (e.g. cf. Section 4.1 in [X2]), the  charged free
fermionic field realization and the charged free bosonic field
realization, we obtain:

{\bf Theorem 7.4}. {\it Assume $\iota\in\mbb{Z}/2$. The family
$({\cal M} ,Y^\iota_{\cal M}(\cdot,z))$ forms an irreducible
module of the vertex algebra $({\cal
V}_\chi(\hat{o}(\vec\ell,\mbb{A})),Y(\cdot,z),\ptl,\vcm)$. Suppose
that $\chi$ is a positive integer. If (4.69) and (5.40) hold, then
the family $({\cal M} ,Y^\iota_{\cal M}(\cdot,z))$ induces an
irreducible module of the quotient simple vertex algebra
$(V_\chi(\hat{o}(\vec\ell,\mbb{A})),Y(\cdot,z),\vcm,\ptl)$ in
Theorem 7.2. When $\lmd\in\G^\chi_{\iota,\es}$ with
$\iota\in\mbb{Z}+1/2$ in (7.49) and $\lmd^{\ast,\es}\in
\in\G^\chi_{\iota,\es}$ with $\iota\in\mbb{Z}$ in (7.53), the
family $({\cal M} ,Y^\iota_{\cal M}(\cdot,z))$ induces an
irreducible module of the quotient simple vertex algebra
$(V_{-\chi}(\hat{o}(\vec\ell,\mbb{A})),Y(\cdot,z),\vcm,\ptl)$ in
Theorem 7.2.}\psp

Recall the general settings in (6.1)-(6.8). Observe
$$\widehat{sp}(\vec\ell,\mbb{A})_-=\sum_{i,j=1}^n\:\sum_{r=0}^{\infty}
\sum_{m=1}^{\infty}
\mbb{C}(t^{-m}\ptl_t^{r+\ell_j}E_{i,j}-(-1)^{p(i)+p(j)+\es}(-\ptl)^rt^{-m}\ptl_t^{\ell_i}E_{j^\ast,i^\ast})
+ \mbb{C}\kappa\eqno(7.59)$$  and
$$\widehat{sp}(\vec\ell,\mbb{A})_+=\sum_{i,j=1}^n\:\sum_{r=0}^{\infty}
\sum_{m=0}^{\infty}
\mbb{C}(t^m\ptl_t^{r+\ell_j}E_{i,j}-(-1)^{p(i)+p(j)+\es}(-\ptl)^rt^m\ptl_t^{\ell_i}E_{j^\ast,i^\ast})
+ \mbb{C}\kappa\eqno(7.60)$$ (cf. (3.54)-(3.57)). The vacuum
module $${\cal
V}_\chi(\widehat{sp}(\vec\ell,\mbb{A}))=U(\widehat{sp}(\vec\ell,\mbb{A})_-)\vcm\eqno(7.61)$$
and
$$\widehat{sp}(\vec\ell,\mbb{A})_+\vcm=\{0\},\qquad\kappa(\vcm)=\chi\vcm,\eqno(7.62)$$
where $\chi\in\mbb{C}$. By similar proof as that of Theorem 7.1,
we have:\psp

{\bf Theorem 7.5}. {\it The module ${\cal V}_\chi(\widehat{sp}
(\vec\ell,\mbb{A}))$ is irreducible if $\chi\not\in\mbb{Z}$. When
$\chi\in\mbb{Z}$,  the module ${\cal
V}(\widehat{sp}(\vec\ell,\mbb{A}))$ has a unique maximal proper
submodule $\bar{\cal V}_\chi(\widehat{sp} (\vec\ell,\mbb{A}))$,
and the quotient
$$V_\chi(\widehat{sp}(\vec\ell,\mbb{A}))={\cal
V}(\widehat{sp}(\vec\ell,\mbb{A}))/\bar{\cal V}_\chi(\widehat{sp}
(\vec\ell,\mbb{A}))\eqno(7.63)$$ is an irreducible $\widehat{sp}
(\vec\ell,\mbb{A})$-module. Assume $\chi\in\mbb{N}$. The submodule
$$\bar{\cal V}_\chi(\widehat{sp}(\vec\ell,\mbb{A}))=
U(\widehat{sp}(\vec\ell,\mbb{A}))(t^{-1}\ptl^{\ell_1}E_{n,1})^{\chi+1}\vcm\eqno(7.64)$$
if $n>1$ and $\es=0$ (cf. (4.13)). When  if $\es=1$ and $n>3$,
$$\bar{\cal V}_\chi(\widehat{sp}(\vec\ell,\mbb{A}))=
U(\widehat{sp}(\vec\ell,\mbb{A}))(t^{-1}(\ptl^{\ell_1}E_{n-1,1}-\ptl_t^{\ell_2}E_{n,2}))^{\chi+1}
\vcm.\eqno(7.65)$$ }\psp

Since $\widehat{sp}(\ell,\mbb{A})$ is a Lie subalgebra of
$\widehat{gl}(\vec\ell,\mbb{A})$, we  view ${\cal
V}_\chi(\widehat{sp}(\vec\ell,\mbb{A}))$ as a subspace of ${\cal
V}_\chi(\widehat{gl}(\vec\ell,\mbb{A}))$. Recall the vertex
algebra $({\cal
V}_\chi(\widehat{gl}(\vec\ell,\mbb{A})),Y(\cdot,z),\ptl,\vcm)$
defined in (6.57)-(6.63). The family
$$({\cal
V}_\chi(\widehat{sp}(\vec\ell,\mbb{A})),Y(|_{{\cal
V}_\chi(\widehat{sp}(\vec\ell,\mbb{A}))},z)|_{{\cal
V}_\chi(\widehat{sp}(\vec\ell,\mbb{A}))},\ptl|_{{\cal
V}_\chi(\widehat{sp}(\vec\ell,\mbb{A}))},\vcm)\eqno(7.66)$$
\newpage

\noindent forms a vertex subalgebra. Let ${\cal M}$ be a weighted
irreducible $\td{gl}(\infty)$-module satisfying (3.18) (also cf.
(3.17)) and $\kappa_0|_{\cal M}=\chi\mbox{Id}_{\cal M}$. Recall
the operator the linear map $Y^{\iota}_{\cal M}(\cdot,z)$ defined
by (6.65)-(6.70). Now we obtain \psp

{\bf Theorem 7.6}. {\it The family (7.66) forms a vertex algebra
and $({\cal M} ,Y^\iota_{\cal M}(|_{{\cal
V}_\chi(\widehat{sp}(\vec\ell,\mbb{A}))},z))$ with
$\iota\not\in\mbb{Z}/2$ forms an irreducible vertex algebra module
of the vertex algebra (7.66). If $\chi\not\in \mbb{Z}$, the vertex
algebra $({\cal
V}_\chi(\widehat{sp}(\vec\ell,\mbb{A})),Y(\cdot,z),\ptl,\vcm)$ is
simple. When $\chi\in\mbb{Z}$, the quotient space
$V_\chi(\widehat{sp}(\vec\ell,\mbb{A}))$ forms a simple vertex
algebra. If $\ell_i=\es$ for $i\in\ol{1,n}$, ${\cal
V}_\chi(\widehat{sp}(\vec\ell,\mbb{A}))$ with
$\chi\in\mbb{C}\setminus\mbb{Z}$ and
$V_\chi(\widehat{sp}(\vec\ell,\mbb{A}))$ with $\chi\in\mbb{Z}$ are
simple vertex operator algebras with the Virasoro element
$$\sum_{i=1}^{[|(n+1)/2|]}(-t^{-1}\ptl_tE_{i,i}+(-1)^\es
(-\ptl_t)^{\dlt_{\es,0}}t^{-1}\ptl_t^\es
E_{i^\ast,i^\ast})\vcm.\eqno(7.67)$$

Suppose that $\chi$ is a positive integer. Assume $\lmd$ is a
weight of $\td{gl}(\infty)$ satisfying (3.59), (6.105), (6.106)
and $\lmd(\kappa_0)=\chi$. Let $\cal M$ be the irreducible highest
weight $\td{gl}(\infty)$-module with highest weight $\lmd$. Then
the family $({\cal M},Y^\iota_{\cal M}(\cdot,z))$ induces an
irreducible module of the simple vertex algebra
$(V_\chi(\widehat{sp}(\vec\ell,\mbb{A})),Y(\cdot,z),\vcm,\ptl)$
when $\iota\not\in\mbb{Z}/2$.

If $\lmd \in \G^\chi$ (cf. (6.143)) and $\cal M$ is the
 highest weight irreducible $\td{gl}(\infty)$-module with highest
weight $\lmd$, then the family $({\cal M},Y^\iota_{\cal
M}(\cdot,z))$ induces an irreducible module of the simple vertex
algebra
$(V_{-\chi}(\widehat{sp}(\vec\ell,\mbb{A})),Y(\cdot,z),\vcm,\ptl)$
 when $\iota\not\in\mbb{Z}/2$. } \psp

Assume $\iota\in \mbb{Z}+1/2$. Recall the Lie algebra ${\cal
L}^{\dg,\vec m}_{\iota_0,\es}$ defined in (4.31) and the highest
weight irreducible ${\cal L}^{\dg,\vec m}_{\iota_0,\es}$-module
${\cal M}^{\dg,\es}_\lmd$ defined in (4.68) with $\tau=\dg$ and
$\lmd(\kappa_0)=\chi$. By (4.76), we define operators
\begin{eqnarray*}& &E_{i,j}^{\iota,\dg}(r,z)\\
 &=&\sum_{l,k\in\mbb{Z}}\la l-m_j-\es+1/2\ra_r((-1)^\es\la
k+m_j+\es-3/2\ra_{\ell_j}\E_{ln+i-1/2,kn-j+1/2}\\ &
&-(-1)^{p(i)+p(j)}\la
l+m_i+1/2\ra_{\ell_i}\E_{(k+\es-1)n-j+1/2,(l+1-\es)n+i-1/2})z^{-l-k-m_i-m_j-r-\es-1}\\
&&+((r+\ell_i)!\Im_{0,r+\ell_i}-r!\ell_i!\Im_{r,\ell_i})
\dlt_{i,j}\kappa_0z^{-\ell_i-r-1}\hspace{6.8cm}(7.68)\end{eqnarray*}
on ${\cal M}^{\dg,\es}_\lmd$ for $i,j\in\ol{1,n}$ and
$r\in\mbb{N}$.

Let $\iota\in\mbb{Z}$.  Recall the Lie algebra ${\cal
L}^{\dg,\vec\ell}_{\iota}$ defined in (5.13) and the highest
weight irreducible ${\cal L}^{\dg,\vec\ell}_{\iota}$-module ${\cal
M}_{\lmd^{\dg,\es}}$ defined in (5.39) with $\tau=\dg$ and
$\lmd^{\dg,\es}(\kappa_0)=\chi$. By (5.49), we define operators
\begin{eqnarray*}E_{i,j}^{\iota,\dg}(r,z)
 &=&\sum_{l,k=0}^{\infty}[\la k\ra_r ((-1)^\es\la-k-1\ra_{\ell_j}
\E_{(-l-1)n+i-1/2,-kn-j+1/2}\\ & &-(-1)^{p(i)+p(j)}\la
-l-1\ra_{\ell_i}\E_{-kn-j+1/2,(-l-1)n+i-1/2})z^{l+k-r} \\ & &+\la
k\ra_r ((-1)^\es\la-k-1\ra_{\ell_j}\E_{ln+i-1/2,-kn-j+1/2}\\
&&-\la
l+\ell_i\ra_{\ell_i}\E_{-kn-j+1/2,ln+i-1/2})z^{-l+k-\ell_i-r-1}\hspace{9cm}\end{eqnarray*}
\begin{eqnarray*}\hspace{2cm}&&+\la-k-\ell_2-1\ra_r((-1)^\es\la
k+\ell_j\ra_{\ell_j}\E_{-(l+1)n+i-1/2,(k+1)n-j+1/2}\\ & &-\la
-l-1\ra_{\ell_i}\E_{(k+1)n-j+1/2,-(l+1)n+i-1/2})z^{l-k-\ell_j-r-1}\\
& &+\la -k-\ell_j-1\ra_r((-1)^\es\la
k+\ell_j\ra_{\ell_j}\E_{ln+i-1/2,(k+1)n-j+1/2}\\ & &-\la
l+\ell_i\ra_{\ell_i}\E_{(k+1)n-j+1/2,ln+i-1/2})z^{-l-k-\ell_i-\ell_j-r-2}\\
& &+ ((r+\ell_i)!\Im_{0,r+\ell_i} -(-1)^\es
r!\ell_i!\Im_{r,\ell_i})\dlt_{i,j}\kappa_0z^{-r-\ell_i-1}\hspace{4.1cm}(7.70)\end{eqnarray*}
on ${\cal M}_{\lmd^{\dg,\es}}$ for $i,j\in\ol{1,n}$ and
$r\in\mbb{N}$.

For convenience, we denote
$${\cal M}=\left\{\begin{array}{ll}{\cal M}^{\dg,\es}_\lmd\;\mbox{in
(4.68)}&\mbox{if}\;\iota\in\mbb{Z}+1/2,\\ {\cal
M}_{\lmd^{\dg,\es}}\;\mbox{in
(5.39)}&\mbox{if}\;\iota\in\mbb{Z}.\end{array}\right.\eqno(7.71)$$
We define linear maps
$$Y^{\iota,\pm}_{\cal M}(\cdot,z):\widehat{sp}(\vec\ell,\mbb{A})_-\rta LM({\cal
M},{\cal M}[z^{-1},z]])\eqno(7.72)$$ by
$$Y^{\iota,\pm}_{\cal M}(t^{-m-1}\ptl_t^{r+\ell_j}E_{i,j}-(-1)^{p(i)+p(j)}\es(-\ptl_t)^r
t^{-m-1}\ptl_t^{\ell_i}E_{j^\ast,i^\ast},z)=\frac{1}{m!}\frac{d^m}{dz^m}
E^{\iota,\dg}_{i,j}(r,z)^{\pm}\eqno(7.73)$$ for $i,j\in\ol{1,n}$
and $r,m\in\mbb{N}$. Now we define a linear map
$$Y^\iota_{\cal M}(\cdot,z):{\cal V}_\chi(\widehat{sp}(\vec\ell,\mbb{A}))\rta LM({\cal M},{\cal M}
[z^{-1},z]])\eqno(7.74)$$ by induction:
$$ Y^\iota_{\cal M}(\vcm,z)=\mbox{Id}_{{\cal
M}},\;\;Y(uv,z)=Y^{\iota,-}_{\cal M}(u,z)Y_{\cal M}
^\iota(v,z)+Y^\iota_{\cal M}(v,z)Y^{\iota,+}_{\cal
M}(u,z)\eqno(7.75)$$ By for $u\in
\widehat{sp}(\vec\ell,\mbb{A})_-$ and $v\in {\cal
V}_\chi(\widehat{sp}(\vec\ell,\mbb{A})).$

By Theorem 4.2, Theorem 5.2, the general theory for vertex
algebras (e.g. cf. Section 4.1 in [X2]), the  charged free
fermionic field realization and the  charged free bosonic field
realization, we obtain:\psp

{\bf Theorem 7.7}. {\it Assume $\iota\in\mbb{Z}/2$. The family
$({\cal M} ,Y^\iota_{\cal M}(\cdot,z))$ forms an irreducible
module of the vertex algebra $({\cal
V}_\chi(\widehat{sp}(\vec\ell,\mbb{A})),Y(\cdot,z),\ptl,\vcm)$.
Suppose that $\chi$ is a positive integer. If (4.69) and (5.40)
hold, then the family $({\cal M} ,Y^\iota_{\cal M}(\cdot,z))$
induces an irreducible module of the quotient simple vertex
algebra
$(V_\chi(\widehat{sp}(\vec\ell,\mbb{A})),Y(\cdot,z),\vcm,\ptl)$.
When $\lmd\in\G^\chi_{\iota,\es}$ with $\iota\in\mbb{Z}+1/2$ in
(7.49) and $\lmd^{\dg,\es}\in \in\G^\chi_{\iota,\es}$ with
$\iota\in\mbb{Z}$ in (7.53), the family $({\cal M} ,Y^\iota_{\cal
M}(\cdot,z))$ induces an irreducible module of the quotient simple
vertex algebra
$(V_{-\chi}(\widehat{sp}(\vec\ell,\mbb{A})),Y(\cdot,z),\vcm,\ptl)$.}

\vspace{0.5cm}

\noindent{\Large \bf References}

\hspace{0.5cm}

\begin{description}

\item[{[A]}] D. Adamovi\'{c}, Representations of the vertex
algebra ${\cal W}_{1+\infty}$ with a negative integer central
charge, {\it Commun. Algebra} {\bf 29} (2001), no.7, 3153-3166.

\item[{[AFMO1]}] H. Awata, M. Fukuma, Y. Matsu and S. Odake,  Quasifinite highest weight
modules over the super $W_{1+\infty}$ algebra, {\it Commun. Math.
Phys.} {\bf 170} (1995), 151-179.

\item[{[AFMO2]}] H. Awata, M. Fukuma, Y. Matsu and S. Odake, Character and determinant formulae
 of Quasifinite representation of the $W_{1+\infty}$ algebra, {\it Commun. Math. Phys.}
 {\bf 172} (1995), 377-400.

\item[{[Bo]}] R. E. Borcherds, Vertex algebras, Kac-Moody algebras, and the Monster,
{\it Proc. Natl. Acad. Sci. USA} {\bf 83} (1986), 3068-3071.

\item[{[BL1]}] C. Boyallian and J. Liberati, Classical Lie
subalgebras of the Lie algebra of matrix differential operators on
the circle, {\it J. Math. Phys.} {\bf 42} (2001), no.8, 3735-3753.

\item[{[BL2]}] C. Boyallian and J. Liberati, On modules over matrix quantum pseudo-differential
 operators ,
 {\it Lett. Math. Phys.} {\bf 60} (2002), no.1, 73-85.

\item[{[BL3]}] C. Boyallian and J. Liberati, Representations of a
symplectic type subalgebra of $W_{1+\infty}$, {\it J. Math. Phys.}
{\bf 44} (2003), no.5, 2192-2205.

\item[{[BKL]}] C. Boyallian, V. Kac and J. Liberati, Finite growth
representations of infinite Lie conformal algebras, {\it J. Math.
Phys.} {\bf 44} (2003), no.2, 754-770.

\item[{[DN1]}] C. Dong and K. Nagatomo, Classification of
irreducible modules for the vextex operator algebra $M(1)^+$, {\it
J. Algebra} {\bf 216} (1999), no. 1, 384-404.

\item[{[DN2]}] C. Dong and K. Nagatomo, Classification of
irreducible modules for the vertex operator algebra $M(1)^+$ II,
higher rank, {\it J. Algebra} {\bf 240} (2001), no. 1, 289-325.

\item[{[FKRW]}] E. Frenkel, V. Kac, A. Radul and W. Wang, ${\cal
W}_{1+\infty}$ and ${\cal W}(gl_N)$ with central charge $N$, {\it
Commun. Math. Phys.} {\bf 170} (1995), 337-357.

\item[{[FLM]}] I. Frenkel, J. Lepowsky and A. Meurman, {\it Vertex Operator
Algebras and the Monster}, Pure and Applied Mathematics {\bf 134},
Academic Press Inc., Boston, 1988.

\bibitem[HKW]{} P. Harpe, M. kervaire and C. Weber, On the Jones polynomial,
{\it L'Enseig. Math.} {\bf 32} (1986), 271-235.

\item[{[J1]}] J. C. Jantzen, Zur charakterformel gewisser darstellungen halbeinfacher grunppen
 und Lie-algebrun, {\it Math. Z.} {\bf 140} (1974), 127-149.

\item[{[J2]}] J. C. Jantzen, Kontravariante formen auf induzierten Darstellungen habeinfacher
Lie-algebren, {\it Math. Ann.} {\bf 226} (1977), no. 1, 53-65.

\item[{[J3]}] J. C. Jantzen, Moduln mit einem h\"{o}chsten gewicht, {\it Lecture Note in Math.}
 {\bf 750}, Springer, Berlin, 1979.

\item[{[K]}] V. G. Kac, {\it Vertex algebras for beginners}, University lectures series,
Vol {\bf 10}, AMS. Providence RI, 1996.

\item[{[KR1]}] V. Kac and A. Radul, Quasifinite highest weight
modules over the Lie algebra of differential operators on the
circle, {\it Commun. Math. Phys.} {\bf 157} (1993), 429-457.

\item[{[KR2]}] V. Kac and A. Radul, Representation theory of the
vertex algebra ${\cal W}_{1+\infty}$, {\it Transf. Groups} {\bf 1}
(1996), 41-70.

\item[{[KWY]}] V. Kac, W. Wang and C. Yan, Quasifinite
representations of classical Lie subalgebras of ${\cal
W}_{1+\infty}$, {\it Adv. Math.} {\bf 139} (1998), 56-140.

\item[{[L]}] W. Li, 2-cocycles on the algebra of differential
operators, {\it J. Algebra} {\bf 122} (1989), 64-80.

\item[{[M1]}] S. Ma, Conformal and Lie superalgebras motivated from
free fermionic fields, {\it J. Phys.} {\bf A36} (2003), no.6,
1759-1787.

\item[{[M2]}] S. Ma, Conformal and Lie superalgebras related to
the differential operators on the circle, {\it Ph.D. Thesis, The
Hong University of Science and Technology}, 2003.

\item[{[V]}] J. van de Leur, The $W_{1+\infty}(gl_s)$-symmetries
of the $s$-component KP hierarchy, {\it J. Math. Phys.} {\bf 37}
(1996), no.5, 2315-2337.

\item[{[X1]}] X. Xu, On spinor vertex operator algebras and their
modules, {\it J. Algebra} {\bf 191} (1997), 427-460.

\item[{[X2]}] X. Xu, {\it Introduction to Vertex Operator Superalgebras and Their Modules},
 Kluwer Academic Publishers, Dordrecht/Boston/London, 1998.

\item[{[X3]}] X. Xu, Skew-symmetric differential operators and
combinatorial identities, {\it Mh. Math.} {\bf 127} (1999),
143-258.

\item[{[X4]}] X. Xu, Simple conformal superalgebras of finite growth, {\it Algebra Colloq.}
{\bf 7}(2000), 205-240.

\item[{[X5]}] X. Xu, Equivalence of conformal superalgebras to
Hamiltonian operators, {\it Algebra Colloq.} {\bf 8} (2001),
63-92.

\end{description}

\end{document}